\documentclass[]{article}
\title{The Problem of Positive Kolmogorov-Sinai entropy for the Standard map} 
\author{Oliver Knill \thanks{Department of Mathematics,
                             University of Texas,
                             Austin, TX 78712, USA (address of 1999)}}
\date{This document replaces an announcement which circulated in 1999. 
In the present document, incorrect parts have been deleted. 
The entropy conjecture is open. The references given in the text might 
be helpful for people trying an operator theoretical or analytic approach 
to this problem.}

\newcommand{\RR}{{\bf R}} \newcommand{\ZZ}{{\bf Z}}
\newcommand{\NN}{{\bf N}} \newcommand{\CC}{{\bf C}}
\newcommand{\DD}{{\bf D}} \newcommand{\TT}{{\bf T}}

\newcommand{\Acal}{\mbox{$\cal A$}}
\newcommand{\Ecal}{\mbox{$\cal E$}}
\newcommand{\Scal}{\mbox{$\cal S$}}
\newcommand{\Xcal}{\mbox{$\cal X$}}
\newcommand{\Ycal}{\mbox{$\cal Y$}}
\newcommand{\Zcal}{\mbox{$\cal Z$}}

\newtheorem{thm}{Conjecture}[section]
\newtheorem{conj}{Conjecture}[section]
\newtheorem{lemma}[thm]{Lemma}
\newtheorem{propo}[thm]{Proposition}
\newtheorem{coro}[thm]{Corollary}
\newcommand{\cc}[2]{\put(#1,#2){\circle*{13}}}
\newcommand{\dd}[2]{\put(#1,#2){\circle*{1}}}
\newenvironment{proof}{\begin{trivlist}\item[]{\em Proof.\/\ }}
                      {\hfill
\setlength{\unitlength}{0.06mm} 
\begin{picture}(100,100)     
\put(0,0){\framebox(100,100){}} 

\dd{10}{90}\dd{40}{10}\dd{81}{40}\dd{5}{81}\dd{33}{5}\dd{77}{33}
\dd{4}{77}\dd{35}{4}\dd{80}{35}\dd{8}{80}\dd{46}{8}\dd{88}{46}
\dd{17}{88}\dd{63}{17}\dd{96}{63}\dd{23}{96}\dd{69}{23}\dd{98}{69}
\dd{24}{98}\dd{68}{24}\dd{97}{68}\dd{21}{97}\dd{62}{21}\dd{91}{62}
\dd{11}{91}\dd{43}{11}\dd{82}{43}\dd{5}{82}\dd{34}{5}\dd{78}{34}

\dd{20}{80}\dd{77}{20}\dd{16}{77}\dd{70}{16}\dd{7}{70}\dd{52}{7}
\dd{95}{52}\dd{32}{95}\dd{85}{32}\dd{23}{85}\dd{79}{23}\dd{18}{79}
\dd{73}{18}\dd{11}{73}\dd{59}{11}\dd{98}{59}\dd{35}{98}\dd{86}{35}
\dd{24}{86}\dd{79}{24}\dd{17}{79}\dd{70}{17}\dd{6}{70}\dd{49}{6}
\dd{93}{49}\dd{30}{93}\dd{83}{30}\dd{21}{83}\dd{77}{21}\dd{14}{77}

\dd{30}{70}\dd{7}{30}\dd{90}{7}\dd{64}{90}\dd{24}{64}\dd{2}{24}
\dd{82}{2}\dd{46}{82}\dd{14}{46}\dd{97}{14}\dd{75}{97}\dd{36}{75}
\dd{10}{36}\dd{95}{10}\dd{76}{95}\dd{38}{76}\dd{13}{38}\dd{100}{13}
\dd{86}{100}\dd{59}{86}\dd{23}{59}\dd{4}{23}\dd{88}{4}\dd{61}{88}
\dd{23}{61}\dd{2}{23}\dd{83}{2}\dd{48}{83}\dd{16}{48}\dd{98}{16}

\dd{40}{60}\dd{30}{40}\dd{37}{30}\dd{57}{37}\dd{69}{57}\dd{65}{69}
\dd{47}{65}\dd{32}{47}\dd{33}{32}\dd{49}{33}\dd{66}{49}\dd{68}{66}
\dd{54}{68}\dd{36}{54}\dd{31}{36}\dd{42}{31}\dd{62}{42}\dd{70}{62}
\dd{61}{70}\dd{41}{61}\dd{30}{41}\dd{36}{30}\dd{55}{36}\dd{69}{55}
\dd{66}{69}\dd{49}{66}\dd{33}{49}\dd{32}{33}\dd{48}{32}\dd{66}{48}

\dd{50}{50}\dd{50}{50}\dd{50}{50}\dd{50}{50}\dd{50}{50}\dd{50}{50}
\dd{50}{50}\dd{50}{50}\dd{50}{50}\dd{50}{50}\dd{50}{50}\dd{50}{50}
\dd{50}{50}\dd{50}{50}\dd{50}{50}\dd{50}{50}\dd{50}{50}\dd{50}{50}
\dd{50}{50}\dd{50}{50}\dd{50}{50}\dd{50}{50}\dd{50}{50}\dd{50}{50}
\dd{50}{50}\dd{50}{50}\dd{50}{50}\dd{50}{50}\dd{50}{50}\dd{50}{50}

\dd{55}{45}\dd{60}{55}\dd{54}{60}\dd{44}{54}\dd{40}{44}\dd{47}{40}
\dd{56}{47}\dd{59}{56}\dd{53}{59}\dd{43}{53}\dd{41}{43}\dd{48}{41}
\dd{58}{48}\dd{59}{58}\dd{51}{59}\dd{42}{51}\dd{41}{42}\dd{50}{41}
\dd{58}{50}\dd{58}{58}\dd{49}{58}\dd{41}{49}\dd{42}{41}\dd{51}{42}
\dd{59}{51}\dd{57}{59}\dd{48}{57}\dd{41}{48}\dd{43}{41}\dd{53}{43}

\dd{60}{40}\dd{70}{60}\dd{63}{70}\dd{43}{63}\dd{31}{43}\dd{35}{31}
\dd{53}{35}\dd{68}{53}\dd{67}{68}\dd{51}{67}\dd{34}{51}\dd{32}{34}
\dd{46}{32}\dd{64}{46}\dd{69}{64}\dd{58}{69}\dd{38}{58}\dd{30}{38}
\dd{39}{30}\dd{59}{39}\dd{70}{59}\dd{64}{70}\dd{45}{64}\dd{31}{45}
\dd{34}{31}\dd{51}{34}\dd{67}{51}\dd{68}{67}\dd{52}{68}\dd{34}{52}

\dd{65}{35}\dd{81}{65}\dd{80}{81}\dd{63}{80}\dd{33}{63}\dd{18}{33}
\dd{20}{18}\dd{37}{20}\dd{68}{37}\dd{82}{68}\dd{81}{82}\dd{63}{81}
\dd{33}{63}\dd{18}{33}\dd{18}{18}\dd{35}{18}\dd{66}{35}\dd{82}{66}
\dd{82}{82}\dd{67}{82}\dd{36}{67}\dd{19}{36}\dd{18}{19}\dd{32}{18}
\dd{63}{32}\dd{80}{63}\dd{82}{80}\dd{67}{82}\dd{37}{67}\dd{20}{37}

\dd{70}{30}\dd{93}{70}\dd{10}{93}\dd{36}{10}\dd{76}{36}\dd{98}{76}
\dd{18}{98}\dd{54}{18}\dd{86}{54}\dd{3}{86}\dd{25}{3}\dd{64}{25}
\dd{90}{64}\dd{5}{90}\dd{24}{5}\dd{62}{24}\dd{87}{62}\dd{0}{87}
\dd{14}{0}\dd{41}{14}\dd{77}{41}\dd{96}{77}\dd{12}{96}\dd{39}{12}
\dd{77}{39}\dd{98}{77}\dd{17}{98}\dd{52}{17}\dd{84}{52}\dd{2}{84}

\dd{80}{20}\dd{23}{80}\dd{84}{23}\dd{30}{84}\dd{93}{30}\dd{48}{93}
\dd{5}{48}\dd{68}{5}\dd{15}{68}\dd{77}{15}\dd{21}{77}\dd{82}{21}
\dd{27}{82}\dd{89}{27}\dd{41}{89}\dd{2}{41}\dd{65}{2}\dd{14}{65}
\dd{76}{14}\dd{21}{76}\dd{83}{21}\dd{30}{83}\dd{94}{30}\dd{51}{94}
\dd{7}{51}\dd{70}{7}\dd{17}{70}\dd{79}{17}\dd{23}{79}\dd{86}{23}

\dd{80}{10}\dd{33}{80}\dd{2}{33}\dd{72}{2}\dd{26}{72}\dd{96}{26}
\dd{63}{96}\dd{17}{63}\dd{86}{17}\dd{43}{86}\dd{7}{43}\dd{78}{7}
\dd{31}{78}\dd{1}{31}\dd{73}{1}\dd{27}{73}\dd{98}{27}\dd{67}{98}
\dd{21}{67}\dd{92}{21}\dd{53}{92}\dd{12}{53}\dd{81}{12}\dd{35}{81}
\dd{3}{35}\dd{73}{3}\dd{27}{73}\dd{97}{27}\dd{65}{97}\dd{19}{65}
\end{picture}                 
\end{trivlist}}

\begin{document}
\maketitle

\begin{abstract}
The problem of positive Kolmogorov-Sinai entropy of the Chirikov-Standard map
$T_{\lambda f}: (x,y) \mapsto (2 x-y+\lambda f(x),x)$ with $f(x)=\sin(x)$
with respect to the invariant Lebesgue measure on the two-dimensional is open.
In 1999, we believed to have a proof that the entropy can be bounded below by 
$\log(\lambda/2) - C(\lambda)$ with $C(\lambda)={\rm arcsinh}(1/\lambda)+\log(2/\sqrt{3})$
and that for $\lambda > \lambda_0=(8/(6 - 3\sqrt{3}))^{1/2}=3.1547...$, 
the entropy of $T_{\lambda \sin}$ should be positive. This approach was 
based on an idea of M. Herman using subharmonic estimates.
\end{abstract}

\section{The entropy problem of the Standard map}
The Chirikov-Taylor Standard map
$$  T_{\lambda f}:
  \left( \begin{array}{c} x         \\
                          y         \\ 
  \end{array} \right)
          \mapsto 
  \left( \begin{array}{c} 2x-y+ \lambda f(x)  \\ 
                            x  \\
          \end{array} \right)     \;   $$
with $f(x)=\sin(x)$ and real parameter $\lambda$
is a measure preserving diffeomorphisms on the two-dimensional 
torus $\TT^2=\RR^2/(2 \pi \ZZ)^2$. It is probably the most famous example 
of a symplectic twist map. This map appeared in 1960 in the context 
of electron dynamics in microtrons (see \cite{Chi87}). 
It was first numerically studied by Taylor in 1968
and Chirikov in 1969 (these independent studies are unpublished but see 
\cite{Fro70,Chi79}) and also known under the 
name "kicked rotator" and describing ground states of 
the Frenkel-Kontorova model \cite{KoFr38,AuDa83}.
It is often used to illustrate or motivate various more general mathematical 
theorems in smooth dynamical systems (i.e. \cite{Sinai94})
the calculus of variations (i.e. \cite{Mat84,MaFo94})
perturbation theory or renormalisation group techniques
(for example to understand the break up of invariant tori 
\cite{Mackay,Sti93,Sti97}).  
While it is known that for $\lambda \neq 0$, the map $T_{\lambda f}$
is non-integrable, has positive topological entropy 
and horseshoes (i.e. \cite{Fon90,Ang90,Gel99}), the question, 
whether hyperbolicity can hold on a set of positive Lebesgue measure 
stays open. Similarly, while many Lebesgue measure
preserving diffeomorphisms on the torus are known to be non-ergodic with 
positive topological entropy (i.e. \cite{Tak72,Zeh73}), it is not known 
whether positive metric entropy is dense in the $C^{\infty}$ topology. 
The issue of the positivity of the Lyapunov exponents on 
some set of positive Lebesgue measure for Hamiltonian systems 
has been addressed at various places or reviews 
\cite{Moser,Zeh73,Chi79,Wig81,Mac86,MaMe87,Hal88a,Spe89,Str89,Car91,
Dela91,Kni92,Dua94,You95,Mac95,Via98}. 
According to \cite{Pal93,Dua94}, the particular mathematical problem of
positive entropy of the Chirikov Standard map had been promoted in the
early 80'ies by Sinai. The textbook \cite{Sinai94} states on p. 144 
a conjecture (H2) that the entropy of the Chirikov Standard map 
is positive for all $\lambda>0$ and that the entropy grows to infinity 
for $\lambda \to \infty$. \\

Although numerical experiments show a very clear lower bound 
$\log(|\lambda|/2)$ for the entropy of the Chirikov Standard map 
(e.g. \cite{Fro70,Chi79,PaSe}), it was not known even for one
single value of $\lambda$, whether the metric 
entropy with respect to the invariant Lebesgue 
measure can be positive. The second part of the (H2) conjecture of Sinai. 

\begin{conj}
\label{I}
\fbox{ \parbox{12cm}{
The Kolmogorov-Sinai entropy $\mu(T_{\lambda \sin})$ of the 
Chirikov-Standard map $T_{\lambda \sin}$ with respect to the invariant 
Lebesgue measure is bounded below by $\log(|\lambda|/2)-C(\lambda)$ 
where $C(\lambda)={\rm arcsinh}(1/\lambda)+\log(2/\sqrt{3})$. 
}}
\end{conj}

Numerical experiments suggest that one could probably even hope to get
rid of the $C(\lambda)$-term in the case of the Standard map. 
Similar statements should hold for more general Standard maps and
expect explicit lower bounds 
$\log(|m \lambda|/2)-C(m \lambda)$ if $\lambda \sin(x)$ is replaced by 
$E x+\lambda \sin(m x)$, where the integer $E \in \ZZ$ 
allows additionally to tune the homotopy type of the map. 
This property of positive metric entropy should be stable 
in the real analytic category and dense in the $C^0(\TT)$ topology: 

\begin{conj}
\label{II}
\fbox{ \parbox{12cm}{
There is in $C^0(\TT)$ a $C^0$-dense set of real-analytic Standard maps 
$f: \TT^1 \mapsto \RR$ for which $\mu(T_f)>0$. Let $f$ be in this
dense set. In every Banach space of realanalytic maps, in which $f$ is,
there is an open neighborhood of realanalytic maps $g$, for which $\mu(T_g)>0$ also. 
}} 
\end{conj}

The conjectured sensitive dependence on initial 
conditions could be persistent in the realanalytic category and should be 
obtained by realanalytic, $C^0(\TT)$-small perturbations of integrable maps. 
The stability of positive metric entropy with respect to realanalytic 
perturbations of the map would make the result physically relevant. 
Other results in Hamiltonian dynamics that have both this stability and 
which deal with orbits forming a set of positive probability is the theory 
of Anosov maps (see i.e. \cite{Lan85,Yoc95,KH}) 
(an example is $(x,y) \mapsto (4x+\lambda \sin(x)-y,x)$ for small 
$\lambda$) or KAM perturbation theory (see i.e. \cite{Moser,Her86,KaOr93})
(which applies for example near the integrable SBKP map (\cite{Sur89,Bo+93})
$(x,y) \mapsto (2x-y+4 {\rm arg}(1+\lambda \exp(-i x)),x)$).
It is not known whether there are open sets of
realanalytic Hamiltonian maps and flows for which there is quasiperiodic 
motion on a set of positive Lebesgue measure and simultaneously 
a conjugation to a Markov chain on a different set of positive 
Lebesgue measure. 
\footnote{There is a possibility to achieve such 
mixed behavior on the union of two closed symplectic manifolds, where 
the map is Anosov on one component and integrable on the other. 
For smooth, not realanalytic examples in the 
literature that are unstable under perturbations of the map.} \\

One could ask similar questions to higher-dimensional symplectic maps 
$(x,y) \mapsto (E x-y+ \lambda f(x),x)$ on $\TT^{2d}$, where $E$
is a constant symmetric matrix in $GL(d,\ZZ)$ and $f$ is a vector
valued, real-analytic function on the torus. 
Examples are Froeschle maps, where $E x=2 x$, $f_i(x)=\lambda_i \sin(x_i)+ \sin(\sum_{j} x_j)$
or classes of nearest neighbor coupled map lattices \cite{KaGr87}
on $\TT^{2d}$, 
where $E x=2 x$, $f_i(x)= \lambda \sin(x_i) 
  + \epsilon (\sin(x_{i+1}-x_i) + \sin(x_i-x_{i-1}))$ and $x_{i+d}=x_i$. 
All averaged Lyapunov exponents should be nonzero for large enough $\lambda$. \\

As for any diffeomorphism $T$ on a compact manifold $M$
leaving invariant a smooth measure, the entropy is by the Pesin 
formula \cite{Man81} equal to the sum of the positive 
integrated Lyapunov exponents. In the two-dimensional case, 
the entropy is $\lim_{n \rightarrow \infty} n^{-1}$ 
$\int_M \log||dT^n(x,y)|| \; dx dy$. 
Since an estimate of the norm of the Jacobean of the 
Chirikov-Standard map $T_{\lambda \sin}$ 
shows that the integrated Lyapunov exponent is 
bounded above by 
$\int_{\TT^2} \log||dT_{\lambda \sin}(x,y)|| \; dx dy 
     <  \log(\lambda/2)+c(\lambda)$ with $c(\lambda)=O(1/\lambda)$, 
Conjecture~\ref{I} provides for large $\lambda$ a rather accurate estimate 
for the actual values of the entropy of the Chirikov Standard map and 
shows that the measure of the Pesin region, the set on 
which the Lyapunov exponents 
is positive goes to $1$, when $|\lambda| \to \infty$.  \\

By Pesin theory \cite{Pes77,Mane,KS,Pollicott,KH}, positive metric entropy 
(and in higher dimensions the non-vanishing of all Lyapunov exponents
on a set of positive measure) would imply the existence of invariant sets 
of positive Lebesgue measure on which $T_{\lambda}$ is ergodic and on 
which some iterate of $T_{\lambda}$ is mixing and actually
measure-theoretically conjugated to a Bernoulli shift. Other consequences
are the density of periodic orbits in the Pesin region and shadowing 
properties (see \cite{Kat80,KH}). \\

An obstacle for proving positive entropy estimates for the Standard map 
is Donskaya's observation that elliptic islands can
exist for arbitrary large $\lambda$ \cite{Chi79,Che91,Lichtenberg,Dua94} and 
which make it hard to find invariant cone bundles \cite{Woj85,Woj86} 
in the tangent bundle. In \cite{Dua94} it has been shown 
using renormalisation of coordinates near homoclinic tangencies 
\cite{MoRo97} that for a residual set of large
parameters $\lambda$ one has Standard maps $T_{\lambda \sin}$ which 
are nonergodic. \\

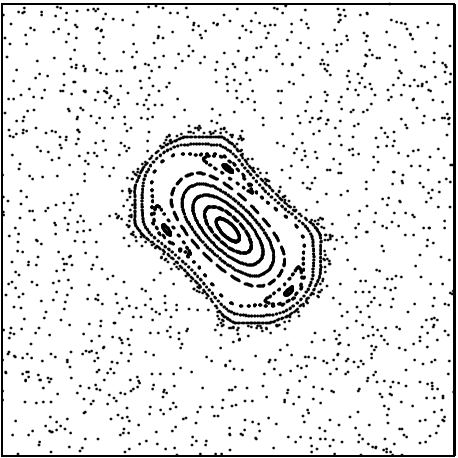
\hbox{
\vbox{ \parbox{5.2cm}{ \vspace{-0cm}
\newcommand{\ee}[2]{\put(#1,#2){\circle*{5}}}
\setlength{\unitlength}{0.06mm} 
                 
}}
\hspace{-8cm}
\vbox{ \parbox{7.5cm}{
Chirikov had already expressed concerns that 
the measure of the stable elliptic component might be close to 
$1$ (\cite{Chi79} p. 333). It has been often asked whether an
area-preserving monotone twist map has a dense set of elliptic islands
in general (see e.g. \cite{Hal88a,Pal91}). It is known that a Baire generic
symplectic non-Anosov $C^1$-diffeomorphism has a dense set of elliptic 
periodic orbits \cite{New75}. 
Some research towards avoiding elliptic islands by
smooth ergodic perturbations of the Chirikov Standard map has been done in
\cite{Sin96}. It is possible to get positive entropy by a smooth 
$C^1$ perturbation of the map \cite{Her94}. Whether this is possible 
with $C^{\infty}$ perturbations is not known \cite{Her98}. 
}}}
\vspace{0.5cm} 

The question whether a dense set in ${\rm Diff}_{\mu}^{\infty}(M)$ of measure 
preserving diffeomorphisms on a manifold $M$ has positive metric entropy has been 
asked in \cite{Her94a}.  \\

It has been conjectured that there exists a set of parameters $\lambda$
with full density at $\infty$ for which the Chirikov Standard map has
no elliptic islands \cite{Car91} (see also \cite{Dua94}). \\

Our proof-attempts of positive metric entropy depended on spectral and 
complex analytic techniques as well as the determinant theory of 
finite von Neumann algebras avoiding the ergodicity question. 
We hoped to obtain many non-ergodic realanalytic Standard maps with positive metric
entropy. The picture above for example shows plots of some orbits of 
the Chirikov Standard map in the case $\lambda=3.4$, where a linearly stable
fixed point $(1/2,1/2)$ coexists with a region with 
positive metric entropy.  \\

There are Standard maps $T_{\lambda \sin}$ arbitrarily close to this map $T_{3.4 \sin}$
for which the Birkhoff normal form at $(1/2,1/2)$ is such that 
the fixed point is surrounded by invariant KAM curves of positive Lebesgue 
measure and for which we still have ergodic components of 
positive Lebesgue measure. Together 
with \cite{Dua94}, we know that there exist for large $\lambda$
a open dense set of parameters
for the Chirikov Standard map which lead to nonergodicity and 
positive Kolmogorov-Sinai entropy.  \\
 
Consider the matrix-valued map 
$$ (x,y) \in \TT^2 \mapsto  A_{E,\lambda f}(x,y)  
               =  \left(
                    \begin{array}{cc}
                          E+\lambda f(x) & -1  \\
                            1              &  0  \\
                    \end{array} 
                  \right)        \;    .   $$ 
Together with a Lebesgue measure-preserving dynamical system 
$T: \TT^2 \to \TT^2$, it defines a cocycle 
$(x,y,n) \mapsto A^n_{E,\lambda f,T}(x,y) 
   = A_{E,\lambda f}(T^{n-1}(x,y)) \circ \cdots 
     \circ A_{E,\lambda f}(x,y)$ 
for which the integrated Lyapunov exponent 
$$ \mu(A_{E,\lambda f,T}) =
     \lim_{n \to \infty} \frac{1}{n} \int_{\TT^2} 
      \log||A^n_{E,\lambda f,T}(x,y)|| \; dx \; dy \;        $$
is defined. For general information on Lyapunov exponents, see
for example \cite{Rue79a,Led84}. \\

Pesin's formula links the entropy of $T_f$ with the Lyapunov exponent 
of $dT_f$ and so with the Lyapunov exponent of the 
cocycle because $A_{E=2,df,T_f}$ 
agrees with the Jacobean $dT_f$ of the map $T_f$.  \\

\section{A nonselfadjoint spectral problem}
Consider the random operator 
\begin{equation}
\label{L(w)}
 (L_{w}(x,y)) u_n = u_{n+1} + u_{n-1}
             + \lambda (w^{-1} \exp(i x_n) + w \exp(-ix_n)) u_n \; , 
\end{equation}
on $l^2(\ZZ,\CC)$, 
where $(x_n,y_n)=T^n(x,y)$ and $w \in \CC$ is a complex parameter. 
For fixed $(x,y) \in \TT^2$, this is a bounded operator on 
$l^2(\ZZ,\CC)$. The term 'random' is used because $(x,y) \mapsto L(x,y)$
is an operator valued random variable. Such operators are in a 
von Neumann algebra with finite trace 
${\rm tr}(K) = \int_{\TT^2} [K(x,y)]_{00} \; dx dy$ (see e.g. 
\cite{MuNe43,Dixmier,Connes}). \\

Similar as in the case, when the complex parameter is the energy,  
the Lyapunov exponent $\mu(w)$ of the transfer cocycle $A_E(w)$ 
is a subharmonic function in $w$. As we will see below, 
the Thouless formula  
$\mu(w)=\log(\det(L_w)) = {\rm tr}(\log|L_w|)$ 
is still true for the in general 
nonselfadjoint operator $L_w$. \\
The problem is to compute $\det(L_w)$. The entropy of the Standard
map is $\log \det L_{w=1}$. \\

The average value of the Lyapunov exponent $\mu(w)$ on the circle
$|w|=1$ can be estimated with the subharmonicity 
argument of Herman \cite{Her83} or with a Jensen formula 
as used by Sorets-Spencer in \cite{SoSp91}. 
In higher dimensions, an adaption \cite{GoSo92}
of \cite{SoSp91} to higher-dimensional cocycles would apply. \\

A theorem of Lax \cite{Lax71} on continuous 
measure-preserving transformations allows to approximate the group 
$\Xcal$ of measure preserving homeomorphisms by finite groups 
$\Ycal_k$ of measure preserving transformations. 
If one looks at the value of the Lyapunov
exponent of the matrix-valued map $A_{E,\lambda \sin,T}$ on $\Ycal_k$, (where
evaluating the Lyapunov exponent is a reliable finite dimensional integration),
one can observe numerically that for every transformation $T$, 
one can find a transformation $\hat{T}$ such that 
$  \mu(A_{E,\lambda \sin,T}) 
 + \mu(A_{E,\lambda \sin,\hat{T}}) \geq 2 \log(\lambda/2)$. 
This relation can be understood in terms of determinants. \\

By looking for a proof of this relation, we were motivated by Aubry duality 
which can be defined if $T$ is in such a finite group $\Ycal_k$. 
Aubry duality is an involutive transformation $L \mapsto \hat{L}$,
which preserves the density of states such that $\log(\det(L_{\lambda}))
= \log(\det( (\lambda/2) \hat{L}_{(4/\lambda)})) 
= \log(\lambda/2) + \log \det(\hat{L}_{(4/\lambda)})$. 
In the Mathieu case $T(x,y)=(x+y,y)$, where $L_{\lambda}
= \hat{L}_{\lambda}$, this leads with 
$\log(\det(\hat{L}_{4/\lambda})) \geq 0$ 
to the estimate $\log \det(L_{\lambda}) \geq \log(\lambda/2)$
which is expected to be true also if $T$ is the Standard map.
Here, we will show that there exists for each transformation 
$T$ and corresponding operator $L=L_{T,\lambda}$
a different transformation $S$ and operator $\hat{L}=L_{S,\lambda}$ such that 
$\log(\det(L \hat{L})) \geq 2 \log(\lambda/2)$.  \\

How is the operator $\hat{L}$ obtained? The group
$R_{\alpha}: (x,y) \mapsto (x+\alpha,y)$ of translations 
on the torus acts by conjugation 
on the whole group $\Ycal$ of measure preserving transformations. The orbits 
of this action are $T_{\alpha} = R_{\alpha} T R_{-\alpha}$. 
The operator $\hat{L}$ will be obtained by changing $T$ to one of the 
conjugates $T_{\alpha}$. The value of the 
Lyapunov exponent $w=\exp(i\alpha) \mapsto \mu(A_{E,\lambda f, T_{\alpha}})$ 
extends to a subharmonic function from the circle $|w|=1$ to the entire complex
plane $\CC$.
If one considers the complex parameter $w$ as a 
'spectral parameter', the Riesz measure $dk$ used to represent
the Lyapunov exponent plays the role of the 
density of states. It is a measure in the complex plane which has in our case
in general positive area in the complex plane. The measure $dk$ has its
support in an annular neighborhood of the unit circle (and an atom at $0$). \\

We use the artificial spectral parameter $w$
because with respect to $z=r e^{ix}$, there is no analyticity of 
$[z^n A^n(z,y)]_{ij}$.
The motivation is trying to generalize results 
for maps which extend analytically to polydiscs like 
for $T(x,y)=(x+\alpha,y)$, where both Aubry duality \cite{AvSi82}
as well as Herman's subharmonicity \cite{Her83} work directly. 
While the map $(t e^{ix}, e^{iy}) 
\mapsto (t e^{i (2x-y + \lambda \sin(x))},e^{ix})$ (or any similar 
complexification attempt) is real-analytic in each of the variables $t,x,y$,
it is definitely not analytic with respect to $z = t e^{ix}$ and 
functions which needed to be subharmonic are not. 
Our first attack on the Standard map using 
plurisubharmonicity was done in the spring 1988, now more than 
eleven years ago. The approach was the observation that the Standard map 
family written in Hamiltonian form 
$(x,y) \mapsto (x+y+\lambda \sin(x),y+\lambda \sin(x))$
is induced on invariant tori of the single analytic map 
$U:(z,w,u,v) \mapsto (z w e^{z-u}, w e^{z-u}, u v e^{u-z}, v e^{u-z})$
on $\CC^4$, and where nonanalyticity manifests itself that
the $U$-invariant two-dimensional tori 
$S_{\lambda}=\{ (z,w,u,v) \; | \; |z|=|u|=\lambda/2, \; |w|=|v|=1 \; 
z = \overline{u}, w=\overline{v}\}$
on which $U$ induces 
$T_{\lambda}$, are not distinguished boundaries of polydiscs. 
Later attempts were to Vlasov-Toda deform the cocycle with the aim 
to estimate the Lyapunov exponent on tori of isospectrally 
deformed operators (on which the Lyapunov exponent is constant) 
\cite{Kni93a,Kni93b} or to fix the dynamical system and to vary the cocycle 
\cite{Kni92}. \\

Because $w \mapsto g_n(w,x,y)=w^n [A_E^n(w,x,y)]_{11}$ is analytic,
the $r=|w|$ dependence of the local Lyapunov exponent can be related
with the average angular dependence of the argument: 
in polar coordinates $w=r e^{i \phi}$, the Cauchy-Riemann differential 
equations 
$$ \frac{d}{d \phi} \frac{1}{r}
                            \int_{\TT^2} \arg(g_n(w,x,y)) \; dx dy 
   = \frac{d}{dr} \int_{\TT^2} \log|g_n(w,x,y)| \; dx dy \;  $$
are valid in any connected region of the $w$- resolvent set. 
In the complement, in the support of the $w$-spectrum, this formula can 
not be used and an integrated version, a Jensen formula replaces it. 
Besides the contribution of the radial change of the argument, 
there is an additional nonnegative subharmonic contribution to 
the Lyapunov exponent. 

\section{More heuristics}
Let $A_{E,\lambda \cos}(z,y) = 
\left(
                    \begin{array}{cc}
                          E+\frac{\lambda}{2} (z+z^{-1}) & -1 \\
                            1                            &  0 \\
                    \end{array} \right) $. 

In situations like the Mathieu case $T(x,y)=(x+y,y)$, where 
$T$ extends to an analytic map on $\DD \times \TT$ with respect to 
$z=r \exp(i x)$, Herman's subharmonicity \cite{Her83}
gives $\mu(A_T)-\log(\lambda/2) \geq 0$. For a general map $T \in \Ycal$ which
does no more commute with $R_{\alpha}: (x,y) \mapsto (x+\alpha,y)$ on $\TT^2$,
the function $w=\exp(i \alpha) \mapsto \mu(A_{E,\lambda \sin,T_{\alpha}})$
with $T_{\alpha} = R_{\alpha} T R_{-\alpha}$ is no more constant. 
In average, we have $\int_{\TT} \mu(A_{E,\lambda \sin,T_{\alpha}}) \; d \alpha
\geq \log(\lambda/2)$. Furthermore $\mu(A_{E,\lambda \sin,T_{\alpha}})$
is bounded above by $\log(\lambda/2) + c(\lambda)$ with 
$c(\lambda)=O(1/\lambda)$. This means that for
large $\lambda$, we have $\mu(A_{E,\lambda \sin,T_{\alpha}}) \geq \log(\lambda/2)$ for a large 
set of $\alpha$'s. We used such a fact in a similar way in
\cite{Kni92}. \\

However, we don't know the specific value of the upper-continuous function 
$\alpha \mapsto \mu(A_{E,\lambda \sin,T_{\alpha}})$ at
the point $\alpha=0$, which we are interested in. It could be zero a priori. 
Numerical computations however show that the minimum of 
$\alpha \mapsto \mu(A_{E,\lambda \sin,T_{\alpha}})$ is not farther
away from the mean value than the maximum.  
There is a heuristic explanation which uses the
Jensen formula in a sector and which will assume $T \in \Xcal$.
This Jensen formula (see Section~\ref{a_jensen_formula})
is essentially an integrated Cauchy-Riemann differential equation
and sharpens the subharmonicity tool. 
We also use the Lax approximation theorem 
in Section~\ref{toral_homeomorphisms_and_lax_approximation}
to explain the now following heuristics.  \\

The Lyapunov exponent of the matrix-valued map
$(x,y) \in \TT^2 \mapsto B(e^{i x})$ with
$$ B(z) = z A_{E,\lambda \cos}(z)  =   z \left(
                    \begin{array}{cc}
                          E+\frac{\lambda}{2} (z+z^{-1}) & -1 \\
                            1                            &  0 \\
                    \end{array} \right) $$
is the same as the Lyapunov exponent of the cocycle
$(x,y) \mapsto A_{E, \lambda \cos}(x)$ because $|z|=1$. \\

Define for $(x,y) \in \TT^{2}$ and $n \in \NN$ the complex function
$$ w \mapsto g_{n,T}(w^{-1} \exp(ix),y)
      = [ B^n_{E, \lambda\cos, T}(w^{-1} \exp(ix),y) ]_{11}  $$
which is analytic in the complex plane. \\

Denote by ${\rm Arg}_{[z_1,z_2]}(g_{n,T})$ the argument change of 
the analytic function $g_{n,T}$ on the line from $z_1$ to $z_2$. 
Unlike ${\rm arg}(g_{n,T})$, the argument change
${\rm Arg}_{[z_1,z_2]}(g_{n,T})$ is a uniquely defined number. \\

Let $T_k$ be a cyclic Lax cube exchange transformation approximating
a homeomorphism $T \in \Xcal$. The map $T_k$ is extends from $\TT^2$
to $\DD \times \TT$ by setting
$T_k(r \exp(i x),y) =(r \exp(i x_1),y_1)$ with $T_k(x,y) = (x_1,y_1)$. \\

In each sector, we could also take the $z$ variable instead of the 
$w$ variable because $z \mapsto A^n(z,y)$ is analytic 
in each sector $\Acal_{s,t,j} = \{ w \in \CC \; | \; |w| \in [s,t], 
{\rm arg}(w) \in [2 \pi j/k,2 \pi (j+1)/k) \; \}$. \\

\vspace{5mm}
\hbox{
\vbox{ \parbox{5.2cm}{ \vspace{-3cm}
 \setlength{\unitlength}{0.06mm} 
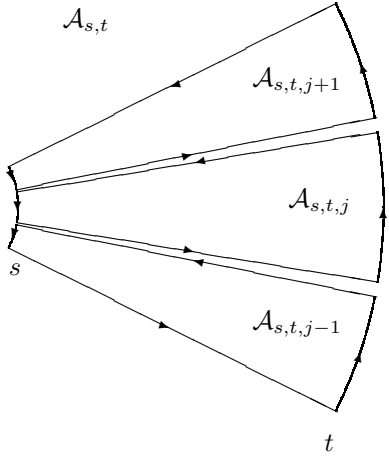
 \begin{picture}(1000,1000)      
 \qbezier(996,34)(1005,89)(1008,144)
 \qbezier(1008,144)(1012,200)(1008,256)
 \qbezier(1008,256)(1005,311)(996,366)
 \qbezier(197,167)(199,178)(200,189)
 \qbezier(200,189)(200,200)(200,211)
 \qbezier(200,211)(199,222)(197,233)
 \qbezier(990,398)(985,424)(978,451)
 \qbezier(978,451)(972,477)(964,502)
 \qbezier(964,502)(956,528)(946,553)
 \qbezier(946,553)(937,578)(926,603)
 \qbezier(926,603)(915,628)(903,652)
 \qbezier(196,239)(195,244)(194,250)
 \qbezier(194,250)(192,255)(191,260)
 \qbezier(191,260)(189,265)(187,270)
 \qbezier(187,270)(186,275)(183,280)
 \qbezier(183,280)(181,285)(179,289)
 \qbezier(990,2)(985,-24)(978,-51)
 \qbezier(978,-51)(972,-77)(964,-102)
 \qbezier(964,-102)(956,-128)(946,-153)
 \qbezier(946,-153)(937,-178)(926,-203)
 \qbezier(926,-203)(915,-228)(903,-252)
 \qbezier(196,161)(195,156)(194,150)
 \qbezier(194,150)(192,145)(191,140)
 \qbezier(191,140)(189,135)(187,130)
 \qbezier(187,130)(186,125)(183,120)
 \qbezier(183,120)(181,115)(179,111)
 \put(197,167){\line(6,-1){800}}
 \put(197,233){\line(6,1){800}}
 \put(196,239){\line(5,1){800}}
 \put(179,289){\line(2,1){730}}
 \put(196,161){\line(5,-1){800}}
 \put(179,111){\line(2,-1){730}}
 \put(592,101){\vector(4,-1){0.2}}
 \put(592,299){\vector(-4,-1){0.2}}
 \put(588,82){\vector(-4,1){0.2}}
 \put(588,318){\vector(4,1){0.2}}
 \put(537,-68){\vector(3,-1){0.2}}
 \put(537,468){\vector(-3,-1){0.2}}
 \put(1010,200){\vector(0,1){0.2}}
 \put(200,190){\vector(0,-1){0.2}}
 \put(187,130){\vector(-1,-3){0.2}}
 \put(190,260){\vector(1,-3){0.2}}
 \put(960,520){\vector(-1,3){0.2}}
 \put(960,-120){\vector(1,3){0.2}}
 \put(300,600){${\cal{A}}_{s,t}$}
 \put(800,200){${\cal{A}}_{s,t,j}$}
 \put(720,470){${\cal{A}}_{s,t,j+1}$}
 \put(720,-70){${\cal{A}}_{s,t,j-1}$}
 \put(180,50){$s$}
 \put(880,-340){$t$}
 \end{picture}                  
}}
\hspace{-8cm}
\vbox{ \parbox{7.5cm}{
Let $\Gamma^+(\alpha)$ be path connecting  
$s e^{i \alpha}$ and $t e^{i \alpha}$ on a radial line segment,
$\Gamma^-$ the reversed path. Let
$\gamma^+_s(\alpha,\beta) = \{ w, |w|=s, {\rm arg}(w) \in [\alpha,\beta] \}$
be a circular arc contained in the circle $\gamma_s=\{|w|=s\}$ and 
$\gamma^-_s(\alpha,\beta)$ the reversed path. 
The difference between the argument change of the function
$g_{n,T}(w)=[B^n_{E,\lambda \sin,T}(w)]_{11}$ along the 
circles $\gamma^-_s$ and $\gamma^+_t$ enclosing the annulus 
$\Acal_{s,t} = \{ |w| \in [s,t] \}$
is the sum of the $k$ angular argument changes 
${\rm Arg}_{\gamma^+_t(2\pi j/k,2\pi (j+1)/k)}(g_{n,T})
-{\rm Arg}_{\gamma^-_s(2\pi j/k,2\pi (j+1)/k)}(g_{n,T})$
plus the sum of the $k$ radial argument changes 
${\rm Arg}_{\Gamma^+(2\pi j/k)}(g_{n,T})
 - {\rm Arg}_{\Gamma^-(2\pi (j+1)/k)}(g_{n,T})$ 
with $j=1, \dots, k$. 
}}} 
\vspace{1.5cm}

The total sum of the radial plus angular argument changes is nonnegative
because it is the sum of the indices of the roots of $z \mapsto g_{n,T}$
inside the annulus $\Acal_{s,t}$, all of which are nonnegative.  \\

The sum of the $k$ radial argument changes is bounded above by a 
$O(1/\lambda)$-term,
because a positive value contributes positively to the 
Lyapunov exponent estimate and the Lyapunov exponent is bounded above
by $\log(\lambda/2) + c(\lambda)$ with $c(\lambda)=O(1/\lambda)$. 
There is no reason, why the sum of the $k$
radial argument changes should differ much 
from its negative. Numerical experiments strongly confirmed that 
they don't. If this sum of the radial argument changes 
\begin{equation}
\label{argumentchange}
 \sum_{j=1}^k  {\rm Arg}_{\Gamma^+(2\pi j/k)}(g_{n,T})
 - {\rm Arg}_{\Gamma^-(2\pi (j+1)/k)}(g_{n,T})  \;  
\end{equation}
(which is averaged over $\TT^2$ and bounded uniformly in $k$) 
would vanish, one would have the same estimates as in the subharmonic
case and
$$  n^{-1} \int_{\TT^2} \log|g_{n,T}(x,y,w)| \; dx dy \geq
    n^{-1} \int_{\TT^2} \log|g_{n,T}(x,y,0)| \; dx dy = \log(\lambda/2)  \;  $$
for all $|w|=1$. This would imply $\mu(A_{E,\lambda \cos,T}) 
\geq \log(\lambda/2)$. \\

While the sum of the radial argument changes~(\ref{argumentchange})
does not vanish for a general 
$T \in \Ycal$, numerical data makes ~(\ref{argumentchange})
appear to be small. 
Indeed, we will see that for every transformation $T \in \Ycal$, there is 
an other transformation $S_{\alpha}=R_{\alpha} T R_{-\alpha}$
such that the Lyapunov exponent of $A_T$ plus the Lyapunov exponent 
of $A_S$ beats $2 \log(\lambda/2)$. This will imply that the 
fluctuations of 
$$ \alpha \mapsto \mu_n(\alpha) 
   = n^{-1} \int_{\TT^2} \log|[A^n(e^{i (\alpha+x)},y)]_{11}| \; dx dy $$
around its mean 
$\int_{\TT} \mu_n(\alpha) \; d\alpha \geq \log(\lambda/2)$ 
can not get too large. \\

We have not been able to prove such a statement using complex analytic
methods only. 

\section{The spectrum of periodic difference operators} 
This  section revies some facts about not necessarily selfadjoint 
periodic Jacobi matrices and higher order difference operators. 
Unlike in the selfadjoint case, when the spectrum is confined to the real 
axes, the topology of the spectrum can then be more interesting. 
All we need to know for our purposes is that the spectrum 
is contained in a nowhere dense set. \\

While nonselfadjoint differential operators have been studied quite a bit
by Russian mathematicians, nonselfadjoint higher order difference operators 
have appeared less frequently in the literature. But they do occur: 
these operators appear naturally when
studying the stability of stationary solutions in 
coupled nonharmonic oscillators \cite{Way99,JoWa99}.
For nonselfadjoint Schr\"odinger differential operators,
see \cite{GeWe95} and references. \\

{\bf Jacobi matrices}. 
Three vectors $\underline{a},\underline{b},\underline{c} \in \CC^p$ 
define periodic sequences
$a_n=a_{n+p},b_n=b_{n+p},c_n=c_{n+p} \in l^{\infty}(\ZZ,\CC)$
which can be used to build a $p$-periodic Jacobi matrix 
$$ (Lu)_n=a_n u_{n+1} + c_{n-1} u_{n-1} + b_n u_n   $$
which is a bounded operator on on $l^2(\ZZ,\CC)$. 
Denote by $L_{per}$ the same operator
$L$ but acting on a different Hilbert space, namely the finite dimensional 
space of $p$-periodic sequences equipped with the norm
$(\sum_{i=1}^p |u_i|^2)^{1/2}$. For $w \in \CC$, (a parameter which has
nothing to do with the complex parameter $w$ in the operator (\ref{L(w)}) ), 
this operator is 
given by the matrix
$$ L_{per}(w)=\left( \matrix{
b_1  & a_1  &  0       &\cdot      & 0               & w^{-1} c_p      \cr
c_1  & b_2  &a_2       &\cdot      &\cdot            & 0               \cr
 0   & c_2  &\cdot     &\cdot      &\cdot            &\cdot            \cr
\cdot&\cdot &\cdot     &\cdot      & a_{p-2}         & 0               \cr
 0   &\cdot &\cdot     &c_{p-2}    & b_{p-1}         &a_{p-1}          \cr
w a_p  & 0   &\cdot &0        & c_{p-1}          &b_{p}                \cr
                                             } \right) \;  $$
and satisfies $L_{per}=L_{per}(1)$. 
Denote by $\Delta(z,w)=\det(L_{per}(w)-z)$
the characteristic polynomial of $L_{per}(w)$ and by $\Delta(z)=\Delta(z,1)$
the characteristic polynomial of $L_{per}$. Comparing coefficients in $w$ and 
using $\Delta(z,1)=\Delta(z)$, one gets
$$ \Delta(z,w)=\Delta(z) - a w - c w^{-1} - b   \;  $$
with $a= \prod_{j=1}^p a_j$  and $c= \prod_{j=1}^p c_j$ and $b=-a-c$.
Denote by $\tau$ the shift $u_n \mapsto u_{n+1}$ on $l^2(\ZZ,\CC)$. \\

By the theorem of Cayley-Hamilton, we have
$\Delta(L_{per})=0$.  We deduce from this
that the operator $\Delta(L,\tau^p)$
on $l^2(\ZZ,\CC)$ defined by the functional calculus (plug-in the 
operator $L$ for $z$ and the operator 
$\tau^p$ for $w$ in $\Delta(z,w)$) satisfies
$$ 0 = \Delta(L_{per}) 
     = \Delta(L,\tau^p)=\Delta(L) - a \tau^p - c \tau^{-p}-b  \; . $$
This implies that the operator
$K:=\Delta(L)=a \tau^p + c (\tau^p)^* +b$ on $l^2(\ZZ,\CC)$ 
is a Laplacian with space independent entries.
Conjugated by a Fourier transform $U: l^2(\ZZ,\CC) \rightarrow L^2(\TT,\CC)$,
it is unitarily equivalent to the multiplication operator
$$ (U K U^*)u(\theta) = 
    (a \; e^{ip \theta} + b + c \; e^{-ip \theta}) u(\theta) 
       =: \lambda(\theta) u(\theta) $$
on $L^2(\ZZ,\CC)$ which has as the spectrum on an ellipse $\Ecal$
possibly being degenerated to an interval. 
Since $\Delta(L)=K$, we know by the spectral theorem that
$\sigma(\Delta(L))=\Delta (\sigma(K))$.
The spectrum of $L$ consists of $p$ closed real curves. They are the zeros of
$z \mapsto \Delta(z)-w$, when $w$ is runs through the ellipse $\Ecal$.
The density of states of $L$ is the push-forward measure 
$(\Delta^{-1} \circ \lambda)^* d\theta $
of the Lebesgue measure $d\theta$ on the circle $\TT=\RR/(2 \pi \ZZ)$
under the multi-valued 
map $\phi=\Delta^{-1} \circ \lambda: \TT \mapsto \CC$. 
(See \cite{GKT} for a derivation in the case $a_n=c_n=1$).  \\

Remarks. \\
1)  Since the spectrum of Jacobi matrices is simple,
the curves $z_i(\theta)$ do not intersect transversely. However,
by deforming along a one-parameter family of Jacobi matrices, curves 
could merge or separate. The discriminant 
$ \{ L \; | \; (\Delta_L)'(z)=0
                        {\rm \; for \; some} \; z \in \sigma(L) \}  $
is the set of Jacobi matrices, where the topology of the spectrum 
can change under small parameter variations
by merging or separation of spectral arcs.  \\

2) The Bloch variety of $L$ is the set
$ BV(L)
  =\{ (z,w) \in \CC^* \times \CC^* \; | \; \Delta(z) = a w +b + c w^{-1} \} $,
where $\CC^*$ is the Riemann sphere $\CC \cup \{\infty\}$.
It is the set of all $(z,w)$ such that there exists a 
nontrivial solution  of $Lu=zu$ satisfying
$(\tau^p u)_n=u_{n+p}=w u_n$. 
The spectrum of $L$ is the intersection of $BV(L)$ with $\{|w|=1\}$. \\

3) A theorem of Polya (see \cite{AigZie})
implies that the projection of the spectrum of a periodic (not necessarily
selfadjoint) Jacobi matrix onto any real line in the complex plane 
has Lebesgue measure $\leq 4$. \\

4) In the case when $a_j,b_j,c_j$ are real and $a_n=c_{n}$ and 
$a=\prod_j a_j=1$, the Jacobi matrix $L$ is real and selfadjoint and
$K$ is the free Laplacian $\tau+\tau^*-2$ with spectrum $[-2,2]$ so that
$\sigma(L)=\{z \; | \;  \Delta(z) \in [-2,2] \} $.
The spectrum is then a union of bands on the real line
which can intersect only at the boundaries. \\

5) In the case when $a=\prod_j a_j=0,c=\prod_j b_j=0$, 
the spectrum is $\Delta^{-1}(0)$, the set of
roots of the characteristic polynomial of $L_{per}$ confirming that 
in this case, $L$ and $L_{per}$ have the same spectrum. ($L$ is then a
countable direct sum of operators $L_{per}$.) \\

{\bf Higher order difference operators}. 
For a $p$-periodic operator $L = \sum_{i=-r}^r a_i \tau^i$, 
the spectrum of $L_{per}$ is more complex. The spectral bands
curves can cross and for selfadjoint matrices, the bands 
can overlap. \\

We consider the case $p=k N$, where the operator can be rewritten
as a Jacobi operator 
$L = a \tau + c\tau^* +b$ on $l^2(\ZZ,\CC^N)$,
where $a_j,b_j,c_j$ are $N \times N$ matrices satisfying 
$a_{j+k} = a_j, b_{j+k}=b_j, c_{j+k} = c_j$. 

\begin{lemma}[Spectrum of periodic matrix-valued Jacobi matrices]
\label{matrix-valued-Jacobi}
The spectrum of the Jacobi operator $L$ with $N \times N$-matrix 
valued coefficients is contained in the union of at most $p=k N$ 
closed real curves. 
\end{lemma} 
\begin{proof}
Unlike in the case $N=1$, the Bloch variety 
$$ \{ (z,w) = \Delta(z,w)= {\rm det}(L_{per}(w)-z) = 0  \} \subset \CC^2 $$
can in general no more given explicitly in the form $f(w)=g(z)$. 
The spectrum is still 
the intersection of this complex variety with $\{ |w|=1 \}$, 
the parameter $|w|=1$ labeling Bloch waves. \\

The intersection of the Bloch variety with $|w|=1$
is a set in the complex plane which is contained in a 
union of at most $p=k N$ closed curves.  \\

Let ${\rm Det}$ denote the determinant function on the algebra
of complex $N \times N$ matrices. We write
$\prod_{j=1}^k a_j = a_k \dots a_1$,
also in the noncommutative case. With $a = {\rm Det} \prod_{j=1}^k a_j,
c = {\rm Det} \prod_{j=1}^k c_j$, we have 
$$ \Delta(z,w) = \Delta(z) + a w^N + c w^{-N}  + f(z,w) \; , $$
where the mixed term $f(z,w)$ is a polynomial in $w,w^{-1},z$ which is 
of degree $<N$ in $w,w^{-1}$ and of degree $<p$ in $z$.  \\
\end{proof} 

Examples. \\
1) We are especially interested in 
higher order difference operators which arethe product
$L = L^{(1)} L^{(2)}$ of two Jacobi matrices 
$L^{(i)} = \tau_i+\tau_i^* +V^{(i)}$. \\

2) If $a_n,b_n,c_n$ are diagonal, the operator $L$ is the direct
sum of $N$ one-dimensional Jacobi matrices. There are therefore operators with 
$p = k N$ spectral curves. 

\section{Thouless formula for difference operators}
This section reviews the Thouless formula for operators 
on the strip, that is Jacobi matrices with matrix-valued entries.
We recall the proof in \cite{KoSi88} because the class of 
operators differs a tiny bit from the
operators considered in \cite{KoSi88}: the off diagonal matrices are 
the identity to which some nilpotent part is added. An other change is
that we do not assume selfadjointness of the operator. 
The proof of the Thouless formula however is the same.  \\

The Thouless formula $\mu(A_E) = \int \log|E-E'| \; dk(E')$ relates
the density of states $dk$ of a random Jacobi operator 
$\tau + \tau^* + b$ with the Lyapunov exponent $\mu(A_E)$ of the 
transfer cocycle (see \cite{Cycon,Carmona}). This generalizes to 
operators $\tau+\tau^* +b$ on $l^2(\ZZ,\CC^N)$, where $\tau u_n=u_{n+1}$
and $b(x)$ is a selfadjoint $N \times N$ matrix, if $\mu(A_E)$ is 
the average of the $N$ largest Lyapunov exponent of the $2N \times 2N$-cocycle
\cite{KoSi88}. 
The Thouless formula extends also to the case of operators 
$a \tau + (a \tau)^* +b$, where $a(x,y)$ are unitary \cite{KniII}. 

The next lemma provides a generalization to operators of the type
$L=L^{(1)} L^{(2)}$ with Jacobi matrices 
$$(L^{(k)}(x^{(k)},y^{(k)}) u)_n=
   u_{n+1} + u_{n-1} + V^{(k)}_n(x^{(k)},y^{(k)}) u_n \; , $$
where the later are not necessarily selfadjoint operators on $l^2(\ZZ,\CC)$, 
\begin{equation}
\label{operatoronstrip}
   L(x^{(1)},y^{(1)},x^{(2)},y^{(2)})  = \left( \begin{array}{ccccccc} 
  \cdot &   \cdot   & \cdot    & 1         & 0        &   0     &   0     \\
     1  &   c_{n-2} & b_{n-1}  & a_{n-1}   &  1       &   0     &   0     \\
     0  &       1   &  c_{n-1} & b_{n}     & a_{n}    &  1      &   0     \\
     0  &       0   &   1      & c_{n}     & b_{n+1}  & a_{n+1} &   1     \\
     0  &       0   &   0      &   1       & c_{n+1}  & b_{n+2} & \cdot   \\
     0  &       0   &   0      &   0       &   1      & \cdot   & \cdot   \\
                 \end{array}
          \right)  \; ,
\end{equation} 
where $a_n=a_n({\bf x})
 = V_n^{(1)}(x^{(1)},y^{(1)})+V_{n+1}^{(2)}(x^{(2)},y^{(2)}),
 b_n= b_n({\bf x}) 
  = V_n^{(1)}(x^{(1)},y^{(1)}) V_n^{(2)}(x^{(2)},y^{(2)}) + 2, 
c_n= c_n({\bf x}) 
= V_n^{(1)}(x^{(1)},y^{(1)}) + V_{n-1}^{(2)}(x^{(2)},y^{(2)})$. 
The potentials $V_n^{(j)}(x^{(j)},y^{(j)})$ are obtained from $T_j \in \Ycal$ and
functions $V^{(j)}: \TT^2 \mapsto \RR$ by 
$V_n^{(j)}(x^{(j)},y^{(j)}) = V^{(j)}(T_j^n(x^{(j)},y^{(j)}))$. \\

Each difference operator operator 
$L({\bf x})$ on $l^2(\ZZ,\CC)$ can be written 
as a Jacobi operator 
\begin{equation}
\label{operatoronstrip2} 
\tilde{L}= \tilde{a} \tilde{\tau} + (\tilde{c} \tilde{\tau})^* + \tilde{b}, 
(\tilde{L} u)_{\tilde{n}} = \tilde{a}_{\tilde{n}} u_{\tilde{n}+1} 
             + \tilde{c}_{\tilde{n}-2} u_{\tilde{n}-1} 
             + \tilde{b}_{\tilde{n}} u_{\tilde{n}}  
\end{equation}
on $l^2(\ZZ,\CC^2)$ with $\tilde{\tau} u_{\tilde{n}} = u_{\tilde{n}+1}, 
\tilde{\tau} \tilde{c}({\bf x}) 
  = \tilde{c}(T_1^2(x^{(1)},y^{(1)}),T_2^2(x^{(2)},y^{(2)})) \tilde{\tau}$ 
and with matrix-valued entries 
$$    \tilde{b}_{\tilde{n}} = \left(
            \begin{array}{cc}
                        b_{n-1}  & a_{n-1}   \\
                        c_{n-1}  & b_{n} \\
            \end{array}
            \right) \; , 
      \tilde{a}_{\tilde{n}} = \left(
            \begin{array}{cc}
                        1        & 0      \\
                        a_{n}    & 1      \\
            \end{array}
            \right) \; , 
      \tilde{c}_{\tilde{n}} = \left(
            \begin{array}{cc}
                        1        & c_{n}  \\
                        0        & 1    \\
            \end{array}
            \right) \;  \; , $$
where $n=2 \tilde{n}$. The transfer cocycle 
\begin{equation}
\label{transfercocycle}
  A_E = \left( \begin{array}{cc} 
                 (E-\tilde{b}_{\tilde{n}}) 
           \tilde{a}_{\tilde{n}-2}^{-1}  &  - \tilde{c}_{\tilde{n}-2} \\
           \tilde{a}_{\tilde{n}-2}^{-1}  &    0                       \\ 
                 \end{array}
          \right) 
      = \left( \begin{array}{cccc}
     E-b_{n-1} + a_{n-1} a_{n-2} & -a_{n-1} & -1  & -c_{n-2} \\
      a_{n-2} (b_n-E) - c_{n-1}       & E-b_n     &  0  &  -1      \\
     1                                &  0        &  0  &  0       \\
   -a_{n-2}                           &  1        &  0  &  0       \\
           \end{array}
    \right)
\end{equation} 
satisfies 
$$ A_E \left( \begin{array}{c} \tilde{a}_{n-2} u_{n}     \\
                                               u_{n-1}   \\
              \end{array}
       \right)
     = 
       \left( \begin{array}{c} \tilde{a}_{n} u_{n+1} \\
                                             u_{n}       \\
              \end{array}
       \right) 
 $$
if $(\tilde{L}({\bf x}) u)_n = 
  E u_n$ for $u \in l^2(\ZZ,\CC^2)$.
The density of states of $L$ defined as 
the functional $f \mapsto 
\int_{\TT^4} [f(L({\bf x})]_{00} \; 
     d {\bf x}$ on $C(\CC)$ and the 
density of states of $\tilde{L}$ defined as the functional 
$f \mapsto \int_{\TT^2} 
  {\rm Trace}([f(\tilde{L}({\bf x})]_{00}) \; {\bf x}$
both exist and they are the same probability measure in $\CC$, 
if ${\rm Trace}$ is the normalized trace on $2 \times 2$ matrices satisfying 
${\rm Trace}(1)=1$. \\

Remarks.  \\
1) For $L = L^{(1)} L^{(2)}$, where $L^{(i)}$ are Jacobi operators, there
are other transfer matrices like the product of the individual 
$2 \times 2$ transfer matrices of the operators $L^{(i)}$. These are different
matrix-valued functions in $E$ and the Thouless formula would look different
in this case. \\

2) Note that $dk$ is normalized to be a probability measure. In some
of the literature like for example in \cite{KoSi88}, 
the usual ${\rm Trace}$ for $N \times N$ matrices is
taken which implies that the total mass of the density of states
$dk$ for operators $L$ on the strip becomes $N$. 

\begin{lemma}[Thouless formula for operators on the strip]
\label{thoulessstrip}
Let $dk$ be the density of states of the 
operator~(\ref{operatoronstrip2}). Then
$$ \int_{\CC} \log|E-E'| \; dk(E') = \mu(\wedge^2 A_E)/2  \; , $$
where $\mu(\wedge^2 A_E)/2$ is the arithmetic mean of the 
two largest Lyapunov exponents of the $4 \times 4$ transfer 
cocycle $A_E$ of the operator~(\ref{operatoronstrip2}). 
\end{lemma} 
\begin{proof} 
Small changes are needed to the known proofs \cite{Cycon,Carmona,KoSi88} 
for the same statements for selfadjoint operators $L$, where 
$a_n=c_n=1, n \in \ZZ$.  Write ${\bf x}=
(x^{(1)},y^{(1)},x^{(2)},y^{(2)})$ for a point on $\TT^4$. \\

(i) It is enough to prove the statement for periodic, (nonergodic)
transformations $T \times S$. \\
Proof. Every $T \times S \in \Ycal \times \Ycal$ can be approximated 
by periodic transformations $T_n \times S_n$ using Rohlin's 
lemma \cite{Hal56} (which implies that there exist measurable sets $Y_n$ 
whose Lebesgue measure on $\TT^4$ goes to $1$ for $n \in \infty$ so that 
$(T \times S)^k({\bf x}) = (T_n \times S_n)^k({\bf x})$ 
for $k=0, \dots n-1$).  \\

The Thouless formula is an identity for subharmonic functions 
$f_{T \times S} = g_{T \times S}$, parameterized by $T \times S \in \Ycal \times \Ycal$. 
Let $dk(f)$ denote the Riesz measure of a subharmonic function $f$ it satisfies
$dk(f) = \Delta f$ in the sense of distributions, where $\Delta$ is the Laplacian. 
It follows from the Avron-Simon lemma and the dominated convergence theorem that 
if $T_n \times S_n \to T \times S$ in the uniform topology for measure 
preserving transformations, 
then $dk(f_{T_n \times S_n}) \to dk(f_{T \times S})$ 
weakly. Also $dk(g_{T_n \times S_n}) \to dk(g_{T \times S})$ weakly
holds because $dk(g_T)$ is a weak limit of measures
$dk(g_T^{(k)})=\int_{\TT^4} \sum_j \delta(E_j^{(k)}({\bf x})) \; d{\bf x}$, 
where $E_j^{(k)}({\bf x})$
are the roots of $E \mapsto [A^k_{T,E}({\bf x})]_{11}$ and $T \mapsto 
dk(g_T^{(k)})$ is continuous from $\Ycal$ with the uniform topology 
to the space of measures on $\CC$ with the weak topology. 

By assumption, we have $dk(f_{T_n \times S_n}) = dk(g_{T_n \times S_n})$. 
Therefore $dk(f_{T \times S}) = dk(g_{T \times S})$ which is
$\Delta f_{T \times S}= \Delta g_{T \times S}$. From this, 
$f_{T \times S}=g_{T \times S}$ follows by Weyl's lemma \cite{Ransford}. \\ 
 
(ii) Next, it is sufficient to assume that $E$ is outside 
the spectrum of a periodic operator 
$L$ (which is now contained in a finite union of one-dimensional
curves (Lemma~(\ref{matrix-valued-Jacobi}))) 
as well as that $E$ is outside the spectrum of $L^{(0)}$, 
where $L^{(0)}$ is the Jacobi matrix on $l(\ZZ^+,\CC)$ with zero 
Dirichlet boundary condition at $n=0$. 
The spectrum of $L^{(0)}$ is a finite set of points. \\
Proof.  This is the Craig-Simon subharmonicity argument \cite{CrSi83}. 
Taking balls $B_r(E)$ around a point $E$ in the spectrum, the result holds
for the over $B_r(E)$ averaged Lyapunov exponent. Passing to the limit
$r \to 0$ is possible because of subharmonicity.  \\

(iii) If $\lambda_j(E) = \prod_{k=1}^n (E-\lambda_{jk}), j=1, \dots, N$
is an eigenvalue of the $N \times N$ matrix
$[A_E^n]_{11}$, (where the $2N \times 2N$ matrix $A_E$ is written as a 
$2 \times 2$ matrix of $N \times N$ matrices),
then each root $\lambda_{jk}$ of $E \mapsto \lambda_j(E)$ 
gives rise to an eigenvector 
$u^{(jk)}$ of a truncation $L^{(n)}$ of $L$, where $a_k,b_k,c_k$ 
are put to zero for $k<0$ as well as for $k>n$. The measures 
$dk_n= (2n)^{-1} \sum_{j,k} \delta(\lambda_{jk})$ converge for $n \to \infty$
weakly to the density of states $dk$ which is supported
on finitely many curves.  \\
This is the Thouless-Avron-Simon
lemma \cite{Cycon} which holds in general for random finite-difference 
operators which do not need to be selfadjoint. \\

(iv) The fact that 
$$ \int_{\CC} \log|E-E'| \; dk_n({\bf x},E')
= \sum_{j,k} \log | E-\lambda_{jk}({\bf x}) | 
  = \log| {\rm Det}([A_E^n]_{11}({\bf x}))| $$ 
converges to $\mu(\wedge^2 A_E({\bf x}))/2$ 
for almost all 
${\bf x} \in \TT^4$ in the limit $n \to \infty$ is
seen in the same way as in the Appendix of \cite{KoSi88}: 
If $m^{\pm}({\bf x})
  =a({\bf x}) u^{\pm}
    ( (T \times S) {\bf x}) 
    u^{\pm}({\bf x})^{-1}$ 
are the matrix-valued Titchmarsh-Weyl functions, where 
$L u^{\pm}_n({\bf x}) 
= Eu^{\pm}_n({\bf x})$ and
$\{u^{\pm}_n({\bf x})\}_{n \in \pm \ZZ^+} 
\in l^2(\pm \ZZ^+,M(N,\CC))$, then 
$$(\mp 1/2)  \int_{\TT^4}  
   \log|{\rm Det}(m^{\pm}({\bf x}))| \; 
   d{\bf x}$$
is the arithmetic mean $\mu(\wedge^2 A_E({\bf x})/2$ 
of the two largest Lyapunov exponents of $A_E$. 
Furthermore, the maximal Lyapunov exponent averaged over $\TT^4$ is
$$ - \int_{\TT^4}  \log||m^{+}({\bf x})|| \; 
                        d{\bf x}
   = \int_{\TT^4}  \log||m^{-}({\bf x})|| \;
                        d{\bf x} \; . $$
Because $E$ is not an eigenvalue of $L({\bf x})$, 
there are two $N \times N$-matrices $U,V$ such that the functions
$m^{\pm},m^{(0)} = a u^{(0)}(T \times S) (u^{(0)})^{-1}$ satisfy 
$m^{(0)} = m^-U + m^+ V$. Here $u^{(0)}({\bf x})$
is an other solution of $L({\bf x})u=Eu$ 
satisfying a Dirichlet boundary condition at $0$. The matrix $U$ 
is invertible because else there would exists $a \neq 0$ with 
$Ua=0$ and $m^{(0)} a=m^+ V a$ would be in 
$l^2(\ZZ^+,\CC^2)$ and contradict the fact that 
$E$ was also chosen away from the spectrum of $L^{(0)}$. 
(In the selfadjoint case, this would already have been taken care 
of by choosing ${\rm Im}(E) \neq 0$.) Writing
$$ m_n^{(0)} = (1+m_n^+VU^{-1} (m_n^-)^{-1}) m_n^- U $$
and using $m^{(0)}_n({\bf x}) 
  = [A^n_E]_{11}({\bf x})$,
one concludes
\begin{eqnarray*}
 \log \det (L) &=& \int_{\CC} \log|E-E'| \; dk(E')                           \\
        &=& \lim_{n \to \infty} \int_{\TT^4} 
              \int_{\CC} \log|E-E'| \; 
             dk_n({\bf x},E') \;                                             \\
        &=& \lim_{n \to \infty} \int_{\TT^4} 
          \log {\rm Det}([A^n_E({\bf x})]_{11}) \; 
              d{\bf x}                                                       \\
        &=& \frac{1}{2} \lim_{n \to \infty} \int_{\TT^4}  
          \log|{\rm Det}(m^{(0)}({\bf x})_n)|  \; 
             d{\bf x}                                                        \\
        &=& \frac{1}{2} \int_{\TT^4} 
          \log|{\rm Det}(m^-({\bf x}))| \; d{\bf x}  = \mu(\wedge^2 A_E)/2 \; .
\end{eqnarray*}
\end{proof} 

\section{An upper bound for the Lyapunov exponent}
In this section we determine the constant $C(\lambda)$ in 
Conjecture~(\ref{I}) and formulate a Proposition which outlines, how one 
can improve the constant with more effort. The reason to stick to the constant
$C(\lambda)$ is that it allows explicit expressions. We also show in this
section that $C(\lambda)$ can be replaced by $c(\lambda) = O(\lambda^{-1})$. \\

Because $||A^n|| \leq ||A||^n$ for any matrix $A$, the number 
$$ m(E,\lambda) := \int \log|| A_{E,\lambda \cos} (x,y) || \; dx dy $$
is an upper bound for 
the Lyapunov exponent of the Standard map, rsp. the Lyapunov exponent
of the cocycle $A_E(w)$ for $|w|=1$. Define 
$C(E,\lambda) = m(E,\lambda)- \log(\lambda/2)$
and $C(\lambda)= C(0,\lambda)$. 

\begin{lemma}
\label{normcalculation} 
$C(E,\lambda) = \log(2/\sqrt{3}) + O(1/\lambda)$ for large $\lambda$ and
$C(0,\lambda)>0$ for $\lambda > \lambda_0=2 \sqrt{2/(6 - 3\sqrt{3})}=3.1547...$.
\end{lemma} 
\begin{proof} 
With the norm $|||A|||=\max_{i=1,2} |A e_i|$, where $e_i$ are the basis
vectors, there is an explicit expression for an upper bound of 
the integral $m(E,\lambda)$. Following
\cite{Jit94}, define
$$ M(E,\lambda) = |E + i + \sqrt{(E+i)^2 - \lambda^2} |/\sqrt{3} \;  $$
where the square-root takes the solution with positive imaginary part. \\

For matrices 
$A=\left( \begin{array}{cc} c & 1 \\ -1 & 0 \end{array} \right)$ one has
$|||A||| =| \sqrt{1+c^2} | 
\geq (\sqrt{3}/2) ||A||$. Furthermore, with the analytic function 
$f(z) = z^2 (1+(E-\lambda (z+z^{-1})/2)^2) = \lambda/2 (z-a_-) (z-a_+)$ 
where $|a_-|<1$, one has using the Jensen formula
\begin{eqnarray*}
   \int_{\TT} \log||| A_{E,\lambda \cos} (x) ||| \; dx       
    &=& \frac{1}{2} \int_{\TT} \log|1+(E-\lambda \cos(x))^2| \; dx     \\
    &=& \frac{1}{2} \int_{\{ |z|=|\exp(ix)|=1 \}} \log|f(z)| \; dx     \\
    &=& \frac{1}{2}(\log(f(0)) - \log(|a_-|))                          \\   
    &=& \log(\frac{\sqrt{3}}{{2}} M(E,\lambda)) \; . 
\end{eqnarray*}
Therefore 
$C(E,\lambda) = \log(M(E,\lambda))-\log(\lambda/2)
= \log(2/\sqrt{3}) 
+ \log(|(E+i)/\lambda + \sqrt{(E+i)^2/\lambda^2 - 1}|)$
and $C(0,\lambda)=C(\lambda) = 
\log(1/\lambda + \sqrt{1/\lambda^2 + 1}) +\log(2/\sqrt{3})=
{\rm arcsinh}(1/\lambda)+\log(2/\sqrt{3})$.
For $\lambda = 3.1547.. = 2 \sqrt{2/(6 - 3\sqrt{3})}$,
one has $C(\lambda)=0$. \\
The function $\lambda \mapsto C(\lambda)$, having a positive derivative,
is strictly monotone. 
\end{proof}

\begin{lemma}
We have 
$$   \int_{\TT^2} \log||dT_{E,\lambda \cos}(x,y)|| \; dx dy 
      \leq c(E,\lambda) := \int_{\TT} \log \sqrt{2+(E+\lambda \cos(x))} \; dx
      = O(\lambda^{-1})  \; . $$
\end{lemma} 
\begin{proof}
We use Hadamard's determinant theorem, which tells that for 
a finite matrix $L$, one has $\log \det(L) \leq \sum_{j=1}^n ||a_j||_{2}$, 
where $a_j$ are the rows of the matrix $L=(a_1, \cdots , a_n)$. 
(The relation follows from the fact that if $A=(a_1, \dots, a_n)$, 
then ${\rm det}(A)/\prod_{j=1}^n ||a_j||_{2} = {\rm det}(B)$, where 
$B=(a_1/||a_1||, \dots, a_n/||a_n||)$. The determinant theorem is 
equivalent to ${\rm det}(B)\leq 1$, a relation which can easily be 
seen geometrically (i.e. \cite{Bel43,Dix84}). ) \\

In the random case, where the $a_j$ are vector valued random variables
on $\TT^2$ with Lebesgue measure, 
and $L^{(n)}$ is a finite dimensional $n \times n$ approximation of $L$
and where $n^{-1} \log \det(L^{(n)})$ converges to to the 
Fuglede-Kadison determinant $\log \det(L)$, we obtain from Hadamard's
determinant inequality and Birkhoff's 
ergodic theorem
$$ \log \det(L) \leq \int_{\TT^2} \log ||a(x,y)|| \; dx dy   \; . $$
In the case of random Jacobi operators 
$L = \tau + \tau^* + b(x) = \tau + \tau^* + \lambda \cos(x) + E$, 
where $a(x,y) = ( \dots,0,0,1,b(x),1,0, \dots )$, this gives
\begin{eqnarray*}
 \log \det(L) \leq \int_{\TT} 
   \log \sqrt{2+ (E+\lambda \cos(x))^2} \; dx \;                 
    &=& 
   \log(\lambda/2) 
     + \int_{\TT} \log \sqrt{8/\lambda^2 
       + 4 (E/\lambda + \cos(x))^2} \; dx \;, \\
    &=& \log(\lambda/2) + O(\lambda^{-1}) \; . 
\end{eqnarray*}
\end{proof} 

Improved smaller constants $C_n(\lambda)<C_{n-1}(\lambda)$ with 
$C_1(\lambda) = C(\lambda)$ could be obtained by estimating 
$$ n^{-1} \int_{\TT^2} \log ||A^n(x,y)|| \; dx dy $$
from above, where $n$ is small. But
already for $n=3$, explicit bounds become more difficult due to the
transcendental nature of the integrals.
For $n=2$, one would have to estimate the integral 
$$ \int_{\TT^2} \log|| \left( \begin{array}{cc} 
                                  E - \lambda \cos(y) & -1 \\  
                                            1         &  0 \\
                              \end{array} \right)
                       \left( \begin{array}{cc}
                                  E - \lambda \cos(x) & -1 \\
                                            1         &  0 \\
                              \end{array} \right) || \; dx dy   $$
from above. With a numerical evaluation of this integral,
we get an upper bound for the Lyapunov exponent.

The values of $C_n(\lambda)$ could be estimated computer assisted. We expect
that like this, one should reach $\lambda_n$ close to $2$ for large $n$. 
Maybe one could even get $\lambda_n \to 0$ for $n \to \infty$, if also
the lower bound of the determinant
$\int_{|w|=1} \log {\rm det} L_T K_{S,w}$ could be improved further with the 
Jensen formula. 

\section{Uniform hyperbolicity: a result of Ruelle}
In this section we look at a class of uniformly hyperbolic cocycles
and apply a result of Ruelle \cite{Rue79}
to conclude that the lower bound on the Lyapunov exponent
depends in our case smoothly on parameters. Furthermore, the
result implies that the neighborhood around such a cocycle,
for which one has positive Lyapunov exponents, becomes large
if the Lyapunov exponent is large. \\

Consider a $(d+1) \times (d+1)$- matrix-valued map 
$$ (x,y) \mapsto A(x,y) = \left( \begin{array}{cc} \lambda/2+a(x,y) & b(x,y)  \\ 
                                                     c(x,y) & d(x,y)  \\
                 \end{array} \right) \; , $$
where $a(x,y) \in \CC$, $d(x,y)$ is a $d \times d$ matrix and
$b(x,y)$ and $c(x,y)^T$ are vectors in $\CC^d$. All entries are
assumed to be bounded and measurable functions. 
Together with $T \in \Ycal$, it defines a 
cocycle $(x,y) \mapsto A^n(x,y)$. We look at cocycles in a neighborhood
of $\left( \begin{array}{cc} \lambda/2 & 0   \\
                                     0 & 0   \\
                 \end{array} \right)$, which has one Lyapunov exponent 
$\log(\lambda/2)$ and $(d-1)$ Lyapunov exponents $-\infty$. \\

The maximal Lyapunov exponent in such a uniformly hyperbolic situation 
(which generalizes the logarithm of the 
Perron-Frobenius eigenvalue of a positive matrix) depends 
smoothly on parameters \cite{Rue79}. We formulate this in a special 
case: 

\begin{lemma}[Ruelle]
The maximal Lyapunov exponent $A \mapsto \mu(A)$ of is Fr\'echet 
differentiable and realanalytic in a neighborhood of 
$A_0=\left( \begin{array}{cc}  \lambda/2 & 0 \\ 0 & 0 \\ \end{array} \right) \in 
L^{\infty}(\TT^2,M(d+1,\CC))$, where 
$M(d+1,\CC)$ is the vector space of all complex $(d+1) \times (d+1)$
matrices. If $w \mapsto A_w$ be a complex parameterization of 
the cocycle such that for all $\{ |w| \leq r \}$, the inequalities
$$ |a_w(x,y)|,||b_w(x,y)||,||c_w(x,y)||,||d_w(x,y)|| \leq \epsilon \; $$
hold. 
Then, for fixed $\lambda$ and small enough $|\epsilon|$, or for fixed 
$\epsilon>0$ and large enough $\lambda$, the Lyapunov exponent 
$w \in \{ |w| \leq r \}  \mapsto \mu(A(w))$
is a harmonic function in $w$.
\end{lemma}
\begin{proof} 
If $|a(x,y)|, ||b(x,y)||, ||c(x,y)||, ||d(x,y)|| \leq \epsilon$ and
$\epsilon$ is small enough, the cocycle maps the constant cone
bundle
$$ (x,y) \mapsto C(x,y) = C 
  = \{ c \cdot (1,\eta) \; | \;  ||\eta|| \leq 1, c \in \CC \;  \} $$
strictly into itself. \\

(One could write the complex matrix cocycle as a real
matrix cocycle on $\RR^{2(d+1)}$ for which the cone 
$$ C = \{ c \cdot (1,\eta) \; | \;  ||\eta|| \leq 1, \eta \in \RR^{2d+1},
          c \in \RR \;  \} $$
in $\RR^{2d+2}$ is strictly invariant. 
This is equivalent to the fact that one can
find a basis and $\delta>0$ such that 
$|[A(x)]_{ij}| \geq \delta$.
By (\cite{Rue79} 4.7), one can deal directly with the complex case. ) \\

\cite{Rue79} is written in the case of continuous
cone bundles and homeomorphisms, but remark 4.9 in \cite{Rue79} 
extends this to the measure theoretic case. 
\end{proof} 

From Proposition~4.8 in \cite{Rue79} we obtain especially: 

\begin{coro}
\label{analytic}
For small $r$, the Lyapunov exponent
$(w,z) \mapsto \mu(G(w,z))$
is pluriharmonic (harmonic both in $w$ and $z$) 
on $\{ (w,z) \; | \; 0 \leq |w|<r, 0<|z|<r \;  \}$. 
\end{coro}

\begin{coro}
\label{smallchangeinlargeball}
Given $R>0$ and $\epsilon>0$. Let $B_R(A)$ be the ball of 
radius $R$ around $A$. For large enough $\lambda$, the ball
$B_R(G(w,z))$ is contained in the open subset of uniformly hyperbolic 
cocycles and $|\mu(B) - \mu(A)| \leq \epsilon$ for all $B \in B_R(A)$. 
\end{coro} 
\begin{proof} 
The Fr\'echet derivatives of the function $A \mapsto \mu(A)$
can be computed (see formula (3.6) in \cite{Rue79}) and
satisfy $D_A^n \mu(A) \leq C^n/\lambda^n$, where $C$ is a constant
independent of $\lambda$. \\

Actually, one only needs the first derivative: Bloch's theorem
which says that if a function $f$ is analytic in a disc $B_R$ of radius $R$,
then $f(B_R)$ contains a disc of radius $|D f(a)| R l$, where $l$ is a
constant not depending on anything (see \cite{Andersson} p. 39). 
\end{proof} 

\section{Determinants of random difference operators}
In this section we look at determinants
for classes of random operators. We will see that for a class of 
random operators, which do not need to be invertible, 
the determinant is nonnegative and 
satisfies ${\rm det}(L K) = {\rm det}(L) {\rm det}(K)$. 
This product formula of Fuglede and Kadison \cite{FuKa52} for operators
in a type $II_1$ von Neumann algebra usually assumes the operators $L$ and 
$K$ to be invertible. In the case of random, finite order difference operators, 
we can approximate the operators by periodic operators, which have 
the spectrum on a union of one-dimensional curves. By changing the
energy we can approximate the operators by invertible operators for which the 
formula holds. The formula holds then for the random difference operators also. \\
 
A measure $dk$ in the complex plane is called the density of states of a
(not necessarily selfadjoint) random bounded operator 
$(x,y) \in \TT^2  \mapsto L(x,y) 
       = \sum_{j=-p}^p a^{(j)}(x,y) \tau^j$ on $l^2(\ZZ,\CC^N)$
with $N \times N$ matrices $a^{(j)}$ and $\tau u_n = u_{n+1}$, 
if for every polynomial $f$, one has
$$ \int_{\TT^2} {\rm Trace}([f(L(x,y))]_{00})  \; dx \; dy
 = \int_{\CC} f(E) \; dk(E) \;  $$
where ${\rm Trace}$ is the usual trace for $N \times N$ matrices 
normalized in such a way that ${\rm Trace}(1)=1$. \\

For operators $L$ with a density of states $dk$ one has a trace 
$$ 
  {\rm tr}(L)= \int E \; dk(E) = \int {\rm Trace} [L(x,y)]_{00} \; dx \; dy \; . 
$$

If the Thouless formula holds and the Lyapunov exponents of the transfer
cocycles are all finite, then the measure $dk$ has finite logarithmic 
energy and the potential $E \mapsto {\rm tr}( \log|L-E| )
= \int_{\CC} \log|E-E'| \; dk(E')$ is defined and 
finite for all $E$. \\

In general, if $L$ has a density of states $dk$, 
a determinant ${\rm det}(L)$ of $L$ can be defined by 
$$ \log {\rm det}(L) = {\rm tr}( \log(|L|) )  \; . $$
If the value $\int_{\CC} \log|E| \; dk(E)$ of $dk$ at the point $E=0$ 
should be $-\infty$, one would put ${\rm det}|L|=0$. \\

(The extension of the determinant to noninvertible operators, which we
look at, is what Fuglede and Kadison call
an 'analytic extension of the determinant' to singular operators.  
This is in contrast to the 'algebraic extension' for which 
the determinant is put always to zero if the operator is noninvertible).  \\

Given $T,S \in \Ycal$ and $a^{(j)},b^{(j)} \in L^{\infty}(\TT^2,dx dy)$. 
Two not necessarily selfadjoint
elements $L=\sum_j a^{(j)}(x,y) \tau^j, K=\sum_j b^{(j)}(x,y) \sigma^j$ in 
the crossed products define random operators
$$  (L(x,y) u)_n = \sum_{j=-p}^p a^{(j)}_n(x,y) u_{n+j},  \; \; 
    (K(x,y) u)_n = \sum_{j=-p}^p b^{(j)}_n(x,y) u_{n+j} $$
with $a^{(j)}_n(x,y) = a^{(j)}(T^n(x,y)), b^{(j)}_n(x,y) = b^{(j)}(S^n(x,y))$. 

\begin{lemma}
\label{determinant} 
The operators $L,K,LK$ have all density of states. The Thouless formula
holds for $L,K$ and $L K$ and
${\rm det}(L K) = {\rm det}(L) {\rm det}(K)$. 
\end{lemma} 
\begin{proof} 
(i) The operators $L,K$ can be embedded in a single von Neumann algebra
without changing the determinant. \\
Proof. $L$ is in the crossed product $\Zcal_T$ of $L^{\infty}(\TT^2,dx dy)$ 
with the $\ZZ$-action $\alpha_T: f \mapsto f(T)$ and
$K$ is in the crossed product $\Zcal_S$ of $L^{\infty}(\TT^2,dx dy)$
constructed with the $\ZZ$-action $\alpha_S: f \mapsto f(S)$.
The operator $L K$ is in the crossed product $\Zcal_{T \times S}$ of
$L^{\infty}(\TT^2 \times \TT^2,dx dy \times dx dy)$, constructed with
the $\ZZ$-action $\alpha_{T \times S}$, where 
$T \times S (x_1,y_1,x_2,y_2) = (T(x_1,y_1),S(x_2,y_2))$
on $\TT^2 \times \TT^2$. 
One can embed the two crossed products $\Zcal_T$ and $\Zcal_S$ in
$\Zcal_{T \times S}$.
The density of states of $L K$ exists in the same way as it
exists for $L$ and $K$. Moreover, the density of states of $L$ (rsp. $K$) in 
$\Zcal_{T \times S}$ is the same as the density of states of $L$ in 
$\Zcal_T$ (rsp. in $\Zcal_S$). \\

(ii) For the verification of the Thouless formula, see Section 15). \\

(iii) We can assume without loss of generality that $L$ and 
$K$ are invertible. \\
Proof. By ergodic decomposition (i.e. \cite{DGS}), one can assume $T,S$
to be ergodic. \\
Approximate $T,S\in \Ycal$ by periodic transformations $T_k$
using Rholin's lemma \cite{Hal56}: the measure of 
$\{ (x,y) \in \TT^2, \; T(x,y) 
\neq T_k(x,y) \}$ is then smaller than $1/k$
and $T_k^k=Id$. By ergodic decomposition, the determinant is in 
such a case an integral of determinants of type $II_1$ factors, 
in which case $L_k,K_k$ as well as the product $L_k K_k$ are periodic
operators which have their spectrum contained 
in a finite union of real smooth curves in the complex plane. 
In this case, there are
arbitrarily small complex numbers $E_k$ such that $L_k-E_k,K_k-E_k$ are 
invertible allowing to apply the determinant formula. Now, 
$E \mapsto {\rm det}(L_k-E), E \mapsto {\rm det}(K_k-E)$ as well 
as $E \mapsto {\rm det}( (L_k-E)(K_k-E))$ are continuous. 
$L_k-E_k,K_k-E_k$ are invertible and 
${\rm det}(L_k-E_k) \to {\rm det}(L)$, 
${\rm det}(K_k-E_k) \to {\rm det}(K)$, and
${\rm det}( (L_k -E_k) (K_k - E_k) ) \to {\rm det}((L-E)(K-E))$ for 
$n \to \infty$.   \\

(iv) The product formula.  \\
Proof. We can now invoke the determinant theory for finite type 
von Neumann algebras \cite{FuKa52} (see \cite{DixmierF} 
chapter I) which tells 
that ${\rm det}(L K) = {\rm det}(L) {\rm det}(K)$
for invertible operators $L,K$ in a type $II_1$ von Neumann algebra. 
\end{proof} 

\begin{coro}
\label{determinantcoro}
If $L^{(i)}=a^{(i)} \tau + (c^{(i)} \tau)^* +b^{(i)}$ are
random Jacobi matrices over dynamical systems $T^{(i)}: \TT^2 \mapsto \TT^2$, 
then 
$$  {\rm det}(L^{(1)} L^{(2)}) 
  = {\rm det}(L^{(1)}) \; {\rm det}(L^{(2)}) \; . $$
\end{coro}
\begin{proof} 
Apply Lemma~(\ref{determinant}). 
\end{proof} 

Whatever Lyapunov exponent estimates for an analytic cocycle
can be proven for irrational rotations $x \mapsto x + \alpha$ 
as in \cite{Her83,SoSp91}, they can up to some correction $C(\lambda)$ also 
be proven for this cocycle, where the underlying system is an 
arbitrary measure preserving transformation
$T \in \Ycal$. 
For example, for any realanalytic nonconstant $f$, 
the Lyapunov exponent of $A_{T,\lambda f}$ should be
positive for large enough $\lambda$. 

\section{The Jensen formula in a sector} 
\label{a_jensen_formula}
In this section we prove an integrated
form of the Cauchy-Riemann differential equations. 
Unlike the classical Jensen formula which holds in an annulus
(for a short proof see \cite{Garding}), we
formulate it in a sector. This allows to apply it in piecewise analytic 
situations. 
While we don't use this formula explicitly, we used it to 
explain some heuristics in the introduction. We furthermore think it
could be useful for further improvements on the Lyapunov estimates.  \\

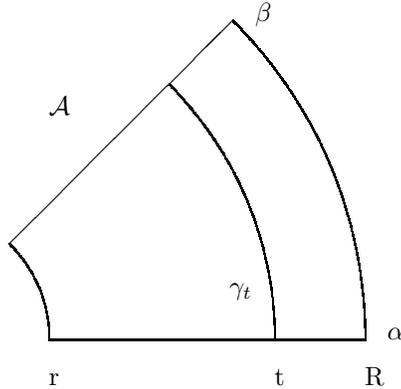
\hbox{ 
\vbox{ \parbox{5.2cm}{ \vspace{-2cm} 
 \setlength{\unitlength}{0.06mm}
 \begin{picture}(1000,1000)
 \qbezier(1000,0)(1000,132)(966,259)
 \qbezier(966,259)(932,386)(866,500)
 \qbezier(866,500)(800,614)(707,707)
 \qbezier(300,0)(300,39)(290,78)
 \qbezier(290,78)(280,116)(260,150)
 \qbezier(260,150)(240,184)(212,212)
 \qbezier(800,0)(800,63)(790,125)
 \qbezier(790,125)(780,187)(761,247)
 \qbezier(761,247)(741,307)(713,363)
 \qbezier(713,363)(684,419)(647,470)
 \qbezier(647,470)(610,521)(566,566)
 \put(300,0){\line(1,0){700}}
 \put(212,212){\line(1,1){495}}
 \put(300,500){${\cal{A}}$}
 \put(700,100){$\gamma_t$}
 \put(1000,-100){R}
 \put(800,-100){t}
 \put(300,-100){r}
 \put(1050,0){$\alpha$}
 \put(757,707){$\beta$}
\end{picture} 
}}
\hspace{-8cm}
\vbox{ \parbox{7.5cm}{
Consider a sector $\Acal=\{ z \in \CC \; | \; |z| \in [r,R], 
{\rm arg}(z) \in [\alpha,\beta] \; \}$ and denote by
$\gamma_t = \{ |z|=t \} \cap \Acal$ a circular arc in $\Acal$. 
Let $g$ be a bounded continuous, complex-valued function on 
$\Acal$ which is analytic and nonzero on $\Acal \setminus \Scal$, 
where $\Scal$ is a finite set of points in $\Acal$. 
Define for $t \in [r,R]$ 
$$ {\rm Arg}_{\gamma_t}(g)
 = {\rm Re} \frac{1}{2\pi i}\int_{\gamma_t} 
                                \frac{g'}{g}(z) \; dz \; .$$
(It is defined except possibly at finitely many $t$, 
where $\gamma_t \cap \Scal$ is nonempty.)
}}}
 
\begin{lemma}[Jensen formula in a sector] 
\label{Jensen}
If $g$ is analytic in $\Acal$ and nonzero in $\Acal \setminus \Scal$, 
where $\Scal$ is a finite set of points, then
$$ \int_\alpha^\beta \log|g(R e^{2\pi i x})| \; dx
 = \int_\alpha^\beta \log|g(r e^{2\pi i x})| \; dx
   +  \int_r^R \frac{dt}{t}{\rm Arg}_{\gamma_t}(g) \; dt\; . $$
\end{lemma}
\begin{proof}
(i) Assume first that $g$ is analytic in $\Acal$ and that it 
has no roots in $\Acal$. Clearly,
$$ \int_\alpha^\beta \log|g(R e^{2\pi i x})| \; dx
 = \int_\alpha^\beta \log|g(r e^{2\pi i x})| \; dx
      + \int_{\Acal} \frac{\partial}{\partial t} \log|g(t e^{2 \pi i x})|
        \; dt \; dx \; . $$
Because $g$ is analytic, the Cauchy-Riemann equations 
$\frac{\partial}{\partial \overline z} g(z)=0$ hold. With
$z=t \exp(2 \pi i x), g'(z) = (\partial/\partial z) g(z)$, we obtain
$$ \frac{\partial}{\partial t} \log|g(t e^{2 \pi i x})| =
   {\rm Re} \left( \frac{g'(t e^{2 \pi i x})}{g(t e^{2 \pi i x})} 
   e^{2 \pi i x} \right) 
     = {\rm Re} \left( \frac{g'(z)}{g(z)} \; \frac{z}{t}  \right) \; . $$
We have therefore, using $2 \pi i dx = dz/z$, 
\begin{eqnarray*}
    \int_\alpha^\beta \log|g(R e^{2\pi i x})| \; dx
  - \int_\alpha^\beta \log|g(r e^{2\pi i x})| \; dx           
   &=& 
     {\rm Re} \frac{1}{2\pi i} 
     \int_r^R \int_{\gamma_t} \frac{g'(z)}{g(z)} \; \frac{dz}{t} \; dt  \\
   &=& \int_r^R {\rm Arg}_{\gamma_t}(g) \; \frac{dt}{t}  \; . 
\end{eqnarray*}

(ii) In general, we partition $\Acal = [r,R] \times [\alpha,\beta]$ into 
$n^2$ small sectors $\Acal_{jk} = [r_j,r_{j+1}] \times [\alpha_k,\alpha_{k+1}]$
and apply the formula to small sectors $\Acal_{jk}$ which do not intersect 
$S$. The assumption on $\Scal$ and $g$
assures that for all the integrals, the contribution of the 
sectors which intersect $\Scal$ can be neglected in the limit 
$n \to \infty$. 
\end{proof}

Remarks. \\
1) The Jensen formula in the annulus \cite{SoSp91,Berenstein} follows: 
if $\Acal$ is the annulus and if $g$ is analytic in 
$\{ {\rm arg}(z) \in (\alpha_i,\alpha_{i+1}), \; |z| \in (r,R) \; \}$ 
and if
$g(z) = \prod_{j=1}^m (z-a_j) h(z)$ such that $h(z) \neq 0$ everywhere 
on $\Acal$, then 
$$ \int_0^{2\pi} \log|g(R e^{2\pi i x})| \; dx
 - \int_0^{2\pi} \log|g(r e^{2\pi i x})| \; dx
 = \int_r^R \frac{dt}{t}{\rm Arg}_{\gamma_t}(h) \; dt 
  + \sum_{j=1}^m \log( \frac{R}{|a_j|} ) \; . $$
2) The Jensen formula for subharmonic functions 
in the annulus $\Acal$ follows: for subharmonic functions 
$z \mapsto \log|g(z)| = \int \log|z-z'| \; dk(z')$ in the disc
$\{ |z|<R \}$, one has
$$ \int_0^{2\pi} \log|g(R e^{2\pi i x})| \; d x
 = \int_0^{2\pi} \log|g(r e^{2\pi i x})| \; d x
      + \int_{\Acal_{r,R}} \log(\frac{R}{|z|}) \; dk(z')  \; . $$
3) Lemma~(\ref{Jensen}) holds for sets $S$ which are
bigger than just a finite point set. For example, $S$ can be a 
subset of a finite union of smooth curves. 

\section{Toral homeomorphisms and Lax approximation}
\label{toral_homeomorphisms_and_lax_approximation}
In this section we review an approximation result
for homeomorphisms on the torus. It allows to approximate the 
group of measure preserving homeomorphisms by a finite group of 
piecewise linear transformations. Investigating the Lyapunov map
$T \mapsto \mu(A_{T,\lambda \sin})$ on finite groups $\Ycal_n$ of 
measure preserving transformations was important for us because the evaluation 
of the Lyapunov exponent is a finite dimensional integration and 
quite reliable. \\
 
On the group $\Ycal$ of all measurable, invertible transformations on 
the torus $\TT^d$ which preserve the Lebesgue measure, one has the metric
$$ \rho(T_1,T_2) = |  {\rm dist}(T_1(x),T_2(x)) |_{\infty} \; , $$
where ${\rm dist}$ is the geodesic distance on the flat torus and 
where $|\cdot|_{\infty}$ is the $L^{\infty}$-norm. \\

A cube exchange transformation on $\TT^d$ is a periodic, piecewise affine
measure-preserving transformation which permutes rigidly all the cubes 
$\prod_{i=1}^d [k_i/n,(k_i+1)/n]$,
where $k_i \in \{0, \dots, n-1 \}$. Every point in $\TT^d$ is periodic. 
Such a transformation is determined by 
a permutation of the set $\{1, \dots, n \}^d$. If the permutation is 
cyclic, the exchange transformation is called cyclic. 
A theorem of Lax \cite{Lax71} states that every $T \in \Xcal$ can be 
approximated in the metric $\rho$ by cube exchange transformations. 
Actually, $T_k$ can be chosen to be cyclic \cite{AlPr93}. 

\begin{lemma}[Lax approximation]
\label{Lax}
For every $\epsilon>0$ and every homeomorphism $T \in \Xcal$, there exists
a cyclic cube exchange transformation $T_k$ which satisfies
$\rho(T_k,T) \leq \epsilon$.
\end{lemma}

Lax's 'book proof' of this result uses Hall's marriage theorem in graph 
theory (for a 'book proof' of the later theorem, see \cite{AigZie}). \\

Remarks.  \\
1)  Periodic approximations of symplectic maps work surprisingly 
well for relatively small $n$ (see \cite{Ran74}). On the Pesin region 
this can be explained in part by the shadowing property \cite{KH}. The 
approximation by cyclic transformations make however long time
stability questions look different \cite{Hal87}. \\

2) It has been measured that certain Lax approximations $T_n$ of the 
Standard map have longest orbits of size $\geq \delta n^2$ with 
$\delta>0$ \cite{ZhVi88}. \\

3) If $T \in \Ycal_k$, the density of states $dk(x,y)$ of a periodic
operator $L(x,y)$ is contained in a finite union of curves.
The density of states of $L$ is the measure
$\int_0^{1/k} \int_0^{1/k}  dk(x,y) \; dx dy$ which has positive area in
the complex plane. For a cyclic $T_k$, the Lyapunov exponent is computed as 
$$ \mu(A_{T_k}) = k^{-2} k^2 \int_0^{1/k} \int_0^{1/k} 
  \log \lambda_{\max}(A^{k^2}(x,y)) \; dx dy  \; , $$
where $\lambda_{max}(A)$ is the maximal eigenvalue of a matrix $A$. The
factor $k^2$ comes because we integrate over a region with Lebesgue measure
$k^{-2}$. 

\section{Thouless formula for nonselfadjoint operators}
\label{thouless_for_nonselfadjoint_operators} 
In this section we look at the Thouless formula for random operators
which can have spectrum with positive area in the complex plane.
Nonselfadjoint operators appear naturally, when 
a problem is parameterized by an external complex parameter like
for operators~(\ref{L(w)}). The proof of the 
Thouless formula which we are looking at in this section 
relies on the fact that for a general parameterization, 
there is an abstract Thouless formula which relies on the Riesz
decomposition theorem in potential theory.  \\

Consider an analytically parameterized, not necessarily selfadjoint,
random Jacobi matrix 
\begin{equation}
\label{Jacobi}
  (L_z(x) u)_n = u_{n+1} + u_{n-1} + f_z(x_n) u_n
\end{equation}
with a complex potential $f_z$ parameterized by a complex parameter $z$.
Examples are the $w$- parameterization 
\begin{equation}
\label{w-parameterisation}
 f_w(x) = (w^{-1} \exp(i x) + w \exp(-i x))/2  
\end{equation}
for $\CC \setminus \{0\}$ and for fixed $w \neq 0$, 
the energy parameterization
\begin{equation}
\label{E-parameterisation}
E \mapsto f_{E,w}(w,x) = f_w(x)-E \; . 
\end{equation}

\begin{propo}[Thouless formula for nonselfadjoint operators]
\label{Thouless}
In the energy parameterization case (\ref{E-parameterisation}), the 
Riesz measure of the Lyapunov exponent $\mu(A_{T,E})$ of the transfer cocycle 
$A_{T,E}$ is the density of states of the operator~(\ref{L(w)}).
\end{propo}
\begin{proof}
For any general analytic parameterization 
$z \mapsto L(z)$, 
the Lyapunov exponent $\mu(A_z)$ of $L(z)$ is a a subharmonic
function of $z$. The Riesz decomposition theorem 
(see i.e. \cite{Hayman,Ransford}) implies the abstract Thouless formula
$\mu(z) = \int \log|z-z'| \; dk(z')  + h(z)$, 
where $z \mapsto h(z)$ is harmonic. \\

The name 'abstract Thouless formula' is from \cite{Her83}. \\

If $z \mapsto \mu(z)$ is harmonic for large $|z|$ and grows like 
$\log|a z|$ for $|z| \mapsto \infty$ then 
$\mu(z) = \log(a) + \int \log|z-z'| \; dk(z')$. \\

The Riesz measure $dk$ of the subharmonic function
$z \mapsto \mu(z)$ satisfies $dk = \Delta \mu$ in the sense of distributions,
where $\Delta$ is the Laplacian on distributions. 
The measure $dk$ coincides in the periodic case with the density of 
states because in this case the Lyapunov exponent is zero exactly on the
spectrum which implies that the density of states is the equilibrium 
measure on the spectrum.  \\

Consider the Riesz measures
$dk_n(z)=\Delta \mu_n(z)$ of the subharmonic functions
$$ z \mapsto \mu_n(z)=n^{-1} \int_{\TT^2} 
     \log|(A^n_{z,T}(x,y))_{11}| \; dx \; dy $$
and let $dk(z)=\Delta \mu(z)$ be the Riesz measure of the Lyapunov exponent 
$z \mapsto \mu(A_{z,T})$.
Because $\mu_n(A_{z,T})(x,y)$ converges by the multiplicative ergodic
theorem to $\mu(A_{z,T})(x,y)$ for almost all $(x,y)$,
the convergence holds also in the sense of distributions.
Then also the Riesz measures $dk_n=\Delta \mu_n$ converge to
$\Delta \mu$ as distributions for $n \to \infty$.
But then, because smooth functions are dense
in all continuous functions, $dk_n=\Delta \mu_n$ converges weakly as 
measures to $dk=\Delta \mu$. \\

In the energy parameterization case~(\ref{E-parameterisation}), 
the density of states of $L$ is defined as the measure $dk$ in the 
parameter plane $\CC$ satisfying 
$\int_{\TT^2} [g(L(x,y))]_{00} \; dx dy = \int g(E') \; dk(E')$
for every continuous function $g$. \\

For $E$ outside the spectrum of $L$, the sequence 
$\int_{\CC} (E'-E)^{-1} \; dk_n(x,y,E')$ converges for $n \to \infty$
uniformly to $\int_{\CC} (E'-E)^{-1} \; dk(x,y,E')$. Because
$$ {\rm tr}((L^{(n)}(x,y)-E)^{-1}) 
  = \int_{\CC} (E'-E)^{-1} \; dk_n(x,y,E') \; , $$
where the Riesz measure $dk_n(x,y) = \Delta \mu_n(x,y)$ is a finite point
measure located on the point spectrum of $L^{(n)}(x,y)$ defined in an 
earlier section, the formula holds also after integration over $\TT^2$ and 
$dk_n$ is the density of states of the random operator
$(x,y) \mapsto L^{(n)}(x,y)$. \\

The functions $g_j(E) = (E-z_j)^{-1}$, where $\{z_j\}_{j \in \NN}$ is a
dense set in the complement of the spectrum of $L$, span a dense set 
in $C(\CC)$. For such functions, we have
$$ \int_{\TT^2} [g_j(L^{(n)}(x,y))]_{00} \; dx dy 
  \mapsto \int_{\TT^2} [g_j(L(x,y))]_{00} \; dx dy $$
by the Avron-Simon lemma. \\

We know therefore
$\int_{\TT^2} [g(L(x,y))]_{00} \; dx dy = \int g(E') \; dk(E')$ 
for all $g \in C(\CC)$ so that $dk$ is the density of states. 
\end{proof}

Examples.\\
1) If $w \mapsto B(w)= w^l A(w)$
is analytic in $\CC$, where $A(w)$ is the transfer cocycle of $L_w$,
then the Lyapunov exponent of $B$ satisfies
$$ \mu(B(w)) = C+\int_{\CC} \log|w-w'| \; dk(w')  \;  $$
where $dk$ has total mass $(l+1)$ and where $C$ is a constant. \\
2) For periodic operators, where the support of the Riesz measure 
is on spectral bands, the Lyapunov exponent of $w \mapsto A(w)$ is
$\mu(w)={\rm Re} [i \arccos({\rm Tr}(A^n)(w)/2)]$. With the parameterization
(\ref{w-parameterisation}), we have
$$ \mu(A(w)) = - \log|w| + \int_{\CC} \log|w-w'| \; dk(w') 
+ \log|\lambda/2| $$
because $w \mapsto \mu(A(w))$ grows like $\log|w \lambda/2|$ for 
$|w| \to \infty$ and $dk(\CC)=2$. 
The Lyapunov exponent of $B(w) = w A(w)$ which is harmonic near $0$
satisfies then 
\begin{equation}
\label{Thouless2}
 \mu(B(w)) = \int_{\CC} \log|w-w'| \; dk(w') + \log|\lambda/2|   \; .
\end{equation} 
We expect the function $w \mapsto \mu(B(w))$, which is harmonic 
near $w=0$ to extend often to a harmonic function in the entire 
complex plane. \\

Remark. By the Hadamard three circle
theorem \cite{Berenstein}, the function $\mu_{max}(r) = \max_{|w|=r} |\mu(w)|$
is a increasing function in $r$ which is convex in $\log(r)$. Aiming
to get rid of the $C(E,\lambda)$ term in Conjecture~(\ref{I}), it would be
good to know $M=\{ |w| = 1  \; | \; \mu(w)=\mu_{max}(|w|) \; \}$ which is 
contained in the subset of the unit circle
on which the potential is $\geq u(0)$. In the case of the Standard map, 
a candidate for a maximum on $|w|=1$ is $w=1$ because the pointwise 
Lyapunov exponent $\mu(x,y,w)$ satisfies $\mu(x,y,w) = \mu(-x,-y,\overline{w})$
so that $\mu(w) = \mu(\overline{w})$. See 18) in the discussion section. 

\section{Log-Holder continuity and capacity}
In this section we consider a potential theoretical property of 
the Riesz measure of the Lyapunov exponent.
This applies to the Riesz measure of the parameterization $w \mapsto L_w$ 
or, for fixed $w$ to the energy parmeterisation $E \mapsto L_w-E$, where 
the Riesz measure $dk$ is the density of states.  \\

As in the selfadjoint case, the fact that the Lyapunov exponent 
is nonnegative implies mild regularity for the measure $dk$. Log-Holder 
continuity is a property of measures with finite potential theoretical
energy. It implies that such a measure is absolutely continuous
with respect to the zero-dimensional Hausdorff measure \cite{Falconer}. 

\begin{lemma}[Log-Holder continuity]
Let $w \mapsto L_w$ be an analytic parameterization of a random 
Jacobi operator. The Riesz measure $dk$ of the nonnegative
subharmonic function 
$w \mapsto \mu(A_T(w))$ is log-Holder continuous: there
exists a constant $C$ such that
$$ dk(B) \leq C \cdot( \log(1/|B|) )^{-1} \;  \;   $$
for every ball $B$ with diameter $|B|$. The constant $C$ only depends
on $L$ and $|E|$.
\end{lemma}
\begin{proof}
The proof is the same as in the
real case \cite{CrSi83,Carmona}: assume $dk$ is contained in a
ball of radius $c$ and assume that $B$ is any ball of radius $r<1$ around
$z \in \CC$. With $A=\{ z' \;  | \;  |z'-z| \geq 1 \; \}$, we have
$$ 0 \leq \int_A \log|z-z'| \; dk(z') \leq \log(1+c+|z|) \; . $$
Therefore, using in the first inequality of the following identities
also that $|z-z'| \leq 1$
on $A^c$ and $\int_{\CC} \log|z-z'| \; dk(z') \geq 0$, we have
$$ 0 \geq \int_{A^c} \log|z-z'| \; dk(z')
     = \int_{\CC} \log|z-z'| dk(z') - \int_{A} \log|z-z'| \; dk(z')
     \geq 0 - \log(1+c+|z|)  \; . $$
For $z_1 \neq z$ with $|z_1-z|=r$, and because $|z-z'|<|z-z_1|$
for $z' \in A^c$
$$ \log|z-z_1| \; dk(B) \geq \int_{A^c} \log|z-z'| \; dk(z')
     \geq - \log(1+c+|z|) $$
so that $0 \leq \log(1+c |z|) +  \log|z-z_1| \; dk(B)$ which is
the claim with $C= \log(1+c |z|)$.
\end{proof}

Example. Consider the case, when $T=T_k$ is
a Lax cube exchange transformation on the torus and where the operator 
$$ L_w(x,y) = \tau+ \tau^* + \lambda/2 (w^{-1} z_j + w z_j^{-1}) $$
with $z_j=z_j(x,y) = \exp(i x_j), (x_j,y_j) = T_k^j(x,y)$
is a periodic, nonselfadjoint Jacobi matrix. 

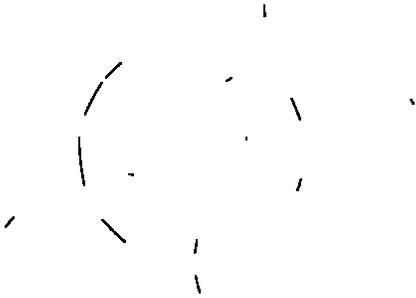
\hbox{
\vbox{ \parbox{6.2cm}{ \vspace{0cm}
 \setlength{\unitlength}{0.03mm}
 \begin{picture}(2000,2000)
 \qbezier(211,708)(191,680)(174,663)
 \qbezier(524,847)(501,960)(503,1059)
 \qbezier(602,1302)(557,1232)(527,1160)
 \qbezier(605,693)(661,633)(704,597)
 \qbezier(620,1326)(663,1369)(685,1388)
 \qbezier(721,897)(733,894)(741,894)
 \qbezier(1023,603)(1019,584)(1014,548)
 \qbezier(1036,372)(1028,400)(1016,447)
 \qbezier(1154,1310)(1163,1316)(1175,1323)
 \qbezier(1244,1057)(1244,1054)(1243,1052)
 \qbezier(1322,1646)(1323,1625)(1324,1599)
 \qbezier(1480,1140)(1464,1187)(1443,1232)
 \qbezier(1484,873)(1477,848)(1470,828)
 \qbezier(1972,1227)(1977,1218)(1983,1209)
 \end{picture}

}}
\hspace{-7cm}
\vbox{ \parbox{6.5cm}{
We are interested in the potential theory of the w-spectrum 
$$ \sigma(L(x,y)) = \{ w \; | \; {\rm tr}(A^n_w(x,y)) 
      \in [-2,2] \; \} \;,$$ 
which supports the 'w-density 
of states' $dk(x,y)$ of $L(x,y)$. In the picture to the left, we 
see an example of the 'w-spectrum' $\sigma(L(x,y))$ for $\lambda=2.1$
where $k=7$. One can see six bands located on the unit circle and eight
bands away from the unit circle. 
}}}

The 'w-density of states' of $L$ is
then $dk = \int_{\TT^2} \; dk(x,y)$ and its support, the 'w-spectrum'
of the random operator $L$ 
is the union $\bigcup_{x,y \in \TT^2} \sigma(L(x,y))$ and has in general 
positive area in the complex plane. \\

The potential theoretical energy 
$$ I(dk(x,y)) = - \int_{\CC} \int_{\CC} \log|z-z'| \; 
                 dk(x,y)(z') \; dk(x,y)(z) $$
is $I(dk(x,y)) = \log(\lambda/2) \; dk(x,y)(\CC) = 2 \log(\lambda/2)$ 
because $dk(x,y)(\CC)=2$ and the Thouless formula assures that 
the Lyapunov exponent $\mu(x,y)(z)$ is related to the potential
$u(x,y)(z) = \int_{\CC} \log|z-z'| \; dk(x,y)(z')$ of $dk(x,y)$
by $\mu(x,y)(z) = \log(\lambda/2) + u(x,y)(z)- \log(w)$ and the Lyapunov 
exponent is zero on the spectrum $\sigma(L(x,y))$. The measure 
$dk(x,y)/2$ is the equilibrium measure on the set $\sigma(L(x,y))$. 
It has the potential theoretical energy $\log(\lambda/2)/2$ 
and the capacity $C(\sigma(L(x,y)))$ is  $\sqrt{\frac{2}{\lambda}}$. 

\section{Harmonic continuation of some singular integrals}
In the case of the Standard map, the
function $\mu(w)=\log(w \; {\det}(L_w))$ with $L_w$ given by (\ref{L(w)})
is for small $|w|$ 
the real part of an analytic function. By the Cauchy-Riemann
differential equations, the ${\rm arg}(w)$-dependence of the function
$w \mapsto \exp(\mu(w))$ is related to the change of the harmonic 
conjugate along the radial direction. The Jensen formula extends
this to the $w$-spectrum of $L_w$, where $\mu$ is 
no more analytic. There are cases (and we expect it to be the case
in the Standard map), in which the harmonic conjugate
can be extended as a harmonic function on a neighborhood of the 
unit disc. Estimating this harmonic conjugate with 
complex analytic methods could lead to possibly better lower bounds 
of the Lyapunov exponent. Alternatively, one could look at the argument 
change (which differs from the harmonic conjugation) and use the 
Jensen formula directly to estimate the oscillation of the Lyapunov 
exponent $w \mapsto \mu(A(w))$ for $w$ the unit circle. 
In this section we formulate the 
problem of computing the harmonic continuation rsp. 
the argument change along radial lines, when the $w$-spectrum is located on the 
unit circle. In the next section, we will see that assuming $dk$ to be 
supported on $\TT$ is no big loss
of generality for estimation purposes 
because projecting the spectrum onto the unit circle 
does not change much the potential (which is the Lyapunov exponent.) \\

A continuous, real-valued periodic function 
$x \mapsto \phi(x)$ on the circle
defines a complex-valued function 
$x \mapsto a(x)=\exp(i (x + \phi(x)))$ and a measure
$dk=h(x) dx$ 
on the circle $\TT$ by pushing forward the Lebesgue measure 
$dk = a^* dx$. Denote by 
$a_n = \int \exp(-i n x) \; h(x) dx$ the Fourier coefficients of this 
measure. 
The function $f(z) = \int_{\TT} dk(x')/(z-\exp(i x'))
= \int_{\TT} dx/(z-a(x))$ defines two 
analytic functions $g_{\pm}(z)$, where $g_+$ is analytic on $|z|<1$ 
and $g_-(z)$ is analytic in $|z|>1$. These functions are defined by their
Taylor expansions at $z_+=0$ and $z_-=\infty$:
\begin{eqnarray*}
 h_+(z) &=& - \int_{\TT} a^{-1}(x) \frac{1}{1-z/a(x)} \; dx  
        = - \sum_{n=0}^{\infty} z^n a_{-(n+1)}, \\
 h_-(z) &=& \int_{\TT} z^{-1} \frac{1}{1-a(x)/z} \; dx    
        = \sum_{n=0}^{\infty} z^{-(n+1)} a_{n}          
        = \sum_{n=-1}^{-\infty} z^n a_{-(n+1)} \; . 
\end{eqnarray*}
The function $h=h_--h_+ = \sum_{n \in \ZZ} z^n a_{-(n+1)}$ is like this
represented by the boundary value of two analytic 
functions $h_{\pm}$. One has for example that
if $h$ is in $L^p$, then $h_{\pm}$ are in the Hardy space
$H^p$. The function $h$ is realanalytic if and only if $h$ is analytic 
in some strip around the unit circle. In this case, we can write
$f(z) = h_+(z) = h(z) + h_-(z)$ for $\{ 1-\epsilon <|z| \leq 1 \}$
and $f(z) = h_-(z) = - h(z) + h_+(z)$ for $\{ 1<|z|<1+\epsilon \}$. 
The functions $h_{\pm}$ can then be analytically continued to 
a larger domain. If $dk$ extends to an entire function,
then both $h_{\pm}$ extend to entire functions.  \\

The integral $H(z)=\int_{\TT} \log(z-a(x)) \; dx$ with $H(0)=0$
can be expanded in the same way at $z=0$ and $\infty$ and 
$$ H=H_- -H_+ 
   = \sum_{n \in \ZZ} z^{n+1} a_{-(n+1)}/(n+1) 
   = \sum_{n \in \ZZ} z^{n} a_{-n}/n \; . $$
The imaginary 
part of $H$ on $|z|=1$ is the average argument change from $0$ to $\infty$ of 
the function $a \mapsto z-a(x)$ along the line going through 
$z=\exp(i \alpha)$. \\

Near $z_-:=0$, the function ${\rm Im}(H_-(z))$ is the argument change
of $\exp(H_-(z))$ from $0$ to $z$. Near $z_+^{-1}:= 0$, the function 
${\rm Im}(H_+(z))$ is the argument change of $\exp(H_+(z))$ from $\infty$
to $z$. \\

Because of the Jensen formula~(\ref{Jensen}), we are interested in the 
argument change $\alpha(z)$ from $0$ to $z$ along a
straight line of the function $\int_{\TT} \log(z-a(x)) \; dx$. 
For each point $a(x)$, the argument change $\alpha(z)$ from $0$ to $z$
is a well defined angle. For $|z|<1$, we have $\alpha(z) = {\rm Im}(H_-(z))$. 
For $|z|>1$ we can write $\alpha(z) = {\rm Im}(H_+(z)) + \beta(z)$, 
where $\beta(z) = \beta(z/|z|)$ is the argument change along the 
straight line from $0$ to $\infty$ which passes through $z$. 
Now $\beta(z) = {\rm Im}(H(z/|z|))$ so that
\begin{eqnarray*}
 \alpha(z) &=&  {\rm Im}(H_+(z)) + H(z/|z|)) \;  , \;     |z| \geq 1   \\
 \alpha(z) &=&  {\rm Im}(H_-(z))             \;  , \;     |z| \leq 1   \; . 
\end{eqnarray*}

\section{Lipshitz approximation of the Lyapunov exponent}
In the theory of Schr\"odinger operators one uses the fact
that the Lyapunov exponent $\mu(A_E)$ of the transfer cocycle $A_E$
of $L$ is the Hilbert transform of 
the integrated density of states (or rotation number \cite{JoMo82,
DeSo83}) $\rho(E)=dk(-\infty,E)$, where $dk$ is the 
density of states of $L$.  For the spectral problem $w \mapsto L_w$, where
the $w$-spectrum is contained in an annulus 
$\Acal_{r}=\{ |z| \in (r,r^{-1}) \; \}, r<1$ containing the unit circle, 
one can ask for the regularity of the 
map $\alpha \mapsto \mu(L(\exp(i \alpha))$. The situation for the 
unit circle in the last section is relevant to this problem.
The idea is to project the measure $dk$ onto the unit circle with the map 
$\pi(z) = {\rm arg}(z)$. While this does not change the Lyapunov 
exponent much, one can use the classical properties for the new 
Lyapunov exponent which is the potential of a measure on the unit circle.
A Lipshitz continuity of the 'integrated density
of states' $\alpha \mapsto dk( \{ z \; | \; {\rm arg}(z) \in (0,\alpha) \} )$
would imply that $\alpha \mapsto \mu(\alpha) = \mu(L(\exp(i \alpha)))$ 
would be close to a Lipshitz continuous map. \\

Let $dk$ be the Riesz measure of the Lyapunov exponent $\mu(w)$
of $w \mapsto L_w$ so that $\mu(w)= \int_{\CC} \log|z-z'| \; dk(z')$. 
Let $\pi^* dk$ be the measure on the unit circle obtained as the 
push-forward of $dk$ under the projection 
$\pi( r \exp(i \phi) ) = \exp(i \phi)$. \\

The fact that the Hilbert transform preserves Lipshitz continuity can be
derived from the Jensen formula in a sector~(\ref{Jensen})
(rsp. the Cauchy-Riemann differential 
equations). Let $H$ be the function as in the last
section with the measure $\pi^* dk$ on the unit circle.
The statement is that if the radial argument change ${\rm Im}(H(z))$ 
of the function $z \mapsto g(z)=\exp(H(z))$
is Lipshitz continuous on the unit circle and the Lipshitz constant is
$C$, then the function ${\rm Re}(H(z))= \int \log|z-z'| \; \pi^* dk(z')$
is Lipshitz continuous on the unit circle with Lipshitz constant $C$. 

\begin{lemma}
If $dk$ is supported on $\Acal_{r} = \{ |z| \in (r,r^{-1}) \}$ and 
$dk$ is the Riesz measure of a Lyapunov exponent,  then 
$$ \int_{\CC} \log|z-z'| \; dk(z') \geq 
   \int_{\CC} \log|z-z'| \; \pi^*(dk)(z') - d(r)  \; . $$
where $d(r) \to 0$ for $r \to 1$. 
\end{lemma}
\begin{proof} 
Given a point $z_0$ with $|z_0|=1$. The set
$Y(z_0,r) = \{ w \in \Acal_{r} \; | \; |w-z_0| > |\pi(w)-z_0| \; \}$ is an 
open set. We have $Y(z_0,r) \subset Y(z_0,r')$ for $r<r'$. Write
$$ \int_{\CC} \log |z-z'| \; dk(z') = 
   \int_{Y(z_0,r)} \log |z-z'| \; dk(z') 
 + \int_{\Acal_r-Y(z_0,r)} \log|z-z'| \; dk(z') \; . $$
We have $\int_{\Acal_r-Y(z_0,r)} \log |z-z'| \; dk(z') 
    \geq \int_{\Acal_r-Y(z_0,r)} \log |z-z'| \; \pi^* dk(z')$ and
$$ \int_{Y(z_0,r)} \log |z-z'| \; dk(z')
    \geq \int_{Y(z_0,r)} \log |z-z'| \; \pi^* dk(z') - d(r) \; ,   $$
where $d(r) \to 0$. 
\end{proof}

A numerical bound on the Lipshitz constant of 
$\alpha \mapsto \int_0^{\alpha} \pi^*(dk)(\beta) \; d\beta
= dk( \{ z \; | \; {\rm arg}(z) \in (0,\alpha) \} )$, 
where $dk$ is the Riesz measure of the Lyapunov exponent 
$w \mapsto \mu(A_{E,\lambda f,T}(w))$ would allow to estimate the Lyapunov 
exponent of the transfer cocycle $A_{E,\lambda f,T}$ of the operator~(\ref{L(w)}).

\section{An estimate for harmonic maps} 
We have numerical evidence that in the case of the Standard map
the Lyapunov exponent $\mu(w)$
is the sum of a positive subharmonic function and a harmonic function 
obtained by harmonic continuation of 
${\rm Re} \int \log(w-w') \; dk(w')+\log|w|$ 
from a neighborhood of $\{ w=0 \}$. 
If this should turn out to be true, 
estimates for harmonic maps could be used to give alternative and maybe
better estimates for the Lyapunov exponent. \\

For a general harmonic function $h$, one can
estimate $\inf_{|z|=r} h(z)$ in terms of $C(r) = \sup_{|z|=r} h(z)$.
We present two ways. The first way uses
Harnack's inequality. A second approach uses the distortion theorem 
for univalent functions. 

\begin{lemma}[Lower bound for harmonic maps]
\label{lowerbound}
Let $h$ be a harmonic map in a disc $\DD(r') = \{ |z| < r' \}$ with 
$r'>r,r'>1$, such that $h(0)=0$ and $h(z) < C(r)$ for all $z \in \DD(r)$. Then 
a) 
$$ h(z) \geq - C(r) \frac{2}{r-1}  \; .  $$
for $|z|=1$ and b) 
$$ h(z) \geq - C(r) (1+ \frac{(1+r^{-1})}{(1-r^{-1})^3} ) $$
for $|z|=1$.
\end{lemma}
\begin{proof}
a) 
$v(z)=C(r)-h(z)$ is nonnegative in the disc $\DD(r) = \{ |z| < r \}$
and $v(0)=C(r)$. Harnack's inequality gives
 $$  C(r)-h(z) = v(z) 
                     \leq C(r) \frac{r+1}{r-1}  \;  $$
so that $h(z)-C(r) \geq - C(r) \frac{r+1}{r-1}$ and 
$h(z) \geq C(r)( 1-\frac{r+1}{r-1} ) = - C(r) \frac{2}{r-1}$.  \\

b) 
The function $C(r)-h(z)$ is positive in $\DD(r)$. 
Let $H(z)$ be the analytic function in $\DD(r)$ 
which has $h$ as its real part and which vanishes at $0$. 
By the Noshiro-Warschawski criterium \cite{Duren}, the 
analytic function
$$ F(z) = \int_0^z (C(r)-H(z)) \; dz       $$
is univalent in $\DD(r)$ because it has a derivative which 
has a positive real part in $\DD(r)$. The function 
$\tilde{F}(z) = F(r z)/(r C(r))$ is in the class
$S$ of functions which are analytic and univalent in the 
unit disc and satisfy $\tilde{F}'(0)=1$. By the distortion 
theorem for maps in the class $S$, we have 
$$ |\tilde{F}'(z)| \leq \frac{(1+R)}{(1-R)^3}    $$
for $|z|=R<1$. Assuming $r>1$ and 
applying this for $R=r^{-1}$, we get for $|z|=R$
$$ |\tilde{F}'(z)| =  |\frac{F'(r z)}{C(r)}| = |1-\frac{H(r z)}{C(r)}| 
       \leq \frac{(1+r^{-1})}{(1-r^{-1})^3} \; . $$
In other words, if $|z|=1$, then 
$$  |1- \frac{H(z)}{C(r)}| \leq \frac{(1+r^{-1})}{(1-r^{-1})^3} $$
which leads to $|H(z)| \leq C(r) (1+\frac{(1+r^{-1})}{(1-r^{-1})^3})$ 
and implies the claim. 
\end{proof} 

\section{Remarks} 

1) {\bf Measure preserving maps}. \\
\hbox{
\vbox{ \parbox{5.2cm}{ \vspace{-3cm}
 \setlength{\unitlength}{0.05mm} 
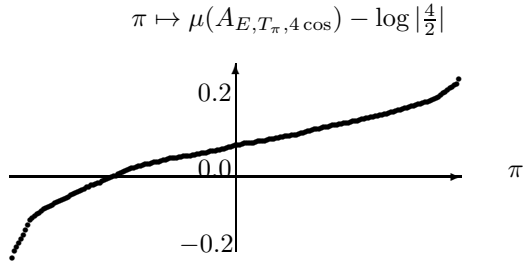
 \begin{picture}(1000,1000)      
 \put(-100,200){\vector(1,0){1200}} 
 \put(500,0){\vector(0,1){500}} 
 \put(1200,200){ $\pi$ } 
 \put(200,600){ $\pi \mapsto \mu(A_{E,T_{\pi}, 4 \cos})-\log|\frac{4}{2}|$ } 
 \put(350,0){$-0.2$ } 
 \put(400,200){$0.0$ } 
 \put(400,400){$0.2$ } 
\cc{-94}{-16}\cc{-88}{1}\cc{-82}{12}\cc{-76}{23}\cc{-70}{33}\cc{-64}{44}
\cc{-58}{58}\cc{-52}{71}\cc{-46}{84}\cc{-40}{91}\cc{-34}{95}\cc{-28}{100}
\cc{-22}{104}\cc{-16}{109}\cc{-10}{113}\cc{-4}{116}\cc{2}{119}\cc{8}{123}
\cc{14}{126}\cc{20}{129}\cc{26}{132}\cc{32}{135}\cc{38}{138}\cc{44}{141}
\cc{50}{144}\cc{56}{147}\cc{62}{151}\cc{68}{154}\cc{74}{157}\cc{80}{160}
\cc{86}{163}\cc{92}{165}\cc{98}{168}\cc{104}{171}\cc{110}{173}\cc{116}{176}
\cc{122}{178}\cc{128}{180}\cc{134}{182}\cc{140}{184}\cc{146}{187}\cc{152}{190}
\cc{158}{193}\cc{164}{196}\cc{170}{198}\cc{176}{201}\cc{182}{203}\cc{188}{206}
\cc{194}{209}\cc{200}{212}\cc{206}{215}\cc{212}{217}\cc{218}{220}\cc{224}{222}
\cc{230}{224}\cc{236}{226}\cc{242}{228}\cc{248}{229}\cc{254}{231}\cc{260}{232}
\cc{266}{234}\cc{272}{235}\cc{278}{236}\cc{284}{238}\cc{290}{239}\cc{296}{240}
\cc{302}{242}\cc{308}{243}\cc{314}{244}\cc{320}{246}\cc{326}{248}\cc{332}{249}
\cc{338}{250}\cc{344}{251}\cc{350}{252}\cc{356}{253}\cc{362}{254}\cc{368}{255}
\cc{374}{255}\cc{380}{256}\cc{386}{257}\cc{392}{258}\cc{398}{259}\cc{404}{261}
\cc{410}{262}\cc{416}{263}\cc{422}{264}\cc{428}{265}\cc{434}{267}\cc{440}{268}
\cc{446}{269}\cc{452}{270}\cc{458}{272}\cc{464}{273}\cc{470}{275}\cc{476}{276}
\cc{482}{278}\cc{488}{280}\cc{494}{281}\cc{500}{283}\cc{506}{284}\cc{512}{285}
\cc{518}{286}\cc{524}{288}\cc{530}{289}\cc{536}{290}\cc{542}{290}\cc{548}{291}
\cc{554}{292}\cc{560}{294}\cc{566}{294}\cc{572}{295}\cc{578}{296}\cc{584}{298}
\cc{590}{299}\cc{596}{300}\cc{602}{301}\cc{608}{302}\cc{614}{303}\cc{620}{304}
\cc{626}{306}\cc{632}{307}\cc{638}{308}\cc{644}{309}\cc{650}{310}\cc{656}{311}
\cc{662}{312}\cc{668}{313}\cc{674}{315}\cc{680}{316}\cc{686}{318}\cc{692}{320}
\cc{698}{321}\cc{704}{323}\cc{710}{325}\cc{716}{327}\cc{722}{328}\cc{728}{330}
\cc{734}{331}\cc{740}{332}\cc{746}{333}\cc{752}{335}\cc{758}{336}\cc{764}{337}
\cc{770}{338}\cc{776}{339}\cc{782}{340}\cc{788}{342}\cc{794}{343}\cc{800}{344}
\cc{806}{345}\cc{812}{347}\cc{818}{348}\cc{824}{349}\cc{830}{351}\cc{836}{353}
\cc{842}{354}\cc{848}{355}\cc{854}{356}\cc{860}{358}\cc{866}{359}\cc{872}{360}
\cc{878}{362}\cc{884}{363}\cc{890}{365}\cc{896}{366}\cc{902}{368}\cc{908}{369}
\cc{914}{371}\cc{920}{373}\cc{926}{375}\cc{932}{376}\cc{938}{378}\cc{944}{380}
\cc{950}{381}\cc{956}{383}\cc{962}{385}\cc{968}{387}\cc{974}{389}\cc{980}{390}
\cc{986}{392}\cc{992}{393}\cc{998}{395}\cc{1004}{397}\cc{1010}{399}\cc{1016}{401}
\cc{1022}{403}\cc{1028}{406}\cc{1034}{409}\cc{1040}{413}\cc{1046}{418}\cc{1052}{422}
\cc{1058}{425}\cc{1064}{431}\cc{1070}{437}\cc{1076}{440}\cc{1082}{444}\cc{1088}{447}
\cc{1094}{459} \put(-150,-170){\small{All cyclic permutations $T_{\pi}$ of $\{1, \cdots, 9 \}$ }. They} 
 \put(-150,-240){\small{are ordered such that the excess to}} 
 \put(-150,-310){\small{$\log|\lambda/2|$, shown here for $\lambda=4$, increases.}} 
 \end{picture}         
}}
\hspace{-8cm}
\vbox{ \parbox{7.5cm}{
The random variables $x_n$ where $T^n(x,y)=(x_n,y_n)$
have as their law the Lebesgue measure. One could enlarge the
class of transformations even so it is not much gain of generality  
\cite{LiTh77}. The conjectured lower bound $\log(\lambda/2)$ 
can be expected to hold in general only for homeomorphisms 
$T \in \Xcal$. The reason is that for some piecewise periodic measure
preserving transformations $T_{\pi}: (x,y) \mapsto
(x+\pi(x)/n \; {\rm mod} \; 1,y)$ with a permutation 
$\pi$ on $\{ 1,2, \dots , n \}$, 
the value of $\mu(A_{T, \lambda \cos})$ can be slightly smaller than
$\log(\lambda/2)$. The Figure to the left shows the distribution of 
the random variable
$\pi \mapsto \mu(A_{T_{\pi},\lambda \cos})$ on the group 
$\Ycal_9 \subset \Ycal$ of all permutations of the 9 annuli
$\{ (j/9,(j+1)/9 ] \times \TT \subset \TT^2 \}_{j=1}^9 $ in the case
$\lambda=4.0$. \\
}}} 

 It would be interesting to know
${\rm inf}_{T \in \Ycal} \mu(A_{T, \lambda \cos})$ for
given $\lambda>0$ and whether the infimum is attained. By upper
continuity of the Lyapunov exponent the 
infimum taken on the discrete set
$\bigcup_n \Ycal_n \subset \Ycal$ is the same as the infimum taken 
on $\Xcal$. \\

2) {\bf Hamiltonian flows}. \\
Like for a large class of monotone twist maps, 
any Standard map is the time-one map of 
a Hamiltonian differential equation with time-dependent periodic 
Hamiltonians \cite{Mos88}. This is also true in higher dimensional cases
\cite{BiPo92,Gol95}. If positive Kolmogorov-Sinai entropy is stable under perturbations of the map,
there would be many real-analytic, time-dependent periodic potentials $V$ 
for which the time-dependent Hamiltonian system $\ddot{x} = V(t,x)$ 
has positive Kolmogorov-Sinai entropy. To see this, we note that the potential 
$V(t,x) = \delta_{t \in \ZZ} \lambda \sin_m(x)$ (which is a distribution 
in $t$) has as a Poincar\'e map the Standard map $T_{\lambda \sin_m}$. The
smoothed-out potentials $V_{\epsilon}(t,x) = \phi_{\epsilon}(t) * V(t,x)$,
where $\phi_{\epsilon}$ is a mollifier function, has a Poincar\'e map 
which is realanalytic in $(x,y)$ and which is arbitrarily close to 
$T_{\lambda \sin_m}$ in a Banach space of realanalytic maps. 
For fixed $\lambda$, the Poincar\'e map has now positive metric entropy for
small $\epsilon$ and large $m$ because the stability result 
as stated in Conjecture~(\ref{I}) extends to general realanalytic 
perturbations of the map (and not only to perturbations of the function 
$f$).  It would be nice to see entropy estimates
for Hamiltonian flows on three dimensional energy surfaces 
of systems $\ddot{x}=V(x)$, where $V(x)$ is 
a potential on $\TT^2$ or to geodesic flows on the torus $\TT^2$ or to 
realanalytic, strictly convex Birkhoff billiards for which we have no
example, for which positive Kolmogorov-Sinai entropy has been proven. \\
Even more off-limit seem problems in celestial mechanics like 
the St\"ormer problem of particles trapped in the van Allen belts 
of the Earth's magnetic dipole field \cite{Bra81} 
(for which Pesin theory in the form \cite{KS} can be shown to apply) 
or particular Newtonian three body problems.  \\ 

3) {\bf Smaller $\lambda$}. \\
Estimating the entropy in the case of small $\lambda$ is also open 
for the Chirikov Standard map.  Experiments indicate that $\mu(\lambda)>0$ for all 
$\lambda>0$. \\

Numerically, we find $\mu(\lambda) \geq \log(\lambda/2)$ for all $\lambda$
and as a "rule of thumb" that 
the averaged Lyapunov exponent $\mu(\lambda)$ satisfies almost but not
exactly the Aubry duality 
$\mu(\lambda) \sim \log(\lambda/2) + \mu(4/\lambda)$. 

\hbox{
\vbox{ \parbox{5.2cm}{ \vspace{-3cm}
 \newcommand{\ff}[2]{\put(#1,#2){\circle*{20}}}  
 \setlength{\unitlength}{0.05mm} 
 \begin{picture}(1000,1000)      
 \put(500,0){\vector(1,0){600}} 
 \put(500,0){\vector(-1,0){600}} 
 \put(500,0){\vector(0,1){500}} 
 \put(1150,-50){ $\lambda$ } 
 \put(-250,-50){ $4/\lambda$ } 
 \put(200,600){ $\mu(A_{\lambda \cos})-\log^+|\frac{\lambda}{2}|$ } 
 \put(490,-60){\small{2}} 
 \put(625,-60){\small{4}} 
 \put(760,-60){\small{6}} 
 \put(895,-60){\small{8}} 
 \put(1030,-60){\small{10}} 
 \put(355,-60){\small{4}} 
 \put(220,-60){\small{6}} 
 \put(85,-60){\small{8}} 
 \put(-50,-60){\small{10}} 
 \put(350,400){\small{0.4}} 
 \ff{-50}{0} 
 \ff{-13}{0} 
 \ff{23}{1} 
 \ff{60}{0} 
 \ff{97}{1} 
 \ff{133}{2} 
 \ff{170}{2} 
 \ff{207}{4} 
 \ff{243}{6} 
 \ff{280}{12} 
 \ff{317}{22} 
 \ff{353}{42} 
 \ff{390}{83} 
 \ff{427}{141} 
 \ff{463}{224} 
 \ff{500}{331} 
 \ff{537}{250} 
 \ff{573}{172} 
 \ff{610}{104} 
 \ff{647}{54} 
 \ff{683}{43} 
 \ff{720}{41} 
 \ff{757}{36} 
 \ff{793}{25} 
 \ff{830}{22} 
 \ff{867}{24} 
 \ff{903}{18} 
 \ff{940}{15} 
 \ff{977}{14} 
 \ff{1013}{12} 
 \ff{1050}{11} 
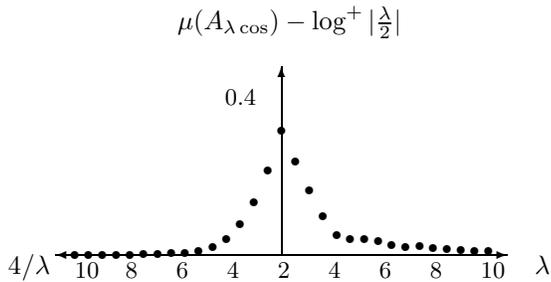
 \put(-150,-170){\small{How close are we to Aubry duality?}} 
 \put(-150,-240){\small{For each $\lambda$ we averaged                                  70x70=4900 different}} 
 \put(-150,-310){\small{orbits of length $10^6$ on the phase space.}} 
 \end{picture}         
}}
\hspace{-8cm}
\vbox{ \parbox{7.5cm}{
\vspace{5mm}
The duality holds in the case of the dynamical system 
$T_0$ (the almost Mathieu case). It does not hold for $T_{\lambda \sin}$
as the Figure to the left shows for $\lambda \in [4/10,2]$ and 
$\lambda \in [2,10]$. 
Indeed, by computing moments of the density of states of $L$ using
a random walk expansion \cite{KniI},
one can check that the discrete random Schr\"odinger operator
$L_{T, \lambda \cos}= \Delta + \lambda \cos(x)$ 
over the dynamical system $T_{\lambda}$
does not have the same density of states as the naive "dual operator"
$\tilde{L} = (\lambda/2) \Delta + 2 \cos(x)
= (\lambda/2) L_{T,(4/\lambda) \cos}$.
(It is not excluded however
that changing $T$ in $\tilde{L}$ could restore a generalized duality.) \\
}}} 
\vspace{5mm}
As in the almost Mathieu case, there is a two-dimensional magnetic
operator involved. The magnetic field at a plaquette $P_{n,m}$
is $\exp(i (x_n-x_{n-1})$ independent of $m$. 
In the Standard map case, the magnetic fields depend on space and are 
correlated. See 25) for more on the Aubry duality. \\

4) {\bf Smooth cocycles}. \\
Results like in \cite{You93,You97} depend 
on the dynamical system, e.g. Diophantine conditions for irrational rotations. 
It is known that positive Lyapunov exponents of $SL(2,\RR)$ cocycles 
can be destroyed away from Anosov systems by measurable 
perturbations of the cocycle \cite{Kni91}, by continuous perturbations 
\cite{Fur97} or using the unproven "Last theorem of Man\'e"
\cite{Man96}, by $C^1$-perturbations of a non-Anosov map. \\

5) {\bf Equilibrium measures}. \\
Via the variational principle (see e.g. \cite{Walters}) metric entropy results
give lower bounds on the topological entrop.
with previous known lower bounds \cite{Dua94,Kni96}. The estimates in
\cite{Kni96} are obtained using a Gronwall argument which estimates the 
spectral radius of the operator $L^{-1}$ (where $L$ is defined by 
a hyperbolic invariant measure of the Chirikov map) and which are far from
optimal (the bounds in \cite{Kni96} are less good than the estimates 
in \cite{Dua94} but are proven when $f$ is a nonconstant Morse function 
and not only in the case of the Standard map where $f=\sin$.) 
The equilibrium measure is a $T$-invariant measure on $\TT^2$
which maximizes the entropy and has as the metric entropy 
the topological entropy.
In the case of the Standard map, it is not known whether there
is an absolutely continuous equilibrium measure. \\

6) {\bf Absolutely continuous spectrum}. \\
\vspace{0mm}
\hbox{
\vbox{ \parbox{5.2cm}{ \vspace{-2cm}
 \setlength{\unitlength}{0.05mm} 
 \begin{picture}(1000,1000)      
 \put(500,0){\vector(1,0){600}} 
 \put(500,0){\line(-1,0){600}} 
 \put(500,0){\vector(0,1){500}} 
 \put(1150,-50){ $E$ } 
 \put(200,600){ $\mu(A_{E,10 \cos})-\log|\frac{10}{2}|$ } 
 \put(490,-80){\small{-2}} 
 \put(718,-80){\small{3.0}} 
 \put(945,-80){\small{8}} 
 \put(248,-80){\small{-7.0}} 
 \put(6,-80){\small{-12}} 
\cc{-94}{694}\cc{-88}{663}\cc{-82}{638}\cc{-76}{609}\cc{-70}{577}\cc{-64}{543}
\cc{-58}{508}\cc{-52}{470}\cc{-46}{435}\cc{-40}{397}\cc{-34}{352}\cc{-28}{309}
\cc{-22}{263}\cc{-16}{207}\cc{-10}{156}\cc{-4}{92}\cc{2}{94}\cc{8}{98}
\cc{14}{104}\cc{20}{102}\cc{26}{101}\cc{32}{98}\cc{38}{96}\cc{44}{93}
\cc{50}{82}\cc{56}{79}\cc{62}{71}\cc{68}{62}\cc{74}{50}\cc{80}{48}
\cc{86}{41}\cc{92}{37}\cc{98}{28}\cc{104}{18}\cc{110}{13}\cc{116}{5}
\cc{122}{11}\cc{128}{18}\cc{134}{13}\cc{140}{16}\cc{146}{21}\cc{152}{25}
\cc{158}{25}\cc{164}{22}\cc{170}{23}\cc{176}{23}\cc{182}{21}\cc{188}{19}
\cc{194}{19}\cc{200}{22}\cc{206}{21}\cc{212}{20}\cc{218}{20}\cc{224}{23}
\cc{230}{21}\cc{236}{20}\cc{242}{16}\cc{248}{16}\cc{254}{11}\cc{260}{14}
\cc{266}{14}\cc{272}{13}\cc{278}{13}\cc{284}{14}\cc{290}{16}\cc{296}{16}
\cc{302}{15}\cc{308}{11}\cc{314}{13}\cc{320}{14}\cc{326}{15}\cc{332}{15}
\cc{338}{16}\cc{344}{13}\cc{350}{13}\cc{356}{12}\cc{362}{12}\cc{368}{13}
\cc{374}{13}\cc{380}{13}\cc{386}{14}\cc{392}{13}\cc{398}{16}\cc{404}{14}
\cc{410}{14}\cc{416}{11}\cc{422}{13}\cc{428}{14}\cc{434}{12}\cc{440}{11}
\cc{446}{14}\cc{452}{11}\cc{458}{11}\cc{464}{14}\cc{470}{11}\cc{476}{11}
\cc{482}{14}\cc{488}{11}\cc{494}{12}\cc{500}{13}\cc{506}{14}\cc{512}{14}
\cc{518}{14}\cc{524}{11}\cc{530}{13}\cc{536}{12}\cc{542}{11}\cc{548}{13}
\cc{554}{10}\cc{560}{12}\cc{566}{15}\cc{572}{11}\cc{578}{12}\cc{584}{11}
\cc{590}{14}\cc{596}{12}\cc{602}{12}\cc{608}{14}\cc{614}{12}\cc{620}{12}
\cc{626}{12}\cc{632}{13}\cc{638}{13}\cc{644}{12}\cc{650}{14}\cc{656}{15}
\cc{662}{15}\cc{668}{15}\cc{674}{14}\cc{680}{15}\cc{686}{14}\cc{692}{12}
\cc{698}{13}\cc{704}{13}\cc{710}{15}\cc{716}{15}\cc{722}{14}\cc{728}{15}
\cc{734}{13}\cc{740}{14}\cc{746}{17}\cc{752}{17}\cc{758}{15}\cc{764}{19}
\cc{770}{18}\cc{776}{19}\cc{782}{19}\cc{788}{20}\cc{794}{22}\cc{800}{20}
\cc{806}{22}\cc{812}{21}\cc{818}{20}\cc{824}{25}\cc{830}{23}\cc{836}{23}
\cc{842}{26}\cc{848}{23}\cc{854}{26}\cc{860}{22}\cc{866}{17}\cc{872}{17}
\cc{878}{11}\cc{884}{5}\cc{890}{13}\cc{896}{19}\cc{902}{25}\cc{908}{34}
\cc{914}{39}\cc{920}{49}\cc{926}{52}\cc{932}{62}\cc{938}{75}\cc{944}{75}
\cc{950}{86}\cc{956}{93}\cc{962}{96}\cc{968}{97}\cc{974}{108}\cc{980}{103}
\cc{986}{101}\cc{992}{100}\cc{998}{95}\cc{1004}{90}\cc{1010}{152}\cc{1016}{209}
\cc{1022}{258}\cc{1028}{306}\cc{1034}{351}\cc{1040}{389}\cc{1046}{433}\cc{1052}{470}
\cc{1058}{509}\cc{1064}{544}\cc{1070}{579}\cc{1076}{609}\cc{1082}{638}\cc{1088}{667}
\cc{1094}{692}\cc{1100}{717} 
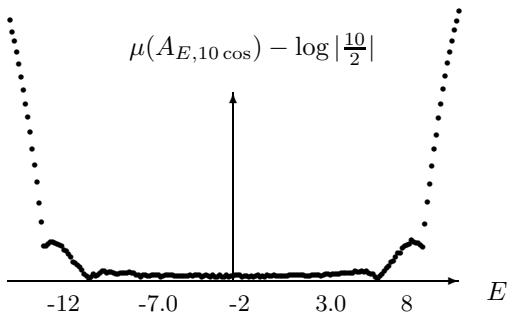
 \put(-150,-170){\small{The Lyapunov exponent of the Chirikov Hamiltonian.}} 
 \put(-150,-240){\small{The graph shows a nonnegative excess to $\log|\lambda/2|$}}
 \put(-159,-310){\small{for all $E$.}} 
 \end{picture}         
}}
\hspace{-8cm}
\vbox{ \parbox{7.5cm}{
Lyapunov exaponents have applications in the theory of discrete random
Schr\"odinger operators \cite{Cycon,Carmona,Pastur}. The
attribute 'random' is used in that theory often in the same way as
'random variable' is used in probability theory. Indeed, random
operators are classes of operator-valued random variables.
There is in general no randomness
assumed in the potential. For example, if $T$ is a translation on the
torus, the random operators are quasiperiodic. 
Using Pastur-Kotani theory, Conjecture~\ref{Igeneral} enlarges
the class of ergodic
discrete one-dimensional random Schr\"odinger operators $L$, 
for which there is no absolutely continuous spectrum.  \\
}}}
\vspace{4mm}

Conjecture~\ref{Igeneral}
implies with \cite{OxUl41} that for a Baire generic ergodic 
$T \in \Xcal$, 
the corresponding ergodic operators have no 
absolutely continuous spectrum for large $\lambda$. 
\footnote{The space $\Xcal$ becomes a complete metric space with
the metric $\rho(T,S) + \rho(T^{-1},S^{-1})$.}
For non-ergodic $T$, the absolutely continuous spectrum is only absent for a 
large set of $x \in \TT^d$ because the Lyapunov exponent can be zero on a 
set of positive measure. Indeed, for invariant measures on KAM tori, the
corresponding ergodic almost periodic Schr\"odinger operator can have some 
absolutely continuous spectrum for small $\lambda$ \cite{KnLa94}. 
\footnote{The results in \cite{KnLa94} have not been written down in a 
final form since both of us got engaged in other projects.}
The size of the spectrum, where the Lyapunov exponent vanishes can be 
estimated if the orbits do not stay too long at the same $x$ value
for a long time \cite{Sur94}.  \\

7) {\bf Discrete Spectrum}. \\
The question of localization of the operators $L_{\lambda \cos}$ 
\cite{Spe89} stays open. 
The conjugation of the map to a Bernoulli shift on some
positive Lebesgue measure suggests that eigenvalues should exist for 
similar reasons as in the case of independent, identically distributed 
potentials. \\

\hbox{
\vbox{ \parbox{5.2cm}{ \vspace{-3.5cm}
 \setlength{\unitlength}{0.05mm} 
 \begin{picture}(1000,1000)      
 \put(600,0){\vector(1,0){700}} 
 \put(-200,0){\vector(1,0){700}} 
 \put(-200,0){\vector(0,1){500}} 
 \put(600,0){\vector(0,1){500}} 
 \put(1250,-100){ $\lambda$ } 
 \put(400,-100){ $\lambda$ } 
 \put(-100,-70){{\small $2.0$}} 
 \put(250,-70){{\small $8.0$}} 
 \put(700,-70){{\small $2.0$}} 
 \put(1050,-70){{\small $8.0$}} 
 \put(-350,350){{\small $0.05$}} 
 \put(500,350){{\small $0.1$}} 
 \put(150,600){$n^{-1} \sum_{k=1}^n |\hat{\mu}_k(L_{\lambda \cos})|^2$} 
\cc{633}{2}\cc{667}{4}\cc{700}{36}\cc{733}{108}\cc{767}{146}
\cc{800}{187}\cc{833}{196}\cc{867}{211}\cc{900}{231}\cc{933}{281}
\cc{967}{310}\cc{1000}{325}\cc{1033}{330}\cc{1067}{311}\cc{1100}{319}
\cc{1133}{338}\cc{1167}{324}\cc{1200}{314}\cc{-167}{0}\cc{-133}{0}\cc{-100}{0}\cc{-67}{1}\cc{-33}{38}
\cc{0}{302}\cc{33}{305}\cc{67}{279}\cc{100}{242}\cc{133}{208}
\cc{167}{178}\cc{200}{154}\cc{233}{134}\cc{267}{118}\cc{300}{104}
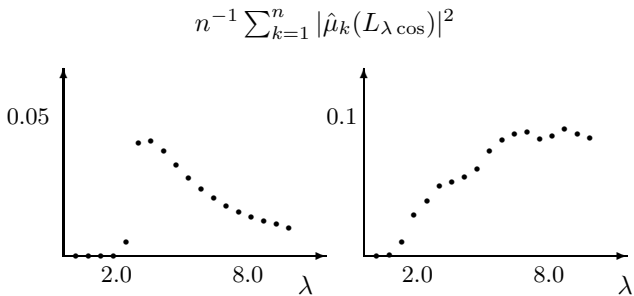
\cc{333}{93}\cc{367}{84}\cc{400}{76} \put(-250,-270){\small{Wiener test for localisation for $\lambda \in [0,10]$. To the left, the}} 
 \put(-250,-340){\small{Mathieu case (n=100'000), to the right, the Chiricov case}} 
 \put(-250,-410){\small{(n=10'000 averaged over 50x50=2500 orbits).}} 
 \end{picture}         
}}
\hspace{-8cm}
\vbox{ \parbox{7.5cm}{
In the Figure to the left, numerical results are shown for
the Fourier coefficients $\hat{\mu}_n=(\phi, \exp(-i n L) \phi)
= \int_{\TT} e^{-in \theta} \; d\mu_{\phi}(\theta)$
of the spectral measure $\mu_{\phi}$
where $\phi = (\dots, 0,1,0, \dots ) \in l^2(\ZZ,\CC)$. We see the 
case of the Standard map operator and Mathieu operator. 
With a Wiener criterion for the existence of discrete 
spectrum, (where we compute a fast discrete quantum evolution using a
related operator which has the same type of spectral measures
\cite{Kni98}), we find numerical evidence for some point spectrum 
for $L(x,y)$ for a set of $(x,y) \in \TT^2$ of positive measure, 
if $\lambda$ is large.  \\

In the picture to the left, we confirm numerically the theoretically 
established localization transition for the Mathieu operator at 
$\lambda=2$ \cite{Jit94,Jit95}. In the picture to the right
we see an indication for some point 
spectrum in the Chirikov case for $\lambda>0$.   \\
}}}
\vspace{1cm}
Note that in the same way as for the Anderson operator,  
the operator $L_T(x,y)$ has for an Anosov $T$ purely 
singular continuous spectrum for a Baire generic set of $(x,y)$. This can
be derived from general principles \cite{Sim95}. Also, for Baire generic
$(T,(x,y)) \in \Xcal \times \TT^2$ and $\lambda>\lambda_0$,
one knows that $L_{T,\lambda \cos}(x,y)$ has purely singular continuous 
spectrum. This means that $L(x,y)$ on $l^2(\ZZ,\CC)$ 
has Baire generically no bound states and no extended states.  \\

8) {\bf The distribution of the Lyapunov exponents}. \\
\hbox{
\vbox{ \parbox{5.2cm}{ \vspace{-3cm}
 \setlength{\unitlength}{0.05mm} 
 \begin{picture}(1000,1000)      
 \put(600,0){\vector(1,0){700}} 
 \put(-200,0){\vector(1,0){700}} 
 \put(-200,0){\vector(0,1){500}} 
 \put(600,0){\vector(0,1){500}} 
 \put(1250,-100){ $t$ } 
 \put(400,-100){ $t$ } 
 \put(200,600){$m( \{ (x,y) \; | \; X_{\lambda}^*(x,y) \leq t \}$)} 
\cc{606}{0}\cc{612}{0}\cc{618}{0}\cc{624}{0}\cc{630}{1}
\cc{636}{1}\cc{642}{1}\cc{648}{1}\cc{654}{1}\cc{660}{1}
\cc{666}{1}\cc{672}{2}\cc{678}{2}\cc{684}{2}\cc{690}{3}
\cc{696}{4}\cc{702}{4}\cc{708}{6}\cc{714}{7}\cc{720}{8}
\cc{726}{10}\cc{732}{12}\cc{738}{14}\cc{744}{16}\cc{750}{20}
\cc{756}{22}\cc{762}{26}\cc{768}{31}\cc{774}{34}\cc{780}{41}
\cc{786}{48}\cc{792}{55}\cc{798}{61}\cc{804}{68}\cc{810}{76}
\cc{816}{86}\cc{822}{95}\cc{828}{104}\cc{834}{115}\cc{840}{127}
\cc{846}{139}\cc{852}{152}\cc{858}{164}\cc{864}{177}\cc{870}{189}
\cc{876}{202}\cc{882}{216}\cc{888}{228}\cc{894}{243}\cc{900}{258}
\cc{906}{272}\cc{912}{286}\cc{918}{298}\cc{924}{313}\cc{930}{327}
\cc{936}{342}\cc{942}{356}\cc{948}{368}\cc{954}{382}\cc{960}{392}
\cc{966}{402}\cc{972}{412}\cc{978}{421}\cc{984}{428}\cc{990}{435}
\cc{996}{441}\cc{1002}{448}\cc{1008}{455}\cc{1014}{460}\cc{1020}{465}
\cc{1026}{470}\cc{1032}{474}\cc{1038}{478}\cc{1044}{481}\cc{1050}{484}
\cc{1056}{486}\cc{1062}{488}\cc{1068}{490}\cc{1074}{492}\cc{1080}{493}
\cc{1086}{495}\cc{1092}{496}\cc{1098}{497}\cc{1104}{497}\cc{1110}{498}
\cc{1116}{499}\cc{1122}{499}\cc{1128}{499}\cc{1134}{499}\cc{1140}{499}
\cc{1146}{499}\cc{1152}{500}\cc{1158}{500}\cc{1164}{500}\cc{1170}{500}
\cc{1176}{500}\cc{1182}{500}\cc{1188}{500}\cc{1194}{500}\cc{1200}{500}
\cc{-194}{121}\cc{-188}{121}\cc{-182}{121}\cc{-176}{121}\cc{-170}{121}
\cc{-164}{121}\cc{-158}{121}\cc{-152}{121}\cc{-146}{121}\cc{-140}{121}
\cc{-134}{121}\cc{-128}{121}\cc{-122}{121}\cc{-116}{121}\cc{-110}{121}
\cc{-104}{121}\cc{-98}{121}\cc{-92}{121}\cc{-86}{121}\cc{-80}{121}
\cc{-74}{121}\cc{-68}{121}\cc{-62}{121}\cc{-56}{121}\cc{-50}{121}
\cc{-44}{121}\cc{-38}{121}\cc{-32}{121}\cc{-26}{121}\cc{-20}{121}
\cc{-14}{121}\cc{-8}{121}\cc{-2}{121}\cc{4}{121}\cc{10}{121}
\cc{16}{121}\cc{22}{121}\cc{28}{121}\cc{34}{121}\cc{40}{121}
\cc{46}{121}\cc{52}{121}\cc{58}{121}\cc{64}{121}\cc{70}{121}
\cc{76}{121}\cc{82}{121}\cc{88}{121}\cc{94}{121}\cc{100}{121}
\cc{106}{121}\cc{112}{121}\cc{118}{121}\cc{124}{121}\cc{130}{121}
\cc{136}{121}\cc{142}{121}\cc{148}{121}\cc{154}{121}\cc{160}{121}
\cc{166}{121}\cc{172}{121}\cc{178}{121}\cc{184}{121}\cc{190}{122}
\cc{196}{122}\cc{202}{122}\cc{208}{122}\cc{214}{122}\cc{220}{122}
\cc{226}{122}\cc{232}{122}\cc{238}{122}\cc{244}{122}\cc{250}{122}
\cc{256}{122}\cc{262}{123}\cc{268}{123}\cc{274}{124}\cc{280}{124}
\cc{286}{124}\cc{292}{124}\cc{298}{125}\cc{304}{126}\cc{310}{126}
\cc{316}{126}\cc{322}{127}\cc{328}{128}\cc{334}{129}\cc{340}{131}
\cc{346}{135}\cc{352}{141}\cc{358}{153}\cc{364}{180}\cc{370}{223}
\cc{376}{286}\cc{382}{367}\cc{388}{450}\cc{394}{492}\cc{400}{500}
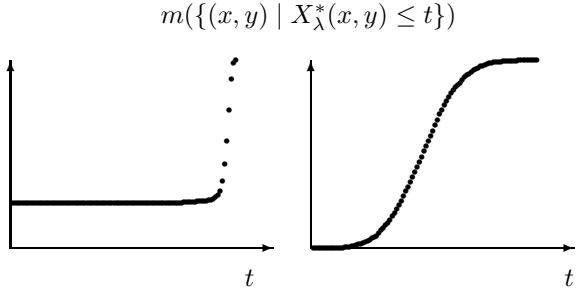
 \put(-150,-270){\small{The distribution functions of the random variable $X_{\lambda}^*$ }} 
 \put(-150,-340){\small{for $\lambda=2,10$. For $\lambda=2$, $m(X_{\lambda}=0)>0$ (KAM). }} 
 \end{picture}         
}}
\hspace{-8cm}
\vbox{ \parbox{7.5cm}{
For every $\lambda$, the pointwise Lyapunov exponent
$(x,y) \mapsto \mu_{\lambda}(x,y)$ is a random 
variable on $\TT^2$ which is $T$-invariant and has a mean 
$\geq \log(|\lambda|/2)-C(\lambda)$.
The numerical experiments indicate that the 
random variables $X_{\lambda}(x,y) = \mu_{\lambda}(x,y)$ have the 
property that the normalizations 
$X^*_{\lambda} 
= (X_{\lambda}-{\rm E}[X_{\lambda}])/\sqrt{{\rm Var}[X_{\lambda}]}$ 
converge in distribution for $\lambda \to \infty$.
The numerical results indicate that such a central limit theorem 
might indeed hold.  
While the tools used here for estimating ${\rm E}[X_{\lambda}]$ do 
not give information about the distribution, 
it is reasonable that the randomness in the map assured
by Pesin theory makes the Lyapunov exponents in the limit
$\lambda \to \infty$ behave like the Lyapunov exponents
of Markovian cocycle, for which central limit theorems 
are known \cite{Ar+86,Bou88}. \\
}}}
\vspace{1mm}

9) {\bf The size of the Pesin region}.  \\
Because the averaged Lyapunov exponent of the Standard map can not be 
bigger than $\log(\lambda/2)+C(2,\lambda)$, the conjectured bound would
lead to the size of the Pesin region, the set where the pointwise 
Lyapunov exponent of the Standard map
is positive, has Lebesgue measure which is bounded below by 
\begin{equation}
\label{Pesinsize}
    \frac{(\log(\lambda/2)-C(\lambda))}{(\log(\lambda/2)+C(2,\lambda))}
 = \frac{(\log(\lambda/2)-{\rm arcsinh}(1/\lambda)-\log(2/\sqrt{3}))}
        {(\log(\lambda/2)+ C(2,\lambda)}  \;.
\end{equation}
For example, for $\lambda=5.42$, the Pesin region would already cover more
than half of the phase space.  (In comparison, it was measured 
numerically that for $\lambda=5.0$, the size of the elliptic regions 
is less than $2$ percent of the phase space \cite{Chi79}). 
It would follow from formula~(\ref{Pesinsize}) 
that the set, where the Lyapunov exponent is zero has measure 
which is smaller than $O(\log(4/3)/\log(\lambda))$ 
for $\lambda \to \infty$.  
If $C(\lambda)$ could be replaced by $c(\lambda)=m(2,\lambda)=
\int_{\TT^2} \log|| dT_{\lambda \cos}(x,y) || \; dx dy = O(1/\lambda)$
the complement of the Pesin region would be the order
$O(1/(\lambda \log(\lambda)))$.  \\
Empirically, the set with stable behavior has been measured 
in \cite{Chi79}. Empirical formulas have to be taken with a grain of 
salt for large $\lambda$
because the size of the elliptic islands is usually so small 
that say for $\lambda=100$ already,
a computer can hardly resolve individual islands. We know 
that from formula~(\ref{Pesinsize}) that the total measure must
then be less then 3/1000.
The size of the elliptic islands of a periodic orbit with 
period $p$ is expected to be of the order $M^{-3}$, where 
$M = \sup_{j<p} ||dT^j(x)||$ \cite{Car91}.  
There are examples in Hamiltonian dynamics where it is possible
to prove that elliptic regions cover a substantial part of the phase space 
\cite{Ne+97}. For other nonergodic Hamiltonian flows, see
\cite{Don96}. The Carleson problem is the question whether 
there are $\lambda$ for which $T_{\lambda \sin}$ is ergodic. \\

10) {\bf Vlasov-Toda deformation of Twist maps}. \\
Other cocycles with the same bounds on the
Lyapunov exponent can be obtained with isospectral Vlasov-Toda deformations 
of the random Jacobi operators $L=\Delta + V
= a \tau + (a \tau)^* + b$ \cite{Kni93a}.  
Assume we deform $L$ and so the transfer cocycle with the first Toda flow 
$\dot{a}=a (b(T)-b), \dot{b}=2a^2-2a^2(T^{-1})$, where 
\begin{eqnarray*}
a(x,y) &=& h_{12}(x,y) = \partial_x \partial_y h(x,y) \\
b(x,y) &=& h_{11}(x,y) + h_{22}(T^{-1}(x,y))
= \partial_x \partial_x h(x,y) + \partial_y \partial_y h(T^{-1}(x,y)) \;. 
\end{eqnarray*}
It is not clear whether the deformed operator $L_t = a_t \tau + (a_t \tau)^* + b_t$ 
is the Hessian of a critical point defining a twist map. 
The question is whether a generating function $h$ can be derived from 
the deformed $a,b$ such that $h$ generates a twist map $T$ satisfying
discrete time Euler-Lagrange equations $h_2(T^{-1}) +  h_1 = 0$. 
Because the deformation keeps the Lyapunov exponents of the cocycle 
constant, the lower bound on the entropy would hold also for the
deformed twist maps. In any case, even if the 'isospectral deformation
of the twist maps' should not exist in general, Toda deformation gives
infinite dimensional manifolds of cocycles over a fixed $T \in \Ycal$
obtained from the cocycle 
$A_{T,\lambda \sin}$ of the Standard map for which one has the same 
Lyapunov exponent. \\

11) {\bf Continuity of the entropy}.  \\
The map $\lambda \mapsto \mu(T_{\lambda \sin})$ is upper-semicontinuous. 
We don't know whether it is continuous. 
On the other hand, the topological entropy 
depends continuously on $\lambda$ because it is always 
a continuous function on $C^{\infty}$ diffeomorphisms on two-dimensional 
manifolds \cite{New89}. \\

\vspace{1.0cm}
\hbox{
\vbox{ \parbox{5.2cm}{ \vspace{-4.0cm}
 \setlength{\unitlength}{0.07mm} 
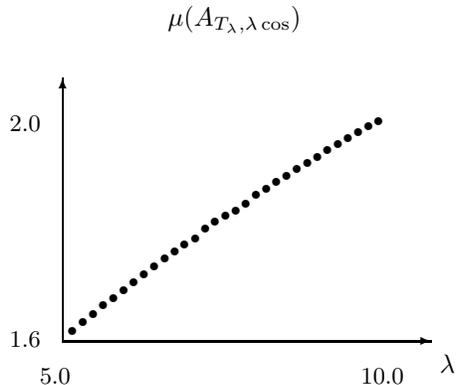
 \begin{picture}(800,1000)      
 \put(0,0){\vector(1,0){700}} 
 \put(0,0){\vector(0,1){500}} 
 \put(700,-60){ $\lambda$ } 
 \put(550,-80){ {\small $10.0$} } 
 \put(-60,-80){ {\small $5.0$} } 
 \put(-100,400){{\small $2.0$}} 
 \put(-100,-10){{\small $1.6$}} 
 \put(200,600){$\mu(A_{T_{\lambda},\lambda \cos})$} 
\cc{19}{20}\cc{39}{36}\cc{58}{52}\cc{77}{68}\cc{97}{83}
\cc{116}{98}\cc{135}{113}\cc{155}{128}\cc{174}{143}\cc{194}{157}
\cc{213}{171}\cc{232}{184}\cc{252}{196}\cc{271}{214}\cc{290}{227}
\cc{310}{238}\cc{329}{249}\cc{348}{261}\cc{368}{279}\cc{387}{290}
\cc{406}{303}\cc{426}{315}\cc{445}{327}\cc{465}{339}\cc{484}{351}
\cc{503}{363}\cc{523}{374}\cc{542}{386}\cc{561}{397}\cc{581}{408}
\cc{600}{419} 
\end{picture}         
}}
\hspace{-9cm} 
\vbox{ \parbox{9.0cm}{
The topological entropy of surface diffeomorphisms
changes only through homoclinic bifurcations \cite{Xia97}. 
Motivated from results on unimodal maps \cite{Dou95}
and consistent with experiments, it is reasonable to ask whether
the topological entropy or the metric 
entropy with respect to the invariant Lebesgue measure 
depends monotonically on $\lambda$. The Figure to the left shows
the numerically computed Kolmogorov-Sinai entropy for $30$ values of 
$\lambda$ in the interval $[5.0,10.0]$ (Lyapunov exponents averaged
over $40 \times 40$ orbits each of length $10^6$.) 
}}}
\vspace{0.5cm}
The dependence of 
the entropy on $\lambda$ even appears to be smooth. Oscillations
indeed become small for larger $\lambda$ because the established 
lower bound is realanalytic in $\lambda$ and the effective value
of the entropy is in a small corridor 
$\log(\lambda/2) + [-c(\lambda),c(\lambda)]$ which becomes 
increasingly narrow for $\lambda \to \infty$.  \\

12) {\bf Nonuniform hyperbolicity and homoclinic tangencies}. \\
Assume $\lambda$ is such that
$T_{\lambda}$ has positive Lyapunov exponents on a set of positive Lebesgue 
measure. The Pesin region of positive Lebesgue measure 
is contained in the closure of 
transverse homoclinic points or in the closure of hyperbolic periodic 
orbits with transverse homoclinic intersections (see \cite{KH}). 
We don't expect however in the case of the Standard map that for some 
$\lambda$, we have uniform hyperbolicity on a set $Y$ of positive Lebesgue 
measure. This is supported by the fact that for large $\lambda$, there 
exists a dense set of parameter values $\lambda$ for which homoclinic
bifurcations \cite{PaTa} happen in the Standard map family 
\cite{Dua94}. Uniform hyperbolicity on $\TT^2$ 
is excluded because the Chirikov Standard map can not be Anosov on the
whole torus $\TT^2$: it would be topologically conjugated to a hyperbolic 
automorphism on the torus \cite{Man74} and by homotopy to 
$(x,y) \mapsto (2x-y,x)$ which is not hyperbolic.
Equivalent to uniform hyperbolicity on a set $Y \subset \TT^2$ is the 
property that there exists an interval $I$ containing $E=0$ such that 
$I$ is disjoint from in the spectrum of $L(x,y)$ for almost all 
$(x,y) \in Y \subset \TT^2$. This is equivalent to the fact that
the random operator $L$ associated to an invariant ergodic set of 
positive Lebesgue measure is invertible in the corresponding crossed 
product algebra. It is an open question whether hyperbolic sets always 
have measure zero or one in the Chirikov Standard map case. This property is 
Baire generic in the class of measure preserving $C^1$-diffeomorphisms
\cite{Man96} and probably always holds for realanalytic diffeomorphisms on 
the torus. Having excluded uniform hyperbolicity almost everywhere,
the positive entropy result could provide a new proof that homoclinic 
tangencies and consequently elliptic islands occur in the Chirikov Standard 
map for a dense set of parameters in $[\lambda_0,\infty)$. \\

13) {\bf Coexistence}.  \\
Conjecture~\ref{II} establishes the existence of explicit real analytic maps
on the torus, where coexistence of elliptic islands and positive metric
entropy holds. Previously known were piecewise smooth maps 
\cite{Woj81,Woj82} or 
$C^{\infty}$ maps \cite{Prz82}. Whether true coexistence 
(in the sense that the Pesin region as well as its complement are dense
\cite{Str89}) 
can hold on open sets of the phase space 
for some Standard maps is not known. According to \cite{Str89},
this coexistence problem was posed already before 1969 by Sinai. 
It was measured numerically in \cite{UmFa85} that the closure of some 
orbits might have fractional box counting dimension. In higher dimensions, 
Mather conjectured \cite{Mat86} that generically one should have 
transitivity in the sense that there exists orbits which are dense in 
the phase space. No mathematical results about these questions are yet 
available. 
It has also been conjectured that transitive components 
in the complement of quasiperiodic sets are trellises \cite{Eas86}, 
closed sets obtained by closing the unstable manifold of some hyperbolic 
periodic point. For the Standard map, where we know now the existence
of ergodic sets of positive measure, one can ask about the nature of 
these ergodic sets (i.e. the Hausdorff dimension of the closure). 
Careful measurements done in \cite{Mei94}
indicate that some ergodic components are obtained as the closure 
of some aperiodic orbits.   \\

14) {\bf The complement of the Pesin region}. \\
An interesting question is what is left on the complement of
the union of the almost periodic KAM set and the Pesin region. 
Is the dynamics weakly mixing but not mixing on a set of positive 
measure? While Pesin theory together with a genericity result for shift 
invariant measures \cite{Kni98a} provide many invariant measures which 
have this property, they are in general supported on sets of zero
Lebesgue measure. This observation by \cite{Aub92} can be strengthened
by estimating the Hausdorff dimension of the measure. 
It was first observed in \cite{Che91} that there are
invariant measures $\mu_{\lambda}$ whose Hausdorff dimension goes to $2$
for $\lambda \to \infty$. This observation is based
on Young's formula \cite{You82} for the Hausdorff dimension $H(\mu)$ of 
$\mu$:
$$   H(\mu) = 2 h_{\mu}(T)/\lambda(\mu)  \; , $$
where $\lambda(\mu)$ is the Lyapunov exponent integrated
with respect to the measure $\mu$ and where $h_{\mu}(T)$ is
the metric entropy of $T$ with respect to the measure $\mu$. \\
While periodic orbits are dense on the Pesin region \cite{Kat80} and on 
the closure of KAM orbits (in the sense that for $x$ in the KAM set, there 
exists periodic $y_n \to x$), the question
whether periodic orbits are dense in $\TT^2$ is open and related to the 
problem to analyze the dynamics on the complement of the union of 
the Pesin and KAM regions. \\

15) {\bf A variational problem}.  \\

\vspace{1.0cm}
\hbox{
\vbox{ \parbox{5.2cm}{ \vspace{-4.0cm}
 \setlength{\unitlength}{0.07mm} 
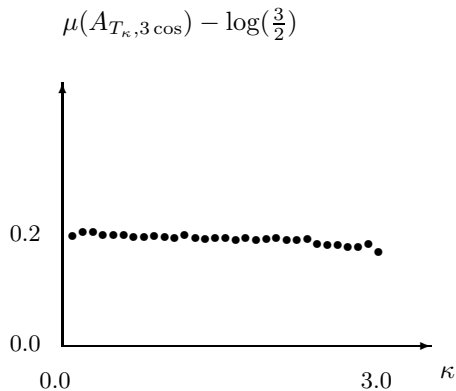
 \begin{picture}(800,1000)      
 \put(0,0){\vector(1,0){700}} 
 \put(0,0){\vector(0,1){500}} 
 \put(700,-60){ $\kappa$ } 
 \put(550,-80){ {\small $3.0$} } 
 \put(-60,-80){ {\small $0.0$} } 
 \put(-100,200){{\small $0.2$}} 
 \put(-100,-10){{\small $0.0$}} 
 \put(0,600){$\mu(A_{T_{\kappa},3 \cos})-\log(\frac{3}{2})$} 
\cc{19}{209}\cc{39}{217}\cc{58}{218}\cc{77}{211}\cc{97}{211}
\cc{116}{212}\cc{135}{208}\cc{155}{207}\cc{174}{209}\cc{194}{208}
\cc{213}{206}\cc{232}{212}\cc{252}{206}\cc{271}{204}\cc{290}{206}
\cc{310}{206}\cc{329}{202}\cc{348}{205}\cc{368}{201}\cc{387}{203}
\cc{406}{205}\cc{426}{201}\cc{445}{202}\cc{465}{203}\cc{484}{195}
\cc{503}{193}\cc{523}{193}\cc{542}{189}\cc{561}{189}\cc{581}{194}
\cc{600}{180} 
 \end{picture}         
}}
\hspace{-9cm} 
\vbox{ \parbox{9.0cm}{
For fixed $\lambda$, is there (a realanalytic) $T \in \Xcal$ (rsp. $\Ycal$) 
for which the Lyapunov exponent of the fixed matrix valued
map $A_{\lambda \cos}$ is maximal? If $T$ is the identity, we 
have $\mu(A_{\lambda \cos,T}) \geq \log|\lambda/2|$. There are many maps
$R_{\alpha} T R_{-\alpha}$ with $R_{\alpha}(x,y) = (x+\alpha,y)$ for
which the Lyapunov exponent is bigger than $\log(\lambda/2)$. 
If we fix the cocycle $A_{\lambda}$, how 
does the Lyapunov exponent depend on $T_{\kappa}$, when we 
deform $T_{\kappa}$ from the Mathieu case $T_0$ to the Standard map case 
$T_{\lambda}$? The graph to the left shows a numerical 
computation in the case $\lambda=3.0$. 
The averaged Lyapunov exponent seems to decrease slightly 
when $\kappa$ is increasing.   \\
}}}
\vspace{0.5cm}
The fact that the Lyapunov exponent is close to the one in the 
Mathieu case $\kappa=0$ is compatible with the empirical fact
of being close to Aubry duality (see 3) 25)). \\
It would be good to know the Taylor expansion of 
$$ \kappa 
 \mapsto 
 \int_{\TT^2} \log( [A^n_{E,\lambda,T_{\kappa}}(w)(x,y)]_{11} ) \; dx dy $$
at $\kappa=0$, for $|w|<1$ which exists and converges for large $\lambda$. 
This requires to understand the dynamics for complex $\kappa$, where 
$(x_n,y_n)$ become complex too. \\

16) {\bf Thermodynamic limit}. \\
For classes of symplectic coupled map lattices like finite 
dimensional versions of \cite{KaGr87}, we could obtain lower bounds 
for the entropy which are independent of the dimension of the map. 
In the thermodynamic limit, these maps define
homeomorphisms $T$ on the compact metric space
$(\TT^{2})^{\ZZ}$ which preserve the product measure. 
While one has no more an Oseledec theorem in this infinite 
dimensional situation (the results in \cite{Rue82,Man83} do not apply), 
the averaged maximal Lyapunov exponent
$\lim_{n \to \infty} n^{-1} \int_{(\TT^2)^{\ZZ}}
      \log ||dT^n({\bf x}|| \; d {\bf x}$ is defined and positive for 
large $\lambda$ and small $\epsilon$. The cocycle $A=dT$ is now
a bounded linear operator on $l^2(\ZZ,\CC^2)$. 
By Vesantini's theorem \cite{Ransford}, the map $w \mapsto 
n^{-1} \log(\rho(A^n(w))$ is subharmonic if 
$\rho(A)$ denotes the spectral radius of $A$. 
The estimates can be done directly for the operator-valued cocycle. 
It is not clear however whether the
metric entropy of the homeomorphism $T$ with invariant product measure
is the limit of the metric entropies in the finite-dimensional 
situations. \\

17) {\bf Response formulas}. \\
The subharmonic estimates could apply for classes of dissipative 
Standard maps on the cylinder. 
We expect the averaged Lyapunov exponent to depend analytically on the 
map in open sets of realanalytic maps because our
results imply that smooth observables see a pretty uniform hyperbolic 
behavior of the dynamics. Indeed, by taking 
periodic Rohlin approximations to the map $T$, we get
$$ n^{-1} \int_{\TT^2} \log ||dT^n_{\lambda \sin}|| \; dx dy 
\geq \log(\lambda/2) - C_n(\lambda)  \; ,           $$ 
where $C_n(\lambda) = C(\lambda) + (\log(\lambda/2)+C(\lambda))/n$. 
The hyperbolicity appears to be uniform when observing the 
dynamics with smooth observables.\\ 
Could it be possible even that Ruelle's
response formula \cite{Rue97} holds for some Standard maps with 
positive Lyapunov exponents? 
This formula says that for an infinitesimal perturbation $T+\delta T$ of 
an Anosov map $T$, the natural $T$-invariant measure $\rho_T$ will 
change to $\rho_{T+\delta T}=\rho_T + \delta \rho_T$ with
$\delta \rho_T(\Phi) = - \sum_{n=0}^{\infty}
\rho_T((\Phi \circ T^n) \cdot {\rm div}(X))$, 
where $X = \delta T \circ T^{-1}$ is the vector field
associated to the change of the map and where $\Phi \in C(\TT^d)$ is
an observable. The response formula would be expected
to hold only for real analytic perturbations of $T$ in $\Xcal$ and for 
a subset of smooth observables $\Phi$ which do nowhere vanish 
on the torus. One has to see this question in the more general context of 
the 'chaotic hypothesis' \cite{Gal98}. \\

18) {\bf Symmetry}. \\
The symmetry $T(-x,-y)=-T(x,y), DT(x,y)=DT(-x,-y)$ which holds in the 
case of the Chirikov Standard map seems to make
the Lyapunov exponent extremal with respect to some parameter perturbations. 
This is supported by some numerical experiments. It is not known 
whether the Lyapunov exponent is continuous as a function of $w$ for example
when $w \mapsto L_w$ is the random operator~(\ref{L(w)}),
so that symmetry does not imply a local extremum. \\

\vspace{7mm}
\hbox{
\vbox{ \parbox{5.2cm}{ \vspace{-2.8cm}        
 \setlength{\unitlength}{0.05mm} 
 \begin{picture}(1000,1000)      
 \put(600,0){\vector(1,0){700}} 
 \put(-200,0){\vector(1,0){700}} 
 \put(-200,0){\vector(0,1){500}} 
 \put(600,0){\vector(0,1){500}} 
 \put(450,-60){ $\epsilon$ } 
 \put(65,-80){ {\small $0$} } 
 \put(370,-80){ {\small $\pi$} } 
 \put(-260,-80){ {\small $-\pi$} } 
 \put(1250,-60){ $\alpha$ } 
 \put(860,-80){ {\small $0$} } 
 \put(1150,-80){ {\small $\pi$} } 
 \put(575,-80){ {\small $-\pi$} } 
 \put(-300,220){{\small $0.2$}} 
 \put(-300,-10){{\small $0.0$}} 
 \put(-300,600){$\mu(A_{T,3 \cos}(x+\epsilon))-\log(\frac{3}{2})$} 
 \put(600,600){$\mu(A_{T_{\alpha},3 \cos})-\log(\frac{3}{2})$} 
\cc{610}{180}\cc{620}{182}\cc{630}{182}\cc{639}{179}\cc{649}{183}
\cc{659}{183}\cc{669}{181}\cc{679}{173}\cc{689}{171}\cc{698}{162}
\cc{708}{158}\cc{718}{153}\cc{728}{154}\cc{738}{152}\cc{748}{145}
\cc{757}{149}\cc{767}{159}\cc{777}{165}\cc{787}{177}\cc{797}{178}
\cc{807}{179}\cc{816}{174}\cc{826}{179}\cc{836}{175}\cc{846}{173}
\cc{856}{163}\cc{866}{158}\cc{875}{160}\cc{885}{177}\cc{895}{180}
\cc{905}{188}\cc{915}{180}\cc{925}{177}\cc{934}{160}\cc{944}{158}
\cc{954}{163}\cc{964}{173}\cc{974}{175}\cc{984}{179}\cc{993}{174}
\cc{1003}{179}\cc{1013}{178}\cc{1023}{177}\cc{1033}{165}\cc{1043}{159}
\cc{1052}{149}\cc{1062}{145}\cc{1072}{152}\cc{1082}{154}\cc{1092}{153}
\cc{1102}{158}\cc{1111}{162}\cc{1121}{171}\cc{1131}{173}\cc{1141}{181}
\cc{1151}{183}\cc{1161}{183}\cc{1170}{179}\cc{1180}{182}\cc{1190}{182}
\cc{1200}{180}\cc{-190}{174}\cc{-180}{176}\cc{-170}{179}\cc{-161}{183}\cc{-151}{186}
\cc{-141}{186}\cc{-131}{182}\cc{-121}{179}\cc{-111}{179}\cc{-102}{183}
\cc{-92}{179}\cc{-82}{175}\cc{-72}{167}\cc{-62}{169}\cc{-52}{164}
\cc{-43}{161}\cc{-33}{159}\cc{-23}{147}\cc{-13}{130}\cc{-3}{108}
\cc{7}{78}\cc{16}{93}\cc{26}{106}\cc{36}{122}\cc{46}{139}
\cc{56}{151}\cc{66}{163}\cc{75}{173}\cc{85}{178}\cc{95}{181}
\cc{105}{180}\cc{115}{181}\cc{125}{178}\cc{134}{173}\cc{144}{163}
\cc{154}{151}\cc{164}{139}\cc{174}{122}\cc{184}{106}\cc{193}{93}
\cc{203}{78}\cc{213}{108}\cc{223}{130}\cc{233}{147}\cc{243}{159}
\cc{252}{161}\cc{262}{164}\cc{272}{169}\cc{282}{167}\cc{292}{175}
\cc{302}{179}\cc{311}{183}\cc{321}{179}\cc{331}{179}\cc{341}{182}
\cc{351}{186}\cc{361}{186}\cc{370}{183}\cc{380}{179}\cc{390}{176}
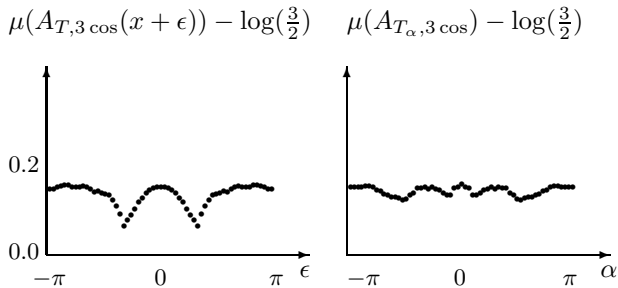
\cc{400}{174} \put(-150,-270){\small{Computing the Lyapunov exponent in two cases }} 
 \put(-150,-340){\small{with lack of symmetry. To the right, we use }} 
 \put(-150,-410){\small{$T_{\alpha}(x,y)=(2x+3 \sin(x)-y+\alpha,x) \in {\cal{X}}$. The }} 
 \put(-150,-480){estimate known for $\epsilon=\alpha=0$ extend.} 
 \end{picture}         
}}
\hspace{-8cm}
\vbox{ \parbox{7.5cm}{
In the Figure to the left, we see two experiments. 
The left graph shows the Lyapunov exponent of the deformed
cocycle $\phi \mapsto A(x+\phi)$, where $A(x)$ is the 
Jacobean cocycle of the Standard map. This is equivalent to a  
deformation of the map $T$ by conjugation
$T_{\phi} = R_{\phi}^{-1} \circ T R_{\phi}$. 
In the right graph, the map $T$
is moved along a path $\alpha \mapsto T_{\alpha} = R_{\alpha} \circ T$
in  $\Xcal$, where $R_{\alpha}(x,y) = (x+\alpha,y)$. For $\phi=\alpha=0$, we 
have the Chirikov Standard map case with $\lambda=3$. In the graph to the 
left, Herman subharmonicity argument shows that 
$\int_{\TT} \mu(A_{T,3 \cos}(\cdot+\alpha)) \; d\alpha -\log(3/2)>0$. 
Nevertheless, in both experiments, we never see a Lyapunov exponent below
$\log(\lambda/2)$. 
}}}
\vspace{1.5cm}

19) {\bf Diffusion and Sinai's (H1) conjecture}.  \\
A discrete Legendre transform brings the Standard map into the Hamiltonian 
form $T:(x,y) \mapsto (x+y + f(x), y+ f(x) )$ which is a map on 
the cylinder, the cotangent bundle $T^* \TT = \TT \times \RR$ of the circle. 
Let $A \subset \TT^2$ be a measurable $T$-invariant set of positive 
Lebesgue measure. 
If $(x_j,y_j)=T^j(x,y)$ is an orbit in the universal cover $\RR^{2d}$,
then $X_j=y_j-y_{j-1}=f(x_j)=\psi(x_j,y_j)$ are random variables on the 
probability space $\Omega$ equipped with the normalized Lebesgue measure.
They have the mean ${\rm E}[X_j]=\int_{\TT^{2}} f(x_j) \; dx dy =0$ and the 
variance ${\rm Var}[X_j]=\int_{\TT^{2}} \psi(x,y)^2 \; dx dy$. Interesting is the 
growth rate of $S_n^2=(\sum_{j=0}^{n-1} X_j)^2 = (y_n-y_0)^2$. 
Using translational invariance 
${\rm E}[X_j X_l] = {\rm E}[X_{j+1},X_{l+1}]$, the variance of $S_n$ is 
$$ {\rm Var}[S_n]= {\rm E}[ (\sum_{j=0}^{n-1} X_j)^2] 
   = n {\rm Var}[X] + \sum_{j=1}^{n-1} (n-j) {\rm E}[X_0 X_j]  
   = n {\rm Var}[X] + \sum_{j=1}^{n-1} (n-j) \hat{\mu}_{\psi}(j)  \; , $$
where $\mu_{\psi}$ is the
spectral measure of $\psi \in L^2(\Omega), \psi(x,y) =f(x)$ 
with respect to the unitary Koopman operator 
$g \mapsto g(\tilde{T})$, where $\tilde{T}$
is the map $T$ induced on $\Omega$. The Fourier
transform of $\mu_{\psi}$ is $\hat{\mu}_{\psi}(j)
= {\rm E}[X_0 X_j]={\rm Cov}[X_0,X_j]$,
a correlation function.
Let $\beta$ be the infimum over all real numbers for which
$\limsup_{n \to \infty}
  n^{-\beta} \sum_{j=1}^{n-1} (n-j) \hat{\mu}_{\psi}(j)$ is finite.
If $\beta=1$, then $S_n$ behaves like a random walk (a case, where the 
$X_j$ are independent) and
$D = \int_{\TT^{2}} V'(x)^2 \; dx +
   \limsup_{n \to \infty} n^{-1} \sum_{j=1}^n (n-j) \hat{\mu}_{\psi}(j)$
is the diffusion constant. \\
It is a conjecture of Sinai (H1) on page 144
in \cite{Sinai94} that there exists 
a set $\Omega$ of positive Lebesgue measure
for which $\beta=1$ if $\lambda$ is larger then $\lambda_{crit}$, where
the last homotopically nontrivial KAM torus disappears. 
If $\hat{\mu}_{\psi}(n)$ would decay fast enough, then
$D = \int_{\TT^{2}} V'(x)^2 \; dx 
+ \sum_{j=0}^{\infty} \hat{\mu}_{\psi}(j)$.
First numerical experiments were done by Chirikov and Hizanidis
\cite{ReWh80}. Numerically, the Standard is reported to show 
such Brownian diffusive behavior for 
large $\lambda$ \cite{Lichtenberg,Leb98}. The fact of having positive 
Lyapunov exponents on a set of positive measure makes it plausible 
that the random variables $X_0,X_n$ get decorrelated for 
$n \to \infty$. 
When computing numerically the first few hundred Fourier 
coefficients for smaller $g$ (like $g=5$) and checking with the 
Wiener theorem, we got the impression that $\mu_{f}$ still 
has some atoms (which prevents decorrelation) if $\Omega = \TT^2$. 
Indeed, the presence of elliptic islands could be responsible for an 
almost periodic component in the Fourier transform of $\mu_{f}$. 
In numerical experiments, one can get rid of this discrete part of 
the spectrum by adding stochastic noise \cite{Lichtenberg}.  \\
In any case, the (H1) conjecture of Sinai would be settled for
$\lambda>\lambda_0$ if one
could show a fast enough decay of correlation of the spectral measure
of $\psi(x,y)=f(x)$ on a mixing component of the Pesin set $\Omega$. 
No exponential decay is necessary. It is enough to establish a power law
decay of correlation $\hat{\mu}_{\psi}(j) = O(j^{-2})$. 
This finite differentiability condition for the spectral measure $\mu_{\psi}$ 
is reasonable since the dynamics on $\Omega$ is conjugated to a Markov chain. 
We even expect many spectral measures to be realanalytic leading to 
exponential decay of correlations. \\

20) {\bf Dissipative Standard maps}. \\
The Lyapunov exponent of dissipative Standard maps
numerically often satisfy the lower bound of the conservative case. 
The Lyapunov exponent can drop however to zero. 
In order that the proof carries over to the dissipative case, it appears
however, that the random variables $x_j$ need a smooth distribution $\mu$
sufficiently close to the uniform distribution. 
Results for Henon maps lead to the expectation that 
for many parameters, there exist invariant SRB measures for 
dissipative Standard maps like 
$$ T_b: (x,y) \mapsto (x+y + \lambda \sin(x), b(y + \lambda \sin(x))) \; ,  $$
with $b<1$. Results like in \cite{Cal91,Mor94} support this. 
Viana conjectured that in general a map with nonzero Lyapunov exponents
almost everywhere in the phase space has an SRB measure \cite{Via98}.
If this conjecture is true, it indicates that for most values of $b$, an
SRB measure should exist. 

The estimates of Lyapunov exponents using spectral methods 
could even carry over to the case $b=0$, where on 
each invariant set $y=\alpha$ we get the 
one-dimensional Arnold family $x \mapsto x+\alpha + \lambda \sin(x)$. 
If such a circle map has a smooth invariant measure sufficiently close
to the Lebesgue measure, the Lyapunov exponent with respect to 
this measure should be $\geq \log(\lambda/2) - C(\lambda)$
for large $\lambda$.  \\

21) {\bf More general stability of positive metric entropy?}  \\
The stability, we have established holds 
for estimates of the entropy which are done using subharmonicity 
rsp. the Jensen
formula. This result provokes the question, whether, in general,
positive metric entropy is an open property in the realanalytic
category: is it true that for any realanalytic, measure preserving 
diffeomorphism $T$ on the torus with positive metric entropy,
there exists a Banach space of 
realanalytic measure-preserving maps such that an open neighborhood 
of $T$ has positive metric entropy? \\
A question related to this stability problem is 
whether the Riesz measure $dk$
of the $w$-parameterization of the operator $L$ for a general twist 
map has the property that the potential 
$\int \log|w-w'| \; dk(w')$ does not fluctuate too much around its
mean on the unit circle $\{ |w|=1 \}$. This could be used to establish
the stability if the entropy is comparable to 
upper bounds of the entropy like in the Standard map case. \\

22) {\bf The Herman spectrum}. \\ 
One can also look at different analytic parameterizations of the cocycle.
The Herman spectrum of a cocycle $A$ 
is the set of complex numbers $z$, such that $R(z) A$ is not uniformly 
hyperbolic, where
$R(z)=\left( \begin{array}{cc} 
                  \frac{z+z^{-1}}{2}  & \frac{z-z^{-1}}{2i} \\
                -\frac{z-z^{-1}}{2i}  & \frac{z+z^{-1}}{2}  \\
             \end{array} \right) \in SL(2,\CC)$ which has the 
property that $R(e^{i \alpha}) = R(\alpha) \in SO(2,\RR)$. 
It is a subset of $\{|z|=1\}$ because for $|z| \neq 1$,
one can find a strict coinvariant cone field in $\CC P^1$. 
There is a measure $\mu$ supported by the Herman spectrum and 
a harmonic function $g$ such that the Lyapunov exponent 
$\mu(z) = \mu(z R(z) A)$ satisfies
$$ \mu(z) = \int \log|z-z'| \; d\mu(z) + g(z)  \; . $$
This abstract Thouless formula \cite{Her83} follows from Riesz theorem and 
$\mu$ is the Riesz measure of the subharmonic function $\mu(z)$. 
The proof of Conjecture~\ref{Igeneral} shows
$\log(\mu(\beta)) > \log( \cos(\beta) \lambda/2)-O(1/\lambda)$ for all $\beta$. 
The Lyapunov exponent is realanalytic and positive outside the 
Herman spectrum. At such points, the directional derivative with 
respect to variations of the angle $\beta$ in $z = r e^{i\beta}$ 
can be computed using a formula of Ruelle \cite{Rue79}
for the Fr\'echet derivative of Lyapunov exponents on the open 
set of uniformly hyperbolic cocycles. One gets
$\frac{d}{d\beta} \mu(\beta)= \int_{\TT^2} \cot(\omega(x,y)) \; dm(x)$,
where $\omega(x)$ is the angle between the stable and unstable directions
$m^+,m^-$ at the point $(x,y)$ (see \cite{Kni93diss}). 
This formula shows that the Lyapunov exponent can change a lot if the
stable and unstable manifolds are close. In the case
$A(x)=\left( \begin{array}{cc} c & b(x)    \\
                               0 & c^{-1}  \\ \end{array} \right)$
with constant $|c| \neq 1$, one obtains
$$ \frac{d}{d\beta} \mu(A(\beta))=\int_{\TT^2} \cot(\omega(x,y)) \; dx dy
                             = \frac{c}{1-c^2} \int b \; dx dy  \; .$$
If $A$ is the cocycle of a Standard map in Hamiltonian form like 
in 19), there is a spectral gap 
containing $\beta=-\pi/4$ in the Herman spectrum, because
$ R(-\pi/4) dT_{\lambda}=
    \left( \begin{array} {cc}
                    2^{-1/2}                              &     0      \\
                    2^{-1/2}+2^{1/2}\lambda \cdot \cos(x) & 2^{1/2}    \\
           \end{array} \right)   $. We compute
$\frac{d}{d\beta} A(-\pi/4) = 1$. 
The point $z=e^{i 0}$ is in the Herman spectrum because it follows
from \cite{Mat68,Rue85} (see also \cite{Kni93diss})
that $z=1$ is outside the Herman spectrum if and only if the map
$T_{\lambda f}$ is Anosov. 
The rotation number of Ruelle \cite{Rue85} defined for lifts of 
$A$ into the universal cover of ${\rm SL}(2,\RR)$ plays the 
role of the integrated density of states in the case of the 
Schroedinger spectrum. It is uniquely defined up to a multiple
of $2 \pi$ if one fixes $\rho(A(0))=0$.  
The rotation number is constant on an interval $I$ if and only 
if $I$ is not in the Herman spectrum \cite{Kni93diss}. Furthermore,
$\rho(\beta) = \int_0^{\beta} \; d\mu(\alpha)$, 
showing that $\mu$ is an 'integrated density of states'. From the stable
and unstable direction fields, one can construct (nonselfadjoint)  
random Jacoby operators $L(z)$ having those direction fields 
$m^{\pm}(z)$ as Titchmarsh-Weyl 
functions (with energy $E=0$). 
They are analytic in $z$ outside the Herman spectrum. 
Call a subset of the unit circle a part of the 'absolutely 
continuous spectrum' of $L(z)$ if there, $m^+(z)= \overline{m^-}(z)$. 
This absolutely continuous spectrum is a subset of the Herman spectrum. 
An adaption of Kotani theory \cite{Cycon} shows that
the 'absolutely continuous spectrum' is the essential closure of the set, 
where the Lyapunov exponent is zero and that in the ergodic case, 
the existence of some absolutely continuous spectrum implies that 
$T$ is 'deterministic' (see \cite{Fur63,Led86,KoSi88}). 
We know therefore for the Standard map 
that there is no absolutely continuous Herman spectrum for almost 
all $(x,y)$ in the Pesin region if $\lambda$ is large. 
We call $z$ an 'eigenvalue' of $L(z)(x,y)$ if $L(z)(x,y)$ has 
an eigenvalue $E=0$. The set of these 'eigenvalues' forms a 
'discrete spectrum' of $L(z)(x,y)$ which is a subset of the 
Herman spectrum and which we expect to be nonempty for 
Lebesgue almost all $(x,y)$ in the Pesin region. \\

24) {\bf Quasiconformality}. \\
The map 
$$  z \mapsto G_{n}(z)=\int_{\TT^{1}} \log|z^n A^n_{E}(z,y)| \; dy  $$
is not conformal for small $|z|$ in general. 
For realanalytic $T$ near a linear automorphism
$(x,y) \mapsto (x+y,y)$ it is quasiconformal for all $n$. 
Question: is there an open set of realanalytic, measure-preserving maps 
$T$ for which $G_n$ are quasiconformal? If yes, it would be important 
to estimate the complex dilation 
$\kappa(z)=\overline{\partial} G_{n}(z)/\partial G_n(z)$ and to 
estimate the Lyapunov exponent in terms of $\kappa$. For
$\kappa=0$, one has conformality and the subharmonic estimates 
of Herman apply. \\

25) {\bf The Aubry duality transform}. \\ 
The Aubry duality transform is an involutive map on a class of 
operators. It provides in 
the Mathieu case an elegant way to estimate the Lyapunov exponent. 
The transform can be defined in more general situations: 
Let $T_k$ be a cyclic interval exchange transformation on $\TT$
and let $L_{\lambda \cos,T_k}$ be the corresponding
Schr\"odinger operator. The density of states of the random operator
$L$ is the same as the density of states of the operator
$$ (L u)_n(\theta) = u_{n+1}(\theta) + u_{n-1}(\theta)
        + \lambda  \cos( T^n \theta) u_n(\theta) $$
on the Hilbert space $H = L^2(\TT \times \ZZ)$.
There is a piecewise smooth potential $V$ such that
$\cos(T^n(\theta)) = V(\theta + n \alpha)$, where $\alpha=1/k$.
Let $V(\theta) = \sum_n V_n \exp(i n \theta)$ be its Fourier series.
The duality transform (see \cite{Gor+97})
$$ (U u)_m(\eta) = \sum_{n \in \ZZ} \int_{\TT}  e^{-(\eta+2\pi m \alpha)n}
              e^{-i m \theta} u_n(\theta) \; d\theta    $$
satisfies
$ \tau U = (\lambda/2) U \sigma$, where $\tau u_n(\theta) = u_{n+1}(\theta)$
and $\sigma u_n(\theta) = \exp( i (\theta+n\alpha)) u_n(\theta)$.
The operator $L$ can be written as 
$L = \tau+\tau^* + \lambda \sum_n V_n \sigma^n$
so that $U^* L U = \sigma+\sigma^* + \lambda/2 \sum_n V_n (\tau^n+(\tau^n)^*)$.
We are interested in the density of states 
$dk_{\epsilon}$ of
$L_{\epsilon} = (2/\lambda) U^* L U = \epsilon (\sigma + \sigma^*) +
 \sum_{n} v_n (\tau^n + \tau_n^*)$ and 
$f(\epsilon)=
  \log(\det(L_{\epsilon}))=\int \log|E| \; dk_{\epsilon}(E)$ which satisfies
$f(0)=0$. Because 
$\log(\det(L_{\lambda \cos,T_k})) \geq \log(\lambda/2) + f(2/\lambda)$, we 
want to estimate $f(\epsilon)$ from below for small $\epsilon$. 
If $T_k(x)=x+1/k$, where $f(\epsilon)$ is the Lyapunov exponent of a
symplectic transfer cocycle, we have 
$\int \log|E| \; dk_{\epsilon}(E) \geq 0$. In general, the Thouless formula
just becomes $\log(\det(L_{\epsilon})) = \lim_{N \to \infty} 
\sum_{k=1}^{2N} \int_{\TT} \log(\lambda_j(x,\epsilon)) \; dx$,
where $\lambda_j(\epsilon,x)$ are the eigenvalues of the truncated 
$N \times N$ matrix $L_{\epsilon,N}(x)$. While $f(0)=0$ and 
a perturbation lemma of Lidskii (\cite{Simon79}) assures that 
$|\lambda_j(\epsilon,x)-\lambda_j(\epsilon,0)| \leq \epsilon$, this does
not allow us to estimate $f(\epsilon)$ from below. The perturbation problem
to estimate $f(\epsilon)$ seems still difficult and it is not clear, whether
the Aubry transform, which transformed the perturbation problem
$V \mapsto V + \epsilon \Delta$ into a perturbation 
$\Delta \mapsto \Delta + \epsilon V$ has made things simpler. \\

26) {\bf Lower bound for topological entropy?} \\
The topological entropy of $T \in \Xcal$ on $\TT^2$ is bounded below by
$\log({\rm sp}(T_*))$, where ${\rm sp}(T_*)$ is the spectral radius of
$T_*: H_*(\TT^2) \mapsto H_*(\TT^2)$ \cite{MiPr77}.
Can the metric entropy of a (smooth) $T \in \Xcal$ 
with respect to the invariant Lebesgue measure become smaller than 
${\rm sp}(T_*)$? \\ 

27) {\bf Conjectures on the Standard map}. \\
The following of conjectures about Standard maps have still to be 
settled. \\

{\bf I)} In \cite{Chi79}, 
the entropy of the Standard map $T_{\lambda \sin}$ was measured 
$\geq \log(\lambda/2)$. Chirikov formulated the possibility that the elliptic
islands might cover arbitrary large regions of the phase space. \\

\begin{center}
\fbox{ \parbox{12cm}{
{\bf I \cite{Kni93diss}:} The Kolmogorov-Sinai entropy of the 
Chirikov Standard map 
is $\geq \log(\lambda/2)$.  \\
}}
\end{center}
\vspace{0.5cm} 

{\bf II)} \cite{Spe89} introduced the problem of 
determining the spectrum of random operators 
$(L(x,y) u)_n = u_{n+1} - 2 u_n + u_{n-1} + \lambda \cos(x_n) u_n$,
where $T^n(x,y)=(x_n(x,y),y_n(x,y))$ is the orbit starting at $(x,y)$. \\

\begin{center}
\fbox{ \parbox{12cm}{
{\bf II \cite{Spe89}:} The operator $L(x,y)$ has some point spectrum for a 
set $(x,y)$ of positive Lebesgue measure if $\lambda$ is large enough.  \\
}}
\end{center}
\vspace{0.5cm}

{\bf III)} \cite{Car91} asked whether there exists a set of 
parameters $\lambda$ with full density at $\infty$ for which 
the Chirikov Standard map has no elliptic islands. We know that we 
have for large $\lambda$ an open dense set of parameters with no 
ergodicity and positive entropy. \\

\begin{center}
\fbox{ \parbox{12cm}{
{\bf III \cite{Car91}:} There are parameters $\lambda$ with full density at 
$\infty$ for which the Chirikov Standard map $T_{\lambda \sin}$ is ergodic. \\
}} 
\end{center}
\vspace{0.5cm}

Heuristic arguments in \cite{GiLa99} predict that the Lebesgue measure 
of the set of parameters $\lambda \in [r,r+1]$
which lead to nonergodic Standard maps $T_{\lambda \sin}$ is of the 
order $O(1/r)$. \\

{\bf IV)} The following problem is related to III): \\

\begin{center} 
\fbox{ \parbox{12cm}{
{\bf IV \cite{Sin96}:}
For large $\lambda$, the Standard map has (in some Baire topology of 
maps) a 
neighborhood in which a residual set of $f$'s gives ergodic maps $T_f$.\\
}}
\end{center} 
\vspace{0.5cm}

{\bf V)} In the textbook \cite{Sinai94} p. 144, 
conjecture (H2) contained as a second part the 
statement that the entropy grows to infinity
for $\lambda \to \infty$. While we solved this part of the problem here, 
the first part of the conjecture is still open: \\

\begin{center} 
\fbox{ \parbox{12cm}{
{\bf V \cite{Sinai94}:} the entropy of the Chirikov Standard map
is positive for all $\lambda>0$. \\
}}
\end{center} 
\vspace{1.0cm} 

One should be able to prove that for any real-analytic non-constant 
periodic $f$, there exists $\lambda_0>0$ such that the Standard map
$T_{\lambda f}$ has positive Kolmogorov-Sinai entropy for all
$\lambda>\lambda_0$ and that for any real-analytic, non-constant, periodic map
$f: \TT^d \to \RR^N$ and every symmetric,
constant matrix $E \in GL(N,\ZZ)$, all the Lyapunov exponents of the
symplectic map $T_{Ex+\lambda f}$ are nonzero on a set of
positive Lebesgue measure for large enough $|\lambda|$.
All these properties are expected to be stable with respect
to realanalytic perturbations of the map $T_{E x+\lambda f}$, and
hold for a fixed cocycle for all $T \in \Ycal$. 

\bibliographystyle{plain}

\end{document}